\input amstex
\documentstyle{amsppt}

\magnification=\magstep1
\NoBlackBoxes
\NoRunningHeads
\topmatter
\title Textile systems on lambda-graph systems 
\endtitle
\author Kengo Matsumoto
\endauthor
\affil Department of Mathematical Sciences \\
 Yokohama City University\\
  Seto 22-2, Kanazawa-ku, Yokohama 236-0027, JAPAN
\endaffil
\abstract
The notions  of symbolic matrix system and 
$\lambda$-graph system for a subshift
are  generalizations of symbolic matrix and $\lambda$-graph (= finite symbolic matrix) for a sofic shift respectively ([Doc. Math. 4(1999), 285-340]).
M. Nasu introduced the notion of textile system for a pair of graph homomorphisms to study automorphisms and endomorphisms of topological Markov shifts ([Mem. Amer. Math. Soc. 546,114(1995)]).
In this paper, we formulate textile systems on $\lambda$-graph systems
and study automorphisms on subshifts.
We will prove that
 for a forward automorphism $\phi$ of a subshift
$(\Lambda,\sigma)$, the automorphisms $\phi^k \sigma^n, k\ge 0, n\ge 1$
can be explicitly realized as a subshift defined by 
certain symbolic matrix systems coming from 
both  the strong shift equivalence representing $\phi$ 
and the subshift  $(\Lambda,\sigma)$. 
As an application of this result, if an automorphism $\phi$ of a subshift 
$\Lambda$ is a simple automorphism, 
the dynamical system
$(\Lambda, \phi \circ \sigma)$ 
is topologically conjugate to the subshift  
$(\Lambda, \sigma).$ 
\endabstract

\endthanks

\endtopmatter


\def\id{{{\operatorname{id}}}}
\def\M{{ {\Cal M} }}
\def\N{{ {\Cal N} }}

\def\A{{ {\Cal A} }}

\def\P{{ {\Cal P} }}
\def\Q{{ {\Cal Q} }}

\def\K{{ {\Cal K} }}
\def\Zp{{ {\Bbb Z}_+ }}

\def\T{{  {\Cal T}_{{\K}^{\M}_{\N}} }}
\def\LCH2{{{{\frak L}^{Ch(D_2)}}}}

\heading 1. Introduction
\endheading

Let $\Sigma$ be a finite set with its discrete topology, 
that is called  an alphabet.
Let $\Sigma^{\Bbb Z}$ be the compact Hausdorff space of all bi-infinite sequences of $\Sigma$.  
One has the homeomorphism $\sigma$ defined by the left-shift 
that sends a point 
$
(\alpha_i)_{i \in \Bbb Z} \in {\Sigma}^{\Bbb Z}
$ 
into the point
$
(\alpha_{i+1})_{i \in \Bbb Z} \in {\Sigma}^{\Bbb Z}.
$ 
A subshift $(\Lambda, \sigma)$ is the topological dynamical system
 that is obtained by restricting 
the shift to a closed shift-invariant subset 
$
\Lambda
$
of $\Sigma^{\Bbb Z}$.
The space 
$
\Lambda \subset \Sigma^{\Bbb Z}
$
is uniquely determined by a set of forbidden words,
such as 
a sequence $
(x_i)_{i \in \Bbb Z} \in \Sigma^{\Bbb Z}
$ 
of $\Sigma$ 
belongs to $\Lambda$ if and only if any word in
the forbidden words can not apper 
as a  subward of
$
(x_i)_{i \in \Bbb Z}
$.
If a subshift is obtained by a finite set of its  forbidden words,
it is said to be a shift of finite type. 
It is well-known that the class of shifts of finite type coincides with the class of topological Markov shifts, 
that are defined by finite square nonnegative matrices. 
For an introduction to the theory of topological Markov shifts 
see [Ki] or [LM].
R. F. Williams [Wi] proved that two shifts of finite type are topologically conjugate if and only if their defining nonnegative matrices are 
strong shift equivalent. 
This result also says a structure of automorphisms 
of topological Markov shifts.
That is, an automorphism is given by a strong shift equivalence from the defining matrix to itself, 
and 
conversely a strong shift equivalence from the defining matrix to itself
gives rise to 
an automorphism of the shift of finite type.
M. Nasu [N] formulated strong shift equivalence between finite symbolic matrices and 
generalized the above Williams's result to sofic shifts.
He  proved that
two sofic shifts are topologically conjugate if and only 
if their canonical symbolic matrices are strong shift equivalent ([N]).

For a subshift $(\Lambda, \sigma)$, 
a homeomorphism $\varphi$ on $\Lambda$ satisfying 
$\varphi \circ \sigma = \sigma \circ \varphi $ 
is called an automorphism of $(\Lambda, \sigma)$.
It is also well-known that if an automorphism $\varphi$ of a subshift
$(\Lambda, \sigma)$ is expansive,
it is topologically conjugate to a subshift [H].
The problem studied in this paper is to {\it subshift-identify} \/
the dynamical system $(\Lambda,\varphi).$ 
Namely, for an expansive  automorphism $\varphi$ of a subshift
$(\Lambda, \sigma),$
the problem is to find, in an  explicit way,  a subshift $(\Lambda_{\varphi},\sigma)$ that is topologically conjugate to $(\Lambda, \varphi)$.
This problem has been studied in the case of topological Markov shifts and sofic shifts.
Boyle and Krieger [BK] proved that for an automorphism $\varphi$ of
topological Markov shift $(\Lambda_A,\sigma_A)$ defined by a nonnegative matrix $A$ and for all integers $n$ greater than a coding bound for $\varphi$ and $\varphi^{-1}$,
the dynamical system $(\Lambda_A, \varphi \sigma_A^n)$ is
topologically conjugate to a topologicsl Markov shift, and specified its dimension triple.  
 M. Nasu [N2] has been introduced the notion of textile system, which is very useful to analyze the automorphisms and 
endomorphisms of topological Markov shifts.
A textile system is defined by an ordered pair of graph homomorphisms
$p$ and $q$ of a directed finite graph $\Gamma$ 
into a directed finite graph $G$, that is written as
$T = (p,q:\Gamma\rightarrow G)$.
Nasu also generaized the fiomulation of the textile systems to textile systems 
on finite labeled graphs.
Among other things, he proved that  
if $\varphi$ is a forward automorphism of a sofic shift 
$(\Lambda_\A,\sigma_\A)$ 
defined by a finite symbolic matrix $\A$ 
and is given by a strong shift equivalence
$$
\A  \overset{\kappa_0}\to{\simeq}\P_1\Q_1,\quad
\Q_1\P_1  \overset{\kappa_1}\to{\simeq}\P_2\Q_2,\quad\cdots,\quad
\Q_{N-1}\P_{N-1}  \overset{\kappa_{N-1}}\to{\simeq}\P_N\Q_N, \quad
\Q_N\P_N  \overset{\kappa_N}\to{\simeq}\A,
$$
then the dynamical system
$(\Lambda_\A, \varphi^k\sigma_\A^n)$ 
is topologically conjugate to the sofic shift define by the symbolic matrix
$\P^k\A^n$ for all $k \ge 0$ and $n \ge 1,$
where
$\P = \P_1\cdots \P_N.$
If in particular $\varphi$ is expansive, the dynamical system
$(\Lambda_\A,\varphi)$ is topologically conjugate to a sofic shift.

In [Ma2], the author has introduced the notion of  
symbolic matrix system and
$\lambda$-graph system
from an idea of  the $C^*$-algebras (cf. [Ma],[Ma3]).
The symbolic matrix system is a generalization of symbolic matrix,
 and $\lambda$-graph system is a generalization of $\lambda$-graph 
 (= finite labeled graph).
 We henceforth denote by $\Zp$ the set of all nonnegative integers and by $\Bbb N$ the set of all positive integers.
 A symbolic matrix  system over alphabet $\Sigma$ consists of
two sequences of rectangular matrices
$(\M_{l,l+1}, I_{l,l+1}), l \in \Bbb \Zp$.
The matrices $\M_{l,l+1},l\in \Zp$ 
 have their entries in formal sums of $\Sigma$ and
the matrices $I_{l,l+1},l\in \Zp$ 
have their entries in $\{0,1\}$.
They satisfy the commutation relations:
$$
I_{l,l+1} \M_{l+1,l+2} = \M_{l,l+1}I_{l+1,l+2}, 
\qquad
l \in \Bbb \Zp.
$$
We further assume that
each row of $I_{l,l+1}$ has at least one $1$ 
and each column of $I_{l,l+1}$ has exactly one $1$.
We denote $(\M_{l,l+1}, I_{l,l+1}), l \in \Bbb \Zp$ 
by $(\M,I)$
or $(\M,I^\M)$.
A $\lambda$-graph system 
$ {\frak L} = (V,E,\lambda,\iota)$
consists of a vertex set 
$V = V_0 \cup V_1\cup V_2\cup\cdots$, an edge set 
$E = E_{0,1}\cup E_{1,2}\cup E_{2,3}\cup\cdots$, 
a labeling
$\lambda: E \rightarrow \Sigma$
and a surjective map
$\iota ( = \iota_{l,l+1}): V_{l+1} \rightarrow V_l$ for each
$l\in \Zp.$
It  naturally arises  from a symbolic matrix system $(\M,I)$.
The edges from a vertex $v_i^l \in V_l$ 
to a vertex  $v_j^{l+1} \in V_{l+1}$ are given by the 
$(i,j)$-component $\M_{l,l+1}(i,j)$
of the matrix $\M_{l,l+1}$.
The matrix $I_{l,l+1}$ 
defines a surjection $\iota_{l,l+1}$ 
from $V_{l+1}$ to $V_l$ for each $l \in \Zp$.
The symbolic matrix systems and the $\lambda$-graph systems 
 are the same objects and give rise to subshifts 
 by gathering all the label sequences
 appearing in the labeled Bratteli diagram.
A canonical method to construct a symbolic matrix system and a 
$\lambda$-graph system
from an arbitrary subshift has been introduced in [Ma2].
The obtained symbolic matrix system and the $\lambda$-graph system are said to be canonical for the subshift.
The notion of strong shift equivalence for 
nonnegative matrices and symbolic matrices
has been generalized to symbolic matrix systems as 
 properly strong shift equivalence.
Two symbolic matrix systems  $(\M,I)$ and $(\M^{\prime},I')$
 are said to be {\it properly strong shift equivalent in 1-step}\/ 
 if
 there exist alphabets 
$C,D$ and
 specifications 
$
 \kappa: \Sigma \rightarrow C D,
$
$
 \kappa': \Sigma' \rightarrow D  C
$
 and increasing sequences $n(l),n'(l)$ on $l \in \Zp$
 such that
 for each $l\in \Zp$, there exist
an $n(l)\times n'(l+1)$ matrix ${\P}_l$ over $C,$ 
an $n'(l)\times n(l+1)$ matrix ${\Q}_l$ over $D,$
an $n(l)\times n(l+1)$ matrix $X_l $ over $ \{0,1\}$
and 
an $n'(l)\times n'(l+1)$ matrix $Y_l$ over $ \{0,1\}$
 satisfying the following equations: 
$$
\align
{\M}_{l,l+1} 
\overset{\kappa}\to {\simeq} & {\P}_{2l}{\Q}_{2l+1},\qquad
{\M}'_{l,l+1} 
\overset{\kappa'}\to {\simeq} {\Q}_{2l}{\P}_{2l+1},\\ 
I_{l,l+1} = & X_{2l}X_{2l+1}, \qquad 
I'_{l,l+1} = Y_{2l}Y_{2l+1} 
\endalign
$$
and
$$
X_l{\P}_{l+1}= {\P}_{l}{Y}_{l+1},\qquad
Y_l{\Q}_{l+1}= {\Q}_{l}{X}_{l+1}. 
$$
This situation is written as
$
(\P,\Q,X,Y): (\M,I)\underset{1-pr}\to {\approx} (\M',I').
$
A finite chain 
of properly strong shift equivalences in 1-step with length $N$ 
  is called a
 properly strong shift equivalence (in N-step). 
Then the previously  mentioned Williams's result and Nasu's result  have been generalized to topological conjugacy between subshifts. 
That is,
if two symbolic matrix systems are properly strong shift equivalent,
then their presented subshifts are topologically conjugate.
Furthermore,   
 two subshifts are 
topologically conjugate 
if and only if their canonical symbolic matrix systems are properly  strong shift equivalent
([Ma2]). 
Hence, in particular, 
a properly strong shift equivalence from a symbolic matrix system to itself 
gives rise to an automorphism of the presented subshift.
And an automorphism of a subshift exactly corresponds to a properly strong shift equivalence 
from the canonical symbolic matrix system of the subshift to itself.

In this paper, 
we will generalize the Nasu's textile systems for graph homomorphisms between finite directed (labeled) graphs to graph homomorphisms between
$\lambda$-graph systems and
generalize the Nasu's formalism for topological Markov shifts and sofic shifts to general subshifts.  
Namely we will formulate textile systems for graph homomorphisms between $\lambda$-graph systems and study automorphisms of general subshifts by using the generalized textile systems. 
Let $(\M,I^\M),(\K,I^\K),(\N,I^\N)$ be symbolic matrix systems and
$
{\frak L}^\K,
{\frak L}^\M,{\frak L}^\N
$
their respect 
$\lambda$-graph systems.
Assume that the vertex sets $V_l^\M$ of ${\frak L}^\M$
and 
the vertex sets $V_l^\N$ of ${\frak L}^\N$
coincide
and that the condition
$I^\M_{l,l+1} =I^\N_{l,l+1}
$ hold for all $l \in \Zp.$
We further assume that the vertex set $V_l^\K$ of ${\frak L}^\K$ 
is identified with the edge set 
$E^\N_{l,l+1}$ of ${\frak L}^\N$ for $l\in \Zp.$
A label preserving graph homomorphism
$p :{\frak L}^{\K}\longrightarrow {\frak L}^{\M}$
compatible to $\iota$
is called
a $\lambda$-graph system homomorphism
if $p(V_{l}^\K) =V_l^\M, l\in \Zp$.
A label preserving graph homomorphism
$q :{\frak L}^{\K}\longrightarrow {\frak L}^{\M}$
compatible to $\iota$
is called
a one-shift $\lambda$-graph system homomorphism
if $q(V_l^\K) =V_{l+1}^\M, l\in \Zp$.
Hence the source map 
$s^\K :E^\K_{l,l+1}\longrightarrow V_l^\K = E^\N_{l,l+l}$
and the terminal map 
$t^\K :E^\K_{l,l+1}\longrightarrow V_{l+1}^\K = E^\N_{l+1,l+2}$
of ${\frak L}^\K$ 
yield a  $\lambda$-graph system homomorphism
and a one-shift $\lambda$-graph system homomorphism respectively.
Then
for a  $\lambda$-graph system homomorphism
$p :{\frak L}^{\K}\longrightarrow {\frak L}^{\M}$
and
a  one-shift $\lambda$-graph system homomorphism  
$q :{\frak L}^{\K}\longrightarrow {\frak L}^{\M},$
the diagram
$$
\CD
                   @. {\frak L}^{\M} @. \\
@.                 @AA{p}A           @. \\
{\frak L}^{\N} @<s^\K<< {\frak L}^{\K} @>t^\K>> {\frak L}^{\N}\\
@.                 @VV{q}V           @. \\
                   @. {\frak L}^{\M} @. \\
\endCD
$$
is called a {\it textile system on}\, $\lambda$-{\it graph systems }
if some further conditions are satisfied.
It is written as $\T$.
This formulation is a generalization of Nasu's sofic textile systems [N2].
We will follow  and generalize Nasu's machinery of [N2] so that the dual of  
$\T$ can be defined and we may consider LR textile systems on $\lambda$-graph systems.
We will prove that
 for a forward automorphism $\phi$ of a subshift
$(\Lambda,\sigma)$, the automorphisms $\phi^k \sigma^n, k\ge 0, n\ge 1$
can be explicitly realized as a subshift defined by 
certain symbolic matrix systems coming from 
both  the strong shift equivalence representing $\phi$ 
and the subshift  $(\Lambda,\sigma)$. 
Suppose that $\Lambda$ is equipped with a metric for which $\sigma$ has $1$
as its expansive constant.
If in particular, $\phi$ is expansive with $\frac{1}{m}$
as its expansive constant for some $m \in \Bbb N$,
the dynamical system
$(\Lambda,\phi)$ 
can be realized as a subshift defined by certain symbolic matrix system coming from the strong shift equivalence representing $\phi$ and the subshift 
$(\Lambda,\sigma)$ (Theorem 7.6). 

We will prove the following 
\proclaim{Theorem(Theorem 7.8)}
Let $(\Lambda, \sigma)(=(\Lambda_{\M}, \sigma_{\M}))$ 
be a subshift presented by a symbolic matrix system 
$(\M,I)$.
Let $\phi$ be a forward automorphism on $(\Lambda, \sigma)$
defined by a properly strong shift equivalence 
$$
(\P^{(j)},\Q^{(j)},X^{(j)},Y^{(j)}):
 (\M^{(j-1)},I^{(j-1)})\underset{1-pr}\to {\approx} (\M^{(j)},I^{(j)}),\quad
 j=1,2,\dots,N
$$
in $N$-step
where $(\M^{(0)},I^{(0)})=(\M^{(N)},I^{(N)}) = (\M,I)$.
Then the dynamical system
$(\Lambda,\phi^k\sigma^n)$
is topologically conjugate to the subshift 
$(\Lambda_{\P^k \M^n},\sigma_{\P^k \M^n})$
presented by the symbolic matrix system
$(\P^k \M^n,I^{kN+n})$
 for
 $k\ge 0, n\ge 1$
defined by
$$
\align
& {({\P}^k{\M}^n)}_{l,l+1}\\ 
=& 
{\P}_{l(kN+n),l(kN+n)+N}{\P}_{l(kN+n) +N,l(kN+n)+2N}\cdots 
{\P}_{l(kN+n)+(k-1)N,l(kN+n)+kN}\cdot\\
& \cdot {\M}_{l(kN+n)+kN,l(kN+n)+kN+1}{\M}_{l(kN+n)+kN+1,l(kN+n)+kN+2}\cdots
{\M}_{(l+1)(kN+n)-1,(l+1)(kN+n)},\\
& I^{kN+n}_{l,l+1} \\ 
=& 
I_{l(kN+n),l(kN+n)+1}I_{l(kN+n)+1,l(kN+n)+2}\cdots 
I_{(l+1)l(kN+n)-1,(l+1)(kN+n)},\qquad l \in \Zp
\endalign
$$
where
$
\P_{l,l+1} 
=\P^{(1)}_{2l}Y^{(1)}_{2l+1}\P^{(2)}_{2l+2}Y^{(2)}_{2l+3}
\cdots \P^{(N)}_{2l+2N-2}Y^{(N)}_{2l+2N-1}
$
and $\P^{(i)}_{2l+2(i-1)},Y^{(i)}_{2l+2i-1}, i=1,\dots,N$ 
are matrices 
appearing in the above properly strong shift equivalence.
\endproclaim
                                  
Namely these automorphisms $\phi^k \sigma^n, k\ge 0, n\ge 1$
are subshift-identified.
As an application of this result, 
if an automorphism $\phi$ of a subshift 
$(\Lambda, \sigma)$ is a simple automorphism, 
that is conjugate to a symbolic automorphism fixing vertices of a 
$\lambda$-graph system,
the dynamical system
$(\Lambda, \phi \circ \sigma^n)$ 
is topologically conjugate to the $n$-th power 
$
(\Lambda, \sigma^n)
$
of the subshift  
$
(\Lambda, \sigma)
$
for $n \in \Bbb Z, n\ne 0$ 
(Theorem 8.2).
\medskip
 
This paper is organized as in the following way.

\medskip

1. Introduction

2. Symbolic matrix systems and $\lambda$-graph systems

3. Textile systems on $\lambda$-graph systems

4. Textile shifts on $\lambda$-graph systems

5. LR textile systems on $\lambda$-graph systems
   
6. LR textile systems and properly strong shift equivalences 

7. Subshift-identifications of automorphisms of subshifts 

8. An application

\heading 2.  Symbolic matrix systems and $\lambda$-graph systems 
\endheading
We call each element of a finite set $\Sigma$ a symbol or a label.
 The transformation $\sigma$ on the infinite product space
$\Sigma^{\Bbb Z}$ 
given by 
$(\sigma({(x_i)}_{i\in \Bbb Z}) = {(x_{i+1})}_{i\in \Bbb Z}$ 
is called the (full) shift.
 Let $\Lambda$ be a shift-invariant closed subset of $\Sigma^{\Bbb Z}$ i.e. $\sigma(\Lambda) = \Lambda$.
 We  write the subshift $(\Lambda,\sigma)$
 as $\Lambda$ for short.
We denote by 
$
{\Lambda}^+
( \subset
\Sigma^{\Bbb N} )
$ 
the set of all right-infinite sequences that appear in $\Lambda$.
   A finite sequence 
$\mu = (\mu_1,...,\mu_k) $ of elements $\mu_j \in \Sigma$ is called a block or a word of length $k$. 
We write the empty symbol $\emptyset$ in  $\Sigma$ 
as $0$. 
We denote by ${\frak S}_{\Sigma}$
the set of all finite formal sums of elements of $\Sigma$.
By a symbolic matrix  $\A$ over $\Sigma$
we mean a finite matrix with entries in ${\frak S}_{\Sigma}.$
A square symbolic matrix $\A$ naturally gives rise to a labeled directed graph, called a $\lambda$-graph,
 which we denote by
$G_{\A}.$
The labeled directed graph defines a subshift over $\Sigma$
which consist of all infinite labeled sequences following the labeled edges
in $G_{\A}.$
Such a subshift is called a sofic shift presented by $G_{\A}$
(cf. [Fi], [Kr], [Kr2], [We], [Ki], [LM]).
If, in particular, different edges have different labels,
the sofic shift is called a topological Markov shift.

Let $\A, \A'$ be 
symbolic matrices
over $\Sigma, \Sigma'$
respectively
 and 
 $\kappa$
a bijection
  from a subset of $\Sigma$ onto a subset of $\Sigma'$.
Following M. Nasu in [N],[N2],  
 we say that $\A$ and $\A'$ are  specified equivalence under specification
 $\kappa$ if $\A'$ can be obtained from $\A$
  by replacing every symbol 
 $a$ appearing in the components of $\A$ by $\kappa(a)$.
 We write it as
 $\A \overset{\kappa}\to{\simeq} \A'$.

Two symbolic matrix systems 
$(\M,I)$ over  $\Sigma$
and
$(\M',I')$ over  $\Sigma'$
are said to be isomorphic 
if the size of $\M_{l,l+1}$ coincides with that of $\M'_{l,l+1}$
for $l \in \Zp$ and there exist a specification $\kappa$ 
from $\Sigma$ to $\Sigma'$ 
and an $m(l) \times m(l)$-square permutation matrix
$S_l$ for each $l \in \Zp $ such that
$$
S_l \M_{l,l+1} \overset{\kappa}\to{\simeq} \M'_{l,l+1}S_{l+1},
\qquad
S_l I_{l,l+1} = I'_{l,l+1}S_{l+1}.
$$
Recall 
that a $\lambda$-graph system 
$\frak L =(V,E,\lambda,\iota)$ over $\Sigma$
is a directed Bratteli diagram with vertex set
$$
V = \cup_{l \in \Zp} V_{l},
$$
and edge set
$$
E = \cup_{l \in \Zp} E_{l,l+1},
$$
that is labeled by a map 
$\lambda ( = \lambda_{l,l+1}) : E_{l,l+1} \rightarrow \Sigma$
with symbols  in  $\Sigma$ for $l \in \Zp$,
and that is supplied with a sequence of surjective maps
$$
\iota ( = \iota_{l,l+1}):V_{l+1} \rightarrow V_l\qquad \text{ for } \quad l \in \Bbb \Zp.
$$
Here the $V_{l},l \in \Zp$
are finite disjoint sets.   
Also  
$E_{l,l+1},l \in \Zp$
are finite disjoint sets.
An edge $e$ in $E_{l,l+1}$ has its source vertex $s(e)$ in $V_{l}$ 
and its terminal vertex  $t(e)$  in
$V_{l+1}$.
Every vertex in $V$ has a successor and  every 
vertex in $V$, 
except the verteces in $V_0$ at level $0$, has a predecessor. 
It is then required that there exists an edge in $E_{l,l+1}$
with label $\alpha$ and its terminal is  $v \in V_{l+1}$
 if and only if 
 there exists an edge in $E_{l-1,l}$
with label $\alpha$ and its terminal is $\iota(v) \in V_{l}.$
A $\lambda$-graph system is said to be essential 
if there is no distinct edges
that have the same source vertices, the same terminal vertices
and the same labels.
Throughout this paper, 
we will treat essential $\lambda$-graph systems.
For 
$u \in V_{l-1}$ and
$v \in V_{l+1},$
we put
$$
\align
E^{\iota}_{l,l+1}(u,v)
& = \{e \in E_{l,l+1} \mid t(e) = v, \iota(s(e)) = u \},\\
E_{\iota}^{l-1,l}(u,v)
& = \{e \in E_{l-1,l} \mid s(e) = u, t(e) = \iota(v) \}.
\endalign
$$ 
Then there exists a bijective map
$\varphi^{\frak L}_{(u,v)}$ 
from
$
E^{\iota}_{l,l+1}(u,v)
$
to
$
E_{\iota}^{l-1,l}(u,v)
$
such that 
$$
\lambda(
\varphi^{\frak L}_{(u,v)}(e) )
= \lambda(e)
\qquad
\text{ for } e \in E^{\iota}_{l,l+1}(u,v).
$$ 
Hence two sets 
$
E^{\iota}_{l,l+1}(u,v)
$
and
$
E_{\iota}^{l-1,l}(u,v)
$
 bijectively correspond
preservings labels
for all pairs $(u,v) \in V_{l-1}\times V_{l+1}$.
We call this property the local property of the $\lambda$-graph system.
We immediately see
\proclaim{Lemma 2.1}
For a $\lambda$-graph system
$\frak L = (V,E,\lambda,\iota)$ over
$\Sigma$,
there exists  
a surjection 
$$
\varphi_l^{\frak L}:E_{l,l+1}\longrightarrow E_{l-1,l}
$$
 for each $l \in \Bbb N$ such that 
$$
\varphi_l^{\frak L}|_{E^{\iota}_{l,l+1}(u,v)} =
\varphi^{\frak L}_{(u,v)}
\qquad
\text{ for } u \in V_{l-1}, v \in V_{l+1}
$$
and
$$
\iota_{l-1,l}(s(e)) = s (\varphi_l^{\frak L}(e)),
\quad
\iota_{l,l+1}(t(e)) = t (\varphi_l^{\frak L}(e))
\qquad \text{ for } e \in E_{l,l+1}.
$$
\endproclaim

We call an edge $e \in E_{l,l+1}$ a $\lambda$-edge and
a connecting finite sequence of $\lambda$-edges a $\lambda$-path.
For 
$u\in V_{l}$ and
$v \in V_{l+1},$
if $\iota(v) = u$,
we say that there exists an $\iota$-edge from
$u$ to
$v.$
Similarly we use the term $\iota$-path.

Two $\lambda$-graph systems
$(V,E,\lambda,\iota)$ over $\Sigma$ 
and
$(V',E',\lambda',\iota')$ over $\Sigma'$ 
are said to be isomorphic
if there exist bijections
$\varPhi_V : V_l \rightarrow V_l'$,
$\varPhi_E : E_{l,l+1} \rightarrow E_{l,l+1}'$
for $l \in \Bbb \Zp,$
and a specification
$\kappa : \Sigma \rightarrow \Sigma'$
such that
$
\varPhi_V(s(e)) = s(\varPhi_E(e)),
$
$ 
\varPhi_V(t(e)) = t(\varPhi_E(e)) 
$
and
$\lambda'(\varPhi_E(e)) = \kappa(\lambda(e)) 
$
for  $e \in E,$ 
and 
$\iota'(\varPhi_V(v)) = \varPhi_V(\iota(v)) 
$
for  $v \in V.$ 
 There exists a bijective correspondence 
between 
the set of all isomorphism classes of symbolic matrix systems 
and 
the set of all isomorphism classes of $\lambda$-graph systems.
We identify isomorphic symbolic matrix systems,  
and similarly isomorphic $\lambda$-graph systems.

A symbolic matrix system
$(\M,I)$ is denoted by 
$(\M,I^\M)$ although the matrices $I_{l,l+1}$ are not determined by the symbolic matrices
$\M_{l,l+1}, l\in \Zp.$
We denote its $\lambda$-graph system by
$
{\frak L}^{\M} =
(V^{\M},E^{\M},\lambda^{\M},\iota^{\M}).
$
The surjections
$
\varphi_l^{{\frak L}^\M}:E_{l,l+1}^\M\longrightarrow E_{l-1,l}^\M, l\in \Zp
$
defined in Lemma 2.1 are denoted by 
$
\varphi_l^\M, l\in \Zp
$.

A $\lambda$-graph (= a finite labeled graph)  defines a $\lambda$-graph system as in the following way.  
Let   
${\Cal G} =(G,\lambda)$ be a  
$\lambda$-graph with underlying finite directed graph $G$ and its labeling $\lambda$.
Let $V^G$
be the vertex set of $G$.
Put
$V_l =  V^G$ for all $l \in \Zp$
and $\iota = \id.$
Write labeled edges from $V_l$ to $V_{l+1}$ for $l \in \Zp$
 following
the directed graph $G$ with labeling $\lambda$.
It is clear  to see that
the resulting labeled Bratteli diagram with $\iota(=\id)$ becomes
a $\lambda$-graph system.
A $\lambda$-graph and also a $\lambda$-graph system 
are said to be left-resolving 
if different edges with the same label must have different terminals.
By the construction
if the $\lambda$-graph is left-resolving, so is the $\lambda$-graph system.
In what follows, we assume that 
a $\lambda$-graph system is left-resolving.

For a $\lambda$-graph system 
$\frak L = (V,E,\lambda,\iota)$
over 
$\Sigma$ and a natural number $N \ge 2$,
the $N$-heigher block ${\frak L}^{[N]}$ 
of $\frak L$ is defined to be  
a $\lambda$-graph system
$(V^{[N]},E^{[N]},\lambda^{[N]},\iota^{[N]})$
over 
$
\Sigma^{[N]}= \undersetbrace{\text{N-times}}\to
{\Sigma \cdots  \Sigma}
$
as follows ([Ma2]):
$$
\align
V_l^{[N]} = \{ ( e_1,e_2,\dots, e_{N-1})& 
\in E_{l,l+1}\times E_{l+1,l+2}\times \cdots \times E_{l+N-2,l+N-1} \mid \\
          t(e_i) & = s(e_{i+1}) \text{ for } i =1,2,\dots, N-2 \}, \\
E_{l,l+1}^{[N]}
          = \{((  e_1,\dots, e_{N-1}),&(f_1, \dots, f_{N-1}))   
                 \in V_{l}^{[N]} \times V_{l+1}^{[N]} \mid \\
         t(e_{N-1}) & = s(f_{N-1}),  e_{i+1}  = f_i  \text{ for } i =1,2,\dots, N-2 \}. 
\endalign
$$
The maps 
$$
s^{[N]}: E_{l,l+1}^{[N]} \rightarrow V_l^{[N]},
\qquad 
t^{[N]}: E_{l,l+1}^{[N]} \rightarrow V_{l+1}^{[N]}
$$
are defined by
$$
\align
s^{[N]}((e_1,\dots, e_{N-1}),(f_1, \dots, f_{N-1})) & 
            = (e_1,\dots,e_{N-1}),\\   
t^{[N]}((e_1,\dots, e_{N-1}),(f_1, \dots, f_{N-1})) & 
            = (f_1,\dots,f_{N-1}).
\endalign
$$
Set
$
V^{[N]} = \cup_{l \in \Zp} V_l^{[N]}
$
and
$
E^{[N]} = \cup_{l \in \Zp} E_{l,l+1}^{[N]}.
$
Hence
$(V^{[N]},E^{[N]},s^{[N]},t^{[N]} )$
is a Bratteli diagram.
A labeling
$\lambda^{[N]}$ on
$(V^{[N]},E^{[N]})$ is defined by   
$$
\lambda^{[N]}((e_1,\dots, e_{N-1}),(f_1, \dots, f_{N-1})) 
 = \lambda(e_1)\lambda(e_2)\dots \lambda(e_{N-1})\lambda(f_{N-1}) \in
 \Sigma^{[N]}
$$   
for
$((e_1,\dots, e_{N-1}),(f_1, \dots, f_{N-1})) \in E^{[N]}.$
A sequence of surjections
$
\iota^{[N]} : V_{l+1}^{[N]} \rightarrow V_l^{[N]}, l \in \Zp
$   
is defined as follows.
As
the $\lambda$-graph system
$(V,E,\lambda,\iota)$ 
is left-resolving,
for 
$(e_1,\dots,e_{N-1}) \in V_{l+1}^{[N]},$
there uniquely exist
$e_i' \in E_{l+i-1,l+i}$
for 
$i=1,2,\dots,N-2$
such that
$$
\iota(s(e_i)) = s(e'_i),\quad
\iota(t(e_i)) = t(e'_i),\quad
\lambda(e_i) = \lambda(e'_i).
$$
As 
$(e'_1,\dots,e'_{N-1}) \in V_{l}^{[N]},$
by setting
$
\iota^{[N]}(e_1,\dots,e_{N-1}) =
(e'_1,\dots,e'_{N-1}),
$ 
we get a $\lambda$-graph system
$
(V^{[N]},E^{[N]},\lambda^{[N]},\iota^{[N]})
$
 over $\Sigma^{[N]}$.
We set 
${\frak L}^{[1]} =\frak L$.
The $N$-heigher 
block 
$(\M^{[N]},I^{[N]})$
of 
a symbolic matrix system
$(\M,I)$ is defined to be the symbolic matrix system for
the $N$-heigher block ${\frak L}^{[N]}$ 
of the  $\lambda$-graph system $\frak L$
for
$(\M,I)$.

\heading 3. Textile systems on $\lambda$-graph systems
\endheading

In what follows, 
let
$(\K,I^{\K})$, $(\M,I^{\M})$ and 
$(\N,I^{\N})$
be  symbolic matrix systems 
over alphabets  
$\Sigma^{\K}$, $\Sigma^{\M}$
and
$\Sigma^{\N}$
respectively.
Let us  consider their respect $\lambda$-graph systems
$
{\frak L}^{\K} =
(V^{\K},E^{\K},\lambda^{\K},\iota^{\K}),
$
$
{\frak L}^{\M} =
(V^{\M},E^{\M},\lambda^{\M},\iota^{\M})
$
and
$
{\frak L}^{\N} =
(V^{\N},E^{\N},\lambda^{\N},\iota^{\N}).
$
We denote by
$(s^\K,t^\K),(s^\M,t^\M)$ and $(s^\N,t^\N)$
their source maps and terminal maps in the $\lambda$-graph systems 
respectively.

A $\lambda$-{\it graph system homomorphism}\,
$p=(p^V,p^E,p^{\Sigma}) :
{\frak L}^{\K}\longrightarrow {\frak L}^{\M}$
consists of sequences of maps
$$
p^V (= p^V_l) : V^{\K}_l\longrightarrow V^{\M}_l,
\qquad
p^E (= p^E_{l,l+1}) : E^{\K}_{l,l+1}\longrightarrow E^{\M}_{l,l+1},
\qquad
l \in \Zp
$$
together with a map
$
p^{\Sigma}: \Sigma^{\K}\longrightarrow \Sigma^{\M}
$
such that 
\roster
\item
$p^V_l(s^{\K}(e)) = s^{\M}(p^E_{l,l+1}(e)),
\quad
p^V_{l+1}(t^{\K}(e)) = t^{\M}(p^E_{l,l+1}(e))
$
for 
$
e \in E_{l,l+1}^{\K},
$
\item
$p^V_l(\iota^{\K}_{l,l+1}(v) = \iota^{\M}_{l,l+1}(p^V_{l+1}(v))
$
for
$v \in V^{\K}_{l+1},
$
\item
$p^{\Sigma}(\lambda^{\K}(e)) = \lambda^{\M}(p^E_{l,l+1}(e))
$
for
$ 
e \in E^{\K}_{l,l+1}.
$
\endroster

A {\it one-shift} $\lambda$-{\it graph system homomorphism}\/
$q=(q^V,q^E,q^{\Sigma}) :
{\frak L}^{\K}\longrightarrow {\frak L}^{\M}$
consists of sequences of maps
$$
q^V (= q^V_{l}) : V^{\K}_{l}\longrightarrow V^{\M}_{l+1},
\qquad
q^E (= q^E_{l,l+1}) : E^{\K}_{l,l+1}\longrightarrow E^{\M}_{l+1,l+2},
\qquad
l \in \Zp
$$
together with a map
$
q^{\Sigma}: \Sigma^{\K}\longrightarrow \Sigma^{\M}
$
such that 
\roster
\item
$q^V_l(s^{\K}(e)) = s^{\M}(q^E_{l,l+1}(e)),
\quad
q^V_{l+1}(t^{\K}(e)) = t^{\M}(q^E_{l,l+1}(e))
$
for 
$
e \in E_{l,l+1}^{\K},
$
\item
$q^V_{l}(\iota^{\K}_{l,l+1}(v) 
= 
\iota^{\M}_{l+1,l+2}(q^V_{l+1}(v))
$
for
$v \in V^{\K}_{l+1},
$
\item
$q^{\Sigma}(\lambda^{\K}(e)) = \lambda^{\M}(q^E_{l,l+1}(e))
$
for
$ 
e \in E^{\K}_{l,l+1}.
$
\endroster
For a  one-shift  $\lambda$-graph system homomorphism 
$q=(q^V,q^E,q^{\Sigma}) :
{\frak L}^{\K}\longrightarrow {\frak L}^{\M}$,
put
$
{p^V_q}_{l} = \iota_{l,l+1}^{\M}\circ q^V_{l}:
 V^{\K}_{l}\longrightarrow V^{\M}_l,
{p^E_q}_{l,l+1} = \varphi_{l+1,l+2}^{\M}\circ q^E_{l,l+1}: 
E^{\K}_{l,l+1}\longrightarrow E^{\M}_{l,l+1},
l \in \Zp,
$
and
$
p^{\Sigma}_q =q^{\Sigma}: 
\Sigma^{\K}\longrightarrow \Sigma^{\M}.
$
Then
$p_q=(p^V_q, p^E_q, p^{\Sigma}_q) :
{\frak L}^{\K}\longrightarrow {\frak L}^{\M}$
is a $\lambda$-graph system homomorphism.

\proclaim{Lemma 3.1}
\roster
\item"(i)"
For a $\lambda$-graph system homomorphism 
$p=(p^V,p^E,p^{\Sigma}) :
{\frak L}^{\K}\longrightarrow {\frak L}^{\M}
$
we have
$
p^E \circ \varphi^{\K}_l = 
\varphi^{\M}_l \circ p^E.
$
\item"(ii)"
For a one-shift $\lambda$-graph system homomorphism 
$q=(q^V,q^E,q^{\Sigma}) :
{\frak L}^{\K}\longrightarrow {\frak L}^{\M}$
we have
$
q^E \circ \varphi^{\K}_l = 
\varphi^{\M}_l \circ q^E.
$
\endroster
\endproclaim
\demo{Proof}
Let $p : {\frak L}^{\K}\longrightarrow {\frak L}^{\M}$ be a $\lambda$-graph system homomorphism. 
For $u \in V_{l-1}^\K, v \in V_{l+1}^\K$
and 
$e \in {E^\K}_{l,l+1}^{\iota^\K}(u,v)$,
it is streightforward to see that 
$
s^\M(p^E(
\varphi^{\K}_l (e))) = 
s^\M(\varphi^{\M}_l ( p^E(e))),
\,
t^\M(p^E(
\varphi^{\K}_l (e))) = 
t^\M(\varphi^{\M}_l ( p^E(e)))
$
and
$
\lambda^\M(p^E(
\varphi^{\K}_l (e))) = 
\lambda^\M(\varphi^{\M}_l ( p^E(e)))
$.
As ${\frak L}^\M$ is essential, one sees that
$
p^E(
\varphi^{\K}_l (e)) = 
\varphi^{\M}_l ( p^E(e)).
$
Hence (i) holds.
The assertion for (ii) is similarly shown to (i).
\qed
\enddemo

We say that 
$
{\frak L}^{\M}
$ and
$
{\frak L}^{\N}
$
{\it form squares }
if
$$
V_l^{\M} = 
V_l^{\N},
\qquad
I_{l,l+1}^{\M} = 
I_{l,l+1}^{\N},
\qquad
l \in \Zp. \tag 3.1
$$
In this case, one may see a square as in the following figure:
$$
\CD
V_l^\N = V_l^\M    @>E^\M_{l,l+1}>> V_{l+1}^\M = V_{l+1}^\N \\
@VV{E^\N_{l,l+1}}V @VV{E^\N_{l+1,l+2}}V \\   
V_{l+1}^\N = V_{l+1}^\M@>E^\M_{l+1,l+2}>> V_{l+2}^\M = V_{l+2}^\N\\
\endCD
\qquad\qquad\quad
\text{ for }
\quad
l \in \Zp.
$$
We will formulate 
textile system on $\lambda$-graph systems
as in the following way.
\noindent
{\bf Definition ( Textile system on $\lambda$-graph systems).}
For $\lambda$-graph systems
$
{\frak L}^\K,
{\frak L}^\M,{\frak L}^\N
$
with
a  $\lambda$-graph system homomorphism
$p :{\frak L}^{\K}\longrightarrow {\frak L}^{\M}$
and
a one-shift $\lambda$-graph system homomorphism  
$q :{\frak L}^{\K}\longrightarrow {\frak L}^{\M},$
the diagram
$$
\CD
                   @. {\frak L}^{\M} @. \\
@.                 @AA{p}A           @. \\
{\frak L}^{\N} @<s^\K<< {\frak L}^{\K} @>t^\K>> {\frak L}^{\N}\\
@.                 @VV{q}V           @. \\
                   @. {\frak L}^{\M} @. \\
\endCD
$$
is called a {\it textile system on} \, $\lambda$-{\it graph systems }\,
if the following six conditions are satisfied.
\roster
\item
$
{\frak L}^{\M}
$ 
and
$
{\frak L}^{\N}
$
form squares.
\item
$V_l^{\K} = E^{\N}_{l,l+1},\qquad l \in \Zp.$
\item
Under the equality (2),
$$
(\iota^{\K}_{l,l+1}:V^{\K}_{l+1}\rightarrow V^{\K}_{l})
=
(\varphi^{\N}_{l+1}:E^{\N}_{l+1,l+2}\rightarrow
E^{\N}_{l,l+1}),\quad l \in \Zp.
$$
\item
Under the equalities (3.1) and (2),
$$
\align
(p^V:V^{\K}_l \rightarrow V^{\M}_{l}) 
= &
(s^{\N}: E^{\N}_{l,l+1}\rightarrow
V^{\N}_{l} ), \\
(q^V:V^{\K}_l \rightarrow V^{\M}_{l+1}) 
= &
(t^{\N}: E^{\N}_{l,l+1}\rightarrow
V^{\N}_{l+1} ), \qquad
l \in \Zp.
\endalign
$$ 
\item
The quadruple:
$$
(s^\K(e),t^\K(e),p^E(e),q^E(e)) 
\in 
V^{\K}_{l}\times V^{\K}_{l+l}\times 
E^{\M}_{l,l+1}\times E^{\M}_{l+1,l+2}
$$ 
determines 
$e \in
 E^{\K}_{l,l+l}.
 $
\item
Under the equality (2),
there exists a specified equivalence
between
$
\Sigma^{\N}\times \Sigma^{\N}\times 
\Sigma^{\M}\times \Sigma^{\M}
$
and
$\Sigma^{\K}$ 
by the correspondence between the symbols:
$$
(\lambda^{\N}(s^\K(e)),\lambda^{\N}(t^\K(e)),
\lambda^{\M}(p^E(e)),\lambda^{\M}(q^E(e))) 
\in 
\Sigma^{\N}\times \Sigma^{\N}\times 
\Sigma^{\M}\times \Sigma^{\M}
$$ 
and
$
\lambda^{\K}(e) \in
 \Sigma^{\K}.
 $
\endroster
A textile system on $\lambda$-graph systems is called 
a {\it textile } $\lambda$-{\it graph system } for short. 
We write the textile $\lambda$-graph system
as
$  
{\Cal T}_{{\K}^{\M}_{\N}} 
=(p,q:{\frak L}^{\K}\rightarrow{\frak L}^{\M})
$,
or simply as $\Cal T$.
In viewing the textile $\lambda$-graph system,
one uses the following square
$$
\CD
 \cdot   @>\lambda^{\M}(p^E(e))>> \cdot \\
@V{\lambda^{\N}(s^\K(e))}VV @VV{\lambda^{\N}(t^\K(e))}V \\   
\cdot @>\lambda^{\M}(q^E(e))>> \cdot \\
\endCD
\qquad\qquad\quad
\text{ for }
\quad
e 
\in E^\K.
$$
\proclaim{Proposition 3.2}
For 
a textile  $\lambda$-graph system
${\Cal T}_{{\K}^{\M}_{\N}} 
=(p,q:{\frak L}^{\K}\rightarrow{\frak L}^{\M})
$,
there exists a $\lambda$-graph system ${\frak L}^{\K^*}$ and 
a textile $\lambda$-graph system
${\Cal T}_{{\K^*}^{\N}_{\M}} 
=(s^\K, t^\K:{\frak L}^{\K^*}\rightarrow{\frak L}^{\N})
$
defined by the diagram:
$$
{\CD
{\Cal T}_{{{\K}^*}^{\N}_{\M}} @. {\frak L}^{\N} @. \\
@.                 @AA{s^\K}A           @. \\
{\frak L}^{\M} @<p<< {\frak L}^{{\K}^*} @>q>> {\frak L}^{\M}\\
@.                 @VV{t^\K}V           @. \\
                   @. {\frak L}^{\N} @. \\
\endCD}
$$
\endproclaim
\demo{Proof}
We define a $\lambda$-graph system
$
{\frak L}^{{\K}^*} =
(V^{{\K}^*},E^{{\K}^*},\lambda^{{\K}^*},\iota^{{\K}^*})
$
over 
$\Sigma^{{\K}^*}$
by setting
$$
V^{\K^*}_l  = E_{l,l+1}^\M,\quad
E^{\K^*}_{l,l+1} = E_{l,l+1}^\K,\quad
\iota^{\K^*}_{l,l+1} = \varphi_{l+1}^\M,\quad
\Sigma^{{\K}^*} = \Sigma^\K \qquad\text{ for }l \in \Zp
$$
and for $e \in  E^{\K^*}_{l,l+1} = E_{l,l+1}^\K$
$$
s^{\K^*}(e) = p^E(e)\in E_{l,l+1}^\M = V^{\K^*}_l,\quad
t^{\K^*}(e) = q^E(e)\in E_{l+1,l+2}^\M = V^{\K^*}_{l+1},\quad
\lambda^{\K^*}(e) = \lambda^{\K}(e).
$$
For 
$u \in V_{l-1}^{\K^*}, v \in V_{l+1}^{\K^*}$
put
$
w=\iota_{l,l+1}^{\K^*}(v). 
$
It then follows that 
$$
\align
E^{l-1,l}_{\iota^{\K^*}}(u,v) 
& = \{ e \in E^{\K^*}_{l-1,l} \mid s^{\K^*}(e) = u, \, t^{\K^*}(e) = w \}\\
& = \{ e \in E^{\K}_{l-1,l} \mid p^E(e) = u, \, q^E(e) = w \},\\
E^{\iota^{\K^*}}_{l,l+1}(u,v) 
& = \{ f \in E^{\iota^{\K^*}}_{l,l+1} \mid 
\iota^{\K^*}_{l-1,l}(s^{\K^*}(f))= u, \, t^{\K^*}(f) = v \}\\
& = \{ f \in E^{\K}_{l,l+1} \mid \varphi_{l}^\M p^E(f) = u, \,
 q^E(f) = v \}.
\endalign
$$
For 
$e \in E^{l-1,l}_{\iota^{\K^*}}(u,v),$
one sees 
$t^\K(e) \in E_{l,l+1}^\N$.
As ${\frak L}^{\N}$ is left-resolving,
there uniquely exists
$v' \in E_{l+1,l+2}^\N =V^\K_{l+1}$ 
such that 
$$
\lambda^\N(v') = \lambda^\N(t^\K(e)),\qquad
t^\N(v') = t^\M(v),\qquad
\varphi_{l+1}^\N(v') = t^\K(e).
$$
For the two vertices
$s^\K(e) \in V_{l-1}^\K, v' \in V_{l+1}^\K$
with
$\iota^\K(v') = t(e) 
$,
by the local property of ${\frak L}^\K$,
there uniquely exists
$f \in E_{l,l+1}^\K$ such that 
$$\iota_{l-1,l}^\K(s^\K(f))=s^\K(e),\qquad t^\K(f) = v',\qquad 
\lambda^\K(f) = \lambda^\K(e).
$$ 
$$
\align
{
\CD
  \cdot  @>p^E(e) = u \in E_{l-1,l}^\M>>  \cdot\\
@V{s^\K(e) \in E^\N_{l-1,l}}VV    @VV{t^\K(e) \in E^\N_{l,l+1}}V \\   
\cdot @>q^E(e) = w \in E_{l,l+1}^\M   >> \cdot\\
\endCD}& \\
&{
\CD
  \cdot  @>p^E(f) = u \in E_{l,l+1}^\M>>  \cdot\\
@V{s^\K(f) \in E^\N_{l,l+1}}VV       @VV{t^\K(f) \in E^\N_{l+1,l+2}}V \\   
\cdot @>q^E(f) = v \in E_{l+1,l+2}^\M   >> \cdot\\
\endCD}
\endalign
$$
Hence one has
$$
\align
&(\lambda^\N(s^\K(f)), \lambda^\N(t^\K(f)),
\lambda^\M(p^E(f)),\lambda^\M(q^E(f))\\
=
&(\lambda^\N(s^\K(e)), \lambda^\N(t^\K(e)),
\lambda^\M(p^E(e)),\lambda^\M(q^E(e)).
\endalign
$$
Since 
${\frak L}^\M$ is left-resolving,
the edge $p^E(f)$, whose label is $\lambda^\M(p^E(e))$ and terminal is the source of $v'\in E_{l+1,l+2}^\N$, is unique,
and also 
the edge $q^E(f)$, whose label is $\lambda^\M(q^E(e))$ and terminal is the terminal of $v'\in E_{l+1,l+2}^\N$, is unique,
Since 
${\frak L}^\N$ is left-resolving,
the edge $s^\K(f)$, whose label is $\lambda^\N(s^\K(e))$ and terminal is the source of $v\in E_{l+1,l+2}^\M$, is unique,
and also 
the edge $t^\K(f)$, whose label is $\lambda^\N(t^\K(e))$ and terminal is the terminal of $v\in E_{l+1,l+2}^\M$, is unique,
Hence the square
$$
\CD
  \cdot  @>p^E(f) = u \in E_{l,l+1}^\M>>  \cdot\\
@V{s^\K(f) \in E^\N_{l,l+1}}VV       @VV{t^\K(f) \in E^\N_{l+1,l+2}}V \\   
\cdot @>q^E(f) = v \in E_{l+1,l+2}^\M   >> \cdot\\
\endCD
$$
is uniquely determined by 
$e \in E_{\iota^{\K^*}}^{l-1,l}(u,v)$
so that
$f\in E_{l,l+1}^{\iota^{\K^*}}(u,v)$
is uniquely determined and hence
$\varphi_l^\K(f) = e.$
Conversely, for an edge
$f\in E_{l,l+1}^{\iota^{\K^*}}$
there uniquely exists
$e \in E_{\iota^{\K^*}}^{l-1,l}(u,v)$
such that
$\varphi_l^\K(f) = e.$
Hence 
$
(V^{{\K}^*},E^{{\K}^*},\lambda^{{\K}^*},\iota^{{\K}^*})
$
satisfies the local property so that 
it yields a $\lambda$-graph system over $\Sigma^{\K}$,
that is written as
$
{\frak L}^{{\K}^*}. 
$

We define a $\lambda$-graph system homomorphism
$
p^*:
{\frak L}^{{\K}^*} \rightarrow {\frak L}^\N
$
by setting
$$
p^*(e) = s^\K(e) \in E_{l,l+1}^\N
\qquad
\text{ for }\quad 
e \in E^{\K^*}_{l,l+1}
$$
and
a one-shift $\lambda$-graph system homomorphism
$
q^*:
{\frak L}^{{\K}^*} \rightarrow {\frak L}^\N
$
by setting
$$
q^*(e) = t^\K(e) \in E_{l+1,l+2}^\N
\qquad
\text{ for }\quad 
e \in E^{\K^*}_{l,l+1}.
$$
Then the diagram bellow 
$$
\CD
@. {\frak L}^{\N} @. \\
@.                 @AA{p^*=s^\K }A           @. \\
{\frak L}^{\M} @<s^{K^*}=p<<{\frak L}^{{\K}^*} @>t^{\K^*}=q>> {\frak L}^{\M}\\
@.                 @VV{q^*=t^\K}V           @. \\
                   @. {\frak L}^{\N} @. \\
\endCD
$$
yields a textile  $\lambda$-graph system
$
{\Cal T}_{{\K^*}^{\N}_{\M}} 
=(s^\K,t^\K:{\frak L}^{\K^*}\rightarrow{\frak L}^{\N})
$.
\qed
\enddemo
We call 
the 
textile $\lambda$-graph system
${\Cal T}_{{\K^*}^{\N}_{\M}} 
=(s^\K, t^\K :{\frak L}^{\K^*}\rightarrow{\frak L}^{\N})
$
the {\it dual} 
of
${\Cal T}_{{\K}^{\M}_{\N}} 
=(p,q:{\frak L}^{\K}\rightarrow{\frak L}^{\M})
$.
It is written as
${{\Cal T}_{{\K}^{\M}_{\N}}}^*$.
It is clear that
${({{\Cal T}_{{\K}^{\M}_{\N}}}^*)}^* = {\Cal T}_{{\K}^{\M}_{\N}}$.

\medskip

For 
${\Cal T}_{{\K}^{\M}_{\N}} 
=(p,q:{\frak L}^{\K}\rightarrow{\frak L}^{\M})
$
and 
$N\in \Bbb N$,
we will define  the $N$-heigher block 
${\Cal T}_{{\K}^{\M}_{\N}}^{[N]}$
of 
$\T$ 
as in the following way.
Let
$
{\frak L}^{\K^{[N]}} 
$
and 
$
{\frak L}^{\M^{[N]}}
$
be the $N$-heigher blocks of
$ 
{\frak L}^{\K}
$
and
$
{\frak L}^{\M}
$
respectively.
For $N \ge 2$, we will define the $\lambda$-graph system
$$
{\frak L}^{\N_{\Cal T}^{[N]}}
=(V^{\N_{\Cal T}^{[N]}}, E^{\N_{\Cal T}^{[N]}}, \lambda^{\N_{\Cal T}^{[N]}},
 \iota^{\N_{\Cal T}^{[N]}})
$$
over 
$
\Sigma^{\N_{\Cal T}^{[N]}} = 
\undersetbrace\text{ $N-1$ times}\to{\Sigma^\K \times\cdots\times \Sigma^\K}
$ 
by setting
$$
\align
V_l^{\N_{\Cal T}^{[N]}}
& = V_l^{\M^{[N]}},\\
E_{l,l+1}^{\N_{\Cal T}^{[N]}}
& = V_l^{\K^{[N]}} \\
&(= \{ (e_1,e_2,\dots,e_{N-1}) \in 
E_{l,l+1}^\K\times
E_{l+1,l+2}^\K\times\cdots\times E_{l+N-2,l+N-1}^\K \mid \\
& \quad t^\K(e_i) = s^\K(e_{i+1}), i=1,2,\dots,N-2\}),\\
\lambda^{\N_{\Cal T}^{[N]}}
& = \undersetbrace\text{ $N-1$ times}\to{\lambda^\K\times \lambda^\K\times\cdots\times\lambda^\K},\\
\iota^{\N_{\Cal T}^{[N]}}
& = \iota^{\M^{[N]}},
\endalign
$$
and
$$
\align
s^{\N_{\Cal T}^{[N]}}(e_1,e_2,\dots,e_{N-1})
& = p^E(e_1) p^E(e_2)\cdots p^E(e_{N-1}),\\
t^{\N_{\Cal T}^{[N]}}(e_1,e_2,\dots,e_{N-1})
& = q^E(e_1) q^E(e_2)\cdots q^E(e_{N-1}),\\
\endalign
$$
for
$(e_1,e_2,\dots,e_{N-1}) \in 
E_{l,l+1}^{\N_{\Cal T}^{[N]}} = V_l^{\K^{[N]}}.
$
The $\lambda$-graph system 
$
{\frak L}^{\N_{\Cal T}^{[N]}}
$
is called the $N$-{\it heigher block of } 
$
{\frak L}^{\N}
$
{\it relative to } $\T$.
For $N=1$, 
we put 
$
{\frak L}^{\N_{\Cal T}^{[1]}} = {\frak L}^\N.
$
The  $\lambda$-graph system homomorphism
$p^{[N]}:{\frak L}^{\K^{[N]} }\rightarrow {\frak L}^{\M^{[N]} }$
and
the one-shift $\lambda$-graph system homomorphism
$q^{[N]}:{\frak L}^{\K^{[N]} }\rightarrow {\frak L}^{\M^{[N]} }$
are defined by
$$
\align
p^{[N]}(e_1,e_2,\dots,e_N)&  =p^E(e_1)p^E(e_2)\cdots p^E(e_N),\\
q^{[N]}(e_1,e_2,\dots,e_N)& =q^E(e_1)q^E(e_2)\cdots q^E(e_N)
\endalign
$$
for
$
(e_1,e_2,\dots,e_N) \in E_{l,l+1}^\K\times E_{l+1,l+2}^\K\times\cdots
\times E_{l+N-1,l+N}^\K.
$
Then we have
\proclaim{Proposition 3.3}
The diagram
$$
\CD
                   @. {\frak L}^{\M^{[N]}} @. \\
@.                 @AA{p^{[N]}}A           @. \\
{\frak L}^{\N_{\Cal T}^{[N]}} @<s^{\K^{[N]}}<< {\frak L}^{\K^{[N]}} 
@>t^{\K^{[N]}}>> {\frak L}^{\N_{\Cal T}^{[N]}}\\
@.                 @VV{q^{[N]}}V           @. \\
                   @. {\frak L}^{\M^{[N]}} @. \\
\endCD
$$
defines a textile $\lambda$-graph system.
\endproclaim
We write the above textile $\lambda$-graph system 
$
{\Cal T}_{{\K^{[N]}}^{\M^{[N]}}_{\N_{\Cal T}^{[N]}}}
=(p^{[N]}, q^{[N]} :{\frak L}^{\K^{[N]} }\rightarrow{\frak L}^{\M^{[N]} })
$
as
${\Cal T}_{{\K}^{\M}_{\N}}^{[N]}
$
and call it the $N$-{\it heigher block }\, of $\T$.

\heading 4. Textile shifts on $\lambda$-graph systems
\endheading
For a $\lambda$-graph system
$\frak L = (V,E,\lambda,\iota)$,
we set 
$$
\align
X_{\frak L} 
& =\{ ({z_{l})}_{l=0}^{\infty} \in \prod_{l=0}^{\infty} E_{l,l+1} \mid
z_l \in E_{l,l+1}, t(z_{l}) = s(z_{l+1}), l=0,1,\dots \},\\
X_{{\frak L}_0} 
& =\{ ({z_{l})}_{l=1}^{\infty} \in \prod_{l=1}^{\infty} E_{l,l+1} \mid
z_l \in E_{l,l+1}, t(z_{l}) = s(z_{l+1}), l=1,2,\dots \}.
\endalign
$$
We define $S:X_{\frak L}  \longrightarrow X_{{\frak L}_0}$
by setting
$$
S(({z_{l})}_{l=0}^{\infty}) = ({z_{l})}_{l=1}^{\infty},
\qquad
({z_{l})}_{l=0}^{\infty}\in X_{\frak L}.
$$ 
For 
a textile $\lambda$-graph system
${\Cal T}_{{\K}^{\M}_{\N}} 
=(p,q:{\frak L}^{\K}\rightarrow{\frak L}^{\M}),
$
there exist maps
$p_X:
X_{{{\frak L}^{\K}}} \longrightarrow 
X_{{\frak L}^{\M}} 
$
and
$q_X:
X_{{{\frak L}^{\K}}} \longrightarrow 
X_{{\frak L}^{\M}_0} 
$
defined by
$
p_X(({z_{l})}_{l=0}^{\infty}) = 
(p^E({z_{l}))}_{l=0}^{\infty}
$
and
$
q_X(({z_{l})}_{l=0}^{\infty}) = 
(q^E({z_{l-1}))}_{l=1}^{\infty}
$
respectively.

Following Nasu's notation,
we say that 
a textile $\lambda$-graph system
$\T$ is {\it nondegenerate} if
both factor maps
$p_X : X_{{\frak L}^\K} \rightarrow X_{{\frak L}^\M}$
and
$q_X : X_{{\frak L}^\K} \rightarrow X_{{\frak L}^\M_0}$
are surjective.
We henceforth assume that textile $\lambda$-graph systems
$\T$ and $\T^*$ are both nondegenerate.

Let $\triangle$ be the lattice of the right-lower half plane:
$
{\triangle} 
 = \{ (i,j) \in {\Bbb Z}^2 \mid i+j \ge 0 \}.
$
A {\it textile edge weaved by } 
${\Cal T}_{{\K}^{\M}_{\N}}$
is a configuration
$$
{(e_{i,j})}_{(i,j) \in \triangle}
$$
such that 
\roster
\item
$e_{i,j}\in E_{i+j,i+j+1}^{\K}$ for $(i,j) \in \triangle$,
\item
${(e_{i,-i+l})}_{l\in \Zp} \in X_{{\frak L}^{\K}}
$
for each $i \in \Zp,$
\item
$p^E(e_{i,j}) =q^E(e_{i-1,j})$ for $i,j \in \Bbb Z$ with $ i+j\ge 1.$
\endroster
That is
a sequence 
$$
{(e_i)}_{i \in \Bbb Z}
$$
such that 
\roster
\item
$e_i = {(e_{i,-i+l})}_{l\in \Zp} \in X_{{\frak L}^{\K}}$ for $ i \in \Bbb Z,$
\item
$S \circ p_X(e_{i}) =q_X(e_{i-1})$ for $i \in \Bbb Z$.
\endroster
A textile edge weaved by $\T$ is regarded as a configuration of concatenated edges of ${\frak L}^\K$ on the lattice $\triangle$ 
of the right-lower half plane as in the following way.  
$$
\matrix
          &        &       &        &       &              &            &
                                            \cdots \qquad\\
          &        &       &        &       &              & e_{-l,l}   &
                                            \cdots \quad e_{-l}\\    
          &        &       &        &       &e_{-(l-1),l-1}&e_{-(l-1),l}&
                                     \quad       \cdots \quad e_{-(l-1)}\\
          &        &       &        & \cdots&\cdots        & \cdots     &
                                            \cdots \qquad \cdot \\ 
          &        &       &e_{-1,1}&e_{-1,2}&e_{-1,3}     &  e_{-1,4}  &
                                            \cdots \quad  e_{-1}\\ 
          &        &e_{0,0}&e_{0,1} &e_{0,2} &e_{0,3}      &  e_{0,4}   &
                                            \cdots \quad  e_0\\ 
          &e_{1,-1}&e_{1,0}&e_{1,1} &e_{1,2} &e_{1,3}      &  e_{1,4}   &
                                            \cdots \quad  e_1\\ 
e_{2,-2}  &e_{2,-1}&e_{2,0}&e_{2,1} &e_{2,2} &e_{2,3}      &  e_{2,4}   &
                                            \cdots \quad  e_2\\ 
\cdots&\cdots&\cdots&\cdots &\cdots &\cdots  & \cdots   & \cdots \quad 
\cdot \\ 
\endmatrix
$$
It is easy to see that 
$
{(e_{i,j})}_{(i,j) \in \triangle}
$
is a textile weaved by  
${\Cal T}_{{\K}^{\M}_{\N}}$
if and only if 
$
{(e_{j,i})}_{(i,j) \in \triangle}
$
is a textile weaved by  
${{\Cal T}_{{\K}^{\M}_{\N}}}^*$.

Consider the set 
$
X({\Cal T}_{{\K}^{\M}_{\N}})
$ 
of all textile edges 
weaved by 
$
{\Cal T}_{{\K}^{\M}_{\N}}
$
$$
X({\Cal T}_{{\K}^{\M}_{\N}})
=
\{
{(e_i)}_{i \in \Bbb Z}\in \prod_{i \in \Bbb Z}X_{{\frak L}^{\K}} \mid
S \circ p_X(e_i) = q_X(e_{i-1}),\, i \in \Bbb Z\}.
$$
\medskip

For $e \in E_{l,l+1}^\K, l \in \Zp$,
there exists a specified equivalence
between
$$
(\lambda^{\N}(s(e)),\lambda^{\N}(t(e)),
\lambda^{\M}(p^E(e)),\lambda^{\M}(q^E(e))) 
\in 
\Sigma^{\N}\times \Sigma^{\N}\times 
\Sigma^{\M}\times \Sigma^{\M}
$$ 
and
$
\lambda^{\K}(e) \in
 \Sigma^{\K}.
$
We may identify them, and assume that 
$$
\Sigma^\K = \{ 
\lambda^{\K}(e) \mid
 e \in E_{l,l+1}^\K, l \in \Zp \}.
 $$
We define
$$
\align
J_{upper}: & \Sigma^\K \rightarrow \Sigma^\M \quad
\text{ by } \quad 
J_{upper}(\lambda^\K(e))= \lambda^\M(p^E(e)),\\
J_{lower}: & \Sigma^\K \rightarrow \Sigma^\M \quad
\text{ by } \quad 
J_{lower}(\lambda^\K(e))= \lambda^\M(q^E(e)),\\
J_{right}: & \Sigma^\K \rightarrow \Sigma^\N \quad
\text{ by } \quad 
J_{right}(\lambda^\K(e))= \lambda^\N(t^\K(e)),\\
J_{left}: & \Sigma^\K \rightarrow \Sigma^\N \quad
\text{ by } \quad 
J_{left}(\lambda^\K(e))= \lambda^\N(s^\K(e)).\\
\endalign
$$
Let 
$
\Lambda_{\M},\Lambda_{\K},\Lambda_{\N} 
$ 
be the two-sided subshifts presented by 
${\frak L}^\M,{\frak L}^\K,{\frak L}^\N $  
respectively.
The above one-block maps 
$J_{left}, J_{right},J_{upper},J_{lower}$
give rise to sliding block codes between the subshifts:
$$
\align
\xi & = {(J_{upper})}_{\infty} : 
\Lambda_{\K} \rightarrow \Lambda_{\M},\\
\eta & = {(J_{lower})}_{\infty}  : 
\Lambda_{\K} \rightarrow \Lambda_{\M},\\
\xi^* & = {(J_{right})}_{\infty}  : 
\Lambda_{\K^*} \rightarrow \Lambda_{\N},\\
\eta^* & = {(J_{left})}_{\infty} : 
\Lambda_{\K^*} \rightarrow \Lambda_{\N}\\
\endalign
$$
respectively. 
We say that 
$\T$ is  1-1 if
 the factor codes
$\xi:\Lambda_{\K} \rightarrow 
\Lambda_{\M}
$
and
$\eta:\Lambda_{\K} \rightarrow 
\Lambda_{\M}
$
are both one-to-one.
Since 
$\T$ is nondegenerate,
the codes 
$\xi:\Lambda_{\K} \rightarrow 
\Lambda_{\M}
$
and
$\eta:\Lambda_{\K} \rightarrow 
\Lambda_{\M}
$
are both
surjective.
Hence, in this case,
we have an automorphism
$
\varphi_{\Cal T} = \eta \circ {\xi}^{-1}
$
on the subshift
$\Lambda_{\M}.
$

Let 
$
{\Cal U}_{{\frak L}^{\K}}
$
be the set of configurations of labels of 
${\frak L}^\K$ on
$X(\T):$
$$
{\Cal U}_{{\frak L}^{\K}}
= \{ {(\lambda^{\K}(e_{i,j}))}_{(i,j) \in \triangle} 
\in \prod \Sb (i,j) \in \triangle  \endSb \Sigma^\K \mid 
{(e_{i,j})}_{(i,j) \in \triangle} \in X({\T}) \}.
$$
Consider the natural product topology on
$\prod \Sb (i,j) \in \triangle  \endSb \Sigma^\K$
and restrict it to 
$
{\Cal U}_{{\frak L}^{\K}}
$
so that 
$
{\Cal U}_{{\frak L}^{\K}}
$ is compact.
For a connected subset $\Omega$ of ${\triangle}$,
we set 
$$
\align
X({\Cal T}_{{\K}^{\M}_{\N}};\Omega)
= \{ {(e_{i,j})}_{(i,j) \in \Omega} 
& \in  \prod_{(i,j) \in \Omega}E_{i+j,i+j+1}^\K \mid \\ 
 p^E(e_{i,j})& =  q^E(e_{i-1,j})\text{ for } (i,j),(i-1,j) \in \Omega,\\
 q^E(e_{i,j})& =  p^E(e_{i+1,j})\text{ for } (i,j),(i+1,j) \in \Omega,\\
 s^\K(e_{i,j})& =  t^\K(e_{i,j-1})\text{ for } (i,j),(i,j-1) \in \Omega,\\
 t^\K(e_{i,j})&=  s^\K(e_{i,j+1})\text{ for } (i,j),(i,j+1) \in \Omega \}
\endalign
$$
and
$$
{\Cal U}_{{\frak L}^\K}(\Omega) = 
\{ {\lambda^\K(e_{i,j})}_{(i,j)\in \Omega}\in 
\prod_{(i,j) \in \Omega} \Sigma^\K \mid
{(e_{i,j})}_{(i,j)\in \Omega}\in X({\Cal T}_{{\K}^{\M}_{\N}};\Omega)
\}.
$$
Hence
${\Cal U}_{{\frak L}^\K}({\triangle}) ={\Cal U}_{{\frak L}^\K}.
$
Put 
$$
{\triangle}_{(n,k)} 
 = \{ (i,j) \in {\triangle} \mid i \le n,\, j \le k \}
\qquad \text{ for }(n,k) \in  {\triangle}.
$$
By noticing the assumption
 that $\T$ and ${\T}^*$ are both nondegenerate, 
 we have 
\proclaim{Lemma 4.1}
For $(n,k) \in {\triangle}$
and
${(a_{i,j})}_{(i,j)\in \triangle_{(n,k)}} \in 
{\Cal U}_{{\frak L}^\K}(\triangle_{(n,k)}),
$
there exists
${(b_{i,j})}_{(i,j)\in \triangle_{(n+1,k+1)}} \in 
{\Cal U}_{{\frak L}^\K}(\triangle_{(n+1,k+1)})
$
such that
$$
b_{i,j} = a_{i,j}\qquad \text{ for all } (i,j) \in\triangle_{(n,k)}.
$$
\endproclaim
\demo{Proof}
For 
${(a_{i,j})}_{(i,j)\in \triangle_{(n,k)}} \in 
{\Cal U}_{{\frak L}^\K}(\triangle_{(n,k)}),
$
take 
$
{(e_{i,j})}_{(i,j)\in \triangle_{(n,k)}}
\in X({\Cal T}_{{\K}^{\M}_{\N}};\triangle_{(n,k)})
$  
such that
$
a_{i,j} = \lambda^\K(e_{i,j})$ for all $(i,j)\in \triangle_{(n,k)}.
$
Since ${\T}^*$ is nondegenerate,
there exist
$e_{i,k+1}\in E^\K_{i,k+1}$ for $i=-k-1,-k,\dots,n-1,n$ 
such that
$$
p^E(e_{i,k+1}) =  q^E(e_{i-1,k+1})
\quad
\text{ and }
\quad
 s^\K(e_{i,k+1}) =  t^\K(e_{i,k})
$$
for
$
 i = -k, -k+1,\dots, n.
$
Hence we have
$
{(e_{i,j})}_{(i,j)\in \triangle_{(n,k+1)}}
\in X({\Cal T}_{{\K}^{\M}_{\N}};\triangle_{(n,k+1)}).
$
Similarly by the condition that $\T$ is nondegenerate,
there exist
$e_{n+1,j}\in E^\K_{n+1,j}$ for 
$j=-n-1,-n,\dots,k,k+1$ 
such that
$$
p^E(e_{n+1,j}) =  q^E(e_{n,j})
\quad
\text{ and }
\quad
 s^\K(e_{n+1,j}) =  t^\K(e_{n+1,j-1})
$$
for
$
 j = -n, -n+1,\dots, k, k+1.
 $
This implies that
$
{(e_{i,j})}_{(i,j)\in \triangle_{(n+1,k+1)}}
\in X({\Cal T}_{{\K}^{\M}_{\N}};\triangle_{(n+1,k+1)})
$
so that by putting
$b_{i,j} = \lambda^\K(e_{i,j})$ for
$
{(e_{i,j})}_{(i,j)\in \triangle_{(n+1,k+1)}}
$
we get the assertion.
\qed
\enddemo
Hence one has
\proclaim{Corollary 4.2}
For $(n,k) \in {\triangle}$
and
${(a_{i,j})}_{(i,j)\in \triangle_{(n,k)}} \in 
{\Cal U}_{{\frak L}^\K}(\triangle_{(n,k)}),
$
there exists
${(\alpha_{i,j})}_{(i,j)\in \triangle} \in 
{\Cal U}_{{\frak L}^\K}
$
such that
$$
\alpha_{i,j} = a_{i,j}\qquad \text{ for all } (i,j) \in\triangle_{(n,k)}.
$$
\endproclaim

\proclaim{Proposition 4.3}
Let ${(\alpha_{i,j})}_{(i,j)\in \triangle}$ be symbols
$\alpha_{i,j} \in \Sigma^\K$ for 
$(i,j)\in \triangle$.
Then
${(\alpha_{i,j})}_{(i,j)\in \triangle} \in {\Cal U}_{{\frak L}^\K}$
if and only if
${(\alpha_{i,j})}_{(i,j)\in \triangle_{(n,k)}} \in 
{\Cal U}_{{\frak L}^\K}(\triangle_{(n,k)})
$
for all $(n,k) \in {\triangle}.$
\endproclaim
\demo{Proof}
The only if part is clear.
Suppose that
${(\alpha_{i,j})}_{(i,j)\in \triangle_{(n,k)}} \in 
{\Cal U}_{{\frak L}^\K}(\triangle_{(n,k)})
$
for all $(n,k) \in {\triangle}.$
By Corollary 4.2, for $(n,k) \in \triangle$
there exists
$\alpha^{(n,k)} ={(\alpha^{(n,k)}_{i,j})}_{(i,j)\in \triangle} \in 
{\Cal U}_{{\frak L}^\K}
$
such that
$$
\alpha^{(n,k)}_{i,j} = \alpha_{i,j} \qquad \text{ for all }  
(i,j) \in\triangle_{(n,k)}.
$$
Since ${\Cal U}_{{\frak L}^\K}$ is compact,
there exists 
$\bar{\alpha}={(\bar{\alpha}_{i,j} )}_{(i,j)\in \triangle}\in
{\Cal U}_{{\frak L}^{\K}}
$
such that
$\bar{\alpha} = \lim_{n,k\to{\infty}}\alpha^{(n,k)}.$
Since
$\bar{\alpha}_{i,j} = \alpha^{(n,k)}_{i,j} = \alpha_{i,j}
$ for all $(i,j) \in \triangle_{(n,k)}$,
one has
$\bar{\alpha}_{i,j} = \alpha_{i,j}
$
for all $(i,j) \in \triangle$
and hence
${(\alpha_{i,j})}_{(i,j)\in \triangle} \in {\Cal U}_{{\frak L}^\K}$.
\qed
\enddemo

\proclaim{Lemma 4.4}
For $\alpha = {(\alpha_{i,j})}_{(i,j)\in \triangle} 
\in 
{\Cal U}_{{\frak L}^{\K}},$
put
$$
{S_R(\alpha)}_{i,j} =  \alpha_{i,j+1},\quad
{S_D(\alpha)}_{i,j} =  \alpha_{i+1,j}\qquad \text{ for } (i,j) \in \triangle.
$$
Then we have
${S_R(\alpha)},{S_D(\alpha)}\in {\Cal U}_{{\frak L}^{\K}}$.
\endproclaim 
\demo{Proof}
For 
$
\alpha = {(\alpha_{i,j})}_{(i,j)\in \triangle} 
\in 
{\Cal U}_{{\frak L}^{\K}},
$
take 
${(e_{i,j})}_{(i,j)\in \triangle} 
\in 
X(\T)$
such that
$
\alpha_{i,j} = \lambda^\K(e_{i,j}),
$
where
$
e_{i,j} \in E_{i+j,i+j+1}^\K.
$
By the map
$\varphi_{i+j,i+j+1}^\K:
E_{i+j,i+j+1}^\K \rightarrow E_{i+j-1,i+j}^\K
$
in Lemma 3,1,
one has 
${\varphi^\K(e_{i,j+1})}_{(i,j)\in \triangle} 
\in 
X(\T).
$
Since
$
\lambda^\K(\varphi^\K(e_{i,j+1})) = 
\lambda^\K(e_{i,j+1}) = \alpha_{i,j+1},
$
one sees that
$S_R(\alpha)\in {\Cal U}_{{\frak L}^{\K}}$.
One may symmetrically prove that
${S_D(\alpha)}\in {\Cal U}_{{\frak L}^{\K}}$
by considering $\varphi^{\K^*}$.
\qed
\enddemo
The assertions above mean that 
${\Cal U}_{{\frak L}^{\K}}$ can be shifted to both left and upper.

We note that $S_R \circ S_D = S_D \circ S_R $ on ${\Cal U}_{{\frak L}^{\K}}$
to set 
$$
{\Cal U}^{\infty}_{{\frak L}^{\K}}
= \bigcap_{n,m=0}^{\infty} S_R^n\circ S_D^m({\Cal U}_{{\frak L}^{\K}}).
$$
Hence one has
$$
S_R({\Cal U}^{\infty}_{{\frak L}^{\K}}
) = {\Cal U}^{\infty}_{{\frak L}^{\K}}
=
S_D({\Cal U}^{\infty}_{{\frak L}^{\K}}
).
$$
A {\it textile label weaved by} $\Cal T=\T$ 
is a two-dimensional configuration 
${(\alpha_{i,j})}_{(i,j)\in {\Bbb Z}^2}$
of $\Sigma^{\K}$ such that
$$
{(\alpha_{i-k,j})}_{(i,j)\in \triangle}\in 
{\Cal U}^{\infty}_{{\frak L}^{\K}}
\qquad\text{ for all }\quad k \in \Bbb Z.
$$
The condition is equivalent to the condition 
$$
{(\alpha_{i,j-k})}_{(i,j)\in \triangle}\in 
{\Cal U}^{\infty}_{{\frak L}^{\K}}
\qquad\text{ for all }\quad k \in \Bbb Z.
$$
Let ${\Cal U}_{\Cal T}$ be the set of all textile label weaved by $\Cal T$. 
We note that

\proclaim{Lemma 4.5}
For 
$
{(\alpha_{i,j})}_{(i,j)\in {\Bbb Z}^2} \in {\Cal U}_{\Cal T}
$
one has
${(\alpha_{i,j})}_{j \in \Bbb Z} \in {\Lambda}_{\K}$
for all $i \in \Bbb Z,$
and
${(\alpha_{i,j})}_{i \in \Bbb Z} \in {\Lambda}_{\K^*}$
for all $j \in \Bbb Z.$
\endproclaim
\demo{Proof}
Fix an arbitrary integer $i \in \Bbb Z$.
For 
$
{(\alpha_{i,j})}_{(i,j)\in {\Bbb Z}^2} \in {\Cal U}_{\Cal T},
$
one sees that
$
{(\alpha_{i,j-k})}_{(i,j)\in \triangle} \in {\Cal U}_{{\frak L}^\K}
$
for $k \in \Zp.$
Hence
there exists
$
{(e_{i,j-k})}_{(i,j)\in \triangle} \in X(\T)
$
such that
$
\alpha_{i,j-k}= \lambda^\K(e_{i,j-k}),
$
so that
$
{(\alpha_{i,j-k})}_{j\in \Bbb Z, (i,j) \in \triangle}
\in \Lambda_{\K}^+
$
for all $k \in \Zp.$
We then have
$
{(\alpha_{i,j})}_{j \in \Bbb Z}
\in \Lambda_{\K}.
$
We similarly  have
$
{(\alpha_{i,j})}_{i \in \Bbb Z} \in \Lambda_{\K^*}
$
for all 
$j \in \Bbb Z.$
\qed
\enddemo

We define a metric $\delta_{\M}$ on 
$\Lambda_\M$ as follows:
for $\alpha = {(\alpha_{i})}_{i\in \Bbb Z},
\alpha' = {({\alpha'}_{i})}_{i\in \Bbb Z}$
in $\Lambda_\M$
$$
\delta(\alpha,\alpha')
=
\cases
0 & \text{ if } \quad \alpha =\alpha', \\
\frac{1}{k+1} & \text{ if }\quad  \alpha \ne\alpha'
\endcases
$$
where
$ k = \min\{|i|  \mid i \in \Bbb Z, \alpha_{i} \ne  {\alpha'}_{i}\}.$
Similarly we define a metric $\delta_{\N}$
on 
${\Lambda}_{\N}.$
We next define a metric $\delta_{\Cal T}$ on 
${\Cal U}_{\Cal T}$ as follows:
for $u = {(\alpha_{i,j})}_{(i,j)\in {\Bbb Z}^2},
u' = {({\alpha'}_{i,j})}_{(i,j)\in {\Bbb Z}^2}$
in ${\Cal U}_{\Cal T}$
$$
\delta_{\Cal T}(u,u')
=
\cases
0 & \text{ if } u=u', \\
\frac{1}{k+1} & \text{ if } u \ne u'
\endcases
$$
where
$
 k = \min\{|i| + |j| \mid i,j \in \Bbb Z, \alpha_{i,j} \ne  {\alpha'}_{i,j}\}.
$

\proclaim{Lemma 4.6}
${\Cal U}_{\Cal T}$ is compact. 
\endproclaim
\demo{Proof}
We first note that the set $X(\T)$ of all textile edges is a compact set in a natural topology of the edge set so that the label set 
${\Cal U}_{{\frak L}^\K}$ and 
${\Cal U}^{\infty}_{{\frak L}^\K}$
are both compact.
Let $\prod\Sb (i,j) \in {\Bbb Z}^2 \endSb\Sigma^\K$ be the set 
$\alpha_{i,j} \in \Sigma^\K, \, (i,j) \in {\Bbb Z}^2$
of all two-dimensional configuration of $\Sigma^\K$,
that is endowed with the topology similarly defined by the above 
$\delta_{\Cal T}$.
Consider the sequence of the following continuous maps
$$
\zeta_k : {(\alpha_{i,j})}_{(i,j) \in {\Bbb Z}^2}
\in \prod \Sb (i,j) \in {\Bbb Z}^2 \endSb \Sigma^\K 
\longrightarrow
{(\alpha_{i-k,j})}_{(i,j) \in \triangle}\in 
\prod \Sb (i,j) \in \triangle \endSb \Sigma^\K, \qquad k \in \Bbb Z.
$$ 
Since we have
$${\Cal U}_{\Cal T} = \bigcap_{k \in \Bbb Z}\zeta_k^{-1}
({\Cal U}_{{\frak L}^{\K}}),
$$
the set ${\Cal U}_{\Cal T}$ is compact.
\qed
\enddemo

Define a one-block code
$$
\Phi_{\Cal T}:{\Cal U}_{\Cal T}\rightarrow \Lambda_{\K}
$$
by setting 
$$
\Phi_{\Cal T}({(\alpha_{i,j})}_{(i,j)\in {\Bbb Z}^2}) = 
{(\alpha_{0,j})}_{j\in \Bbb Z},\qquad {(\alpha_{i,j})}_{(i,j)\in {\Bbb Z}^2}
\in
{\Cal U}_{\Cal T}.
$$
We say that the textile $\lambda$-graph system $\T$ is {\it surjective} \, if
 the map
 $
\Phi_{\Cal T} :
{\Cal U}_{\Cal T} \rightarrow \Lambda_{\K}
$ 
is surjective.
Define the one-block codes
$$
\Theta_{\Cal T} :
{\Cal U}_{\Cal T} \rightarrow \Lambda_{\M},
\qquad
\Theta_{\Cal T}^* :
{\Cal U}_{\Cal T} \rightarrow \Lambda_{\N}
$$
by setting
$$
 \Theta_{\Cal T}({(\alpha_{i,j})}_{(i,j)\in {\Bbb Z}^2}) 
 = ({J_{lower}(\alpha_{0,j}))}_{j \in \Bbb Z},
\qquad 
 \Theta_{\Cal T}^*({(\alpha_{i,j})}_{(i,j)\in {\Bbb Z}^2}) 
 = ({J_{right}(\alpha_{i,0}))}_{i \in \Bbb Z}.
$$
They are continuous in the topology defined by the metric $\delta_{\K}$ on 
${\Cal U}_{\Cal T}$.
Since $\eta:\Lambda_{\K}\rightarrow\Lambda_{\M}$
is always surjective
and
$\Theta_{\Cal T} = \eta\circ \Phi_{\Cal T}$,
if $\T$ is surjective, the map $\Theta_{\Cal T}$ is surjective.
For $k,n \in \Bbb Z,$
the homeomorphism
$$
\sigma^{(k,n)}_{\Cal T}: {\Cal U}_{\Cal T}\longrightarrow
{\Cal U}_{\Cal T}
$$
is
defined by 
$$
\sigma^{(k,n)}_{\Cal T}({(\alpha_{i,j})}_{(i,j)\in {\Bbb Z}^2})
=({\alpha_{i+k,j+n})}_{(i,j)\in {\Bbb Z}^2}
\qquad \text{ for }
\quad
{(\alpha_{i,j})}_{(i,j)\in {\Bbb Z}^2} \in {\Cal U}_{\Cal T}.
$$
The dynamical system
$$
({\Cal U}_{\Cal T},
\sigma^{(k,n)}_{\Cal T})$$ 
is called the $(k,n)$-textile shift on $\lambda$-graph systems.
\proclaim{Lemma 4.7}
\roster
\item"(i)"
If  $\T$ is 1-1 and surjective,
then 
$
\Theta_{\Cal T} :
{\Cal U}_{\Cal T} \rightarrow \Lambda_{\M}
$
is a homeomorphism such that 
$
\Theta_{\Cal T}\circ \sigma^{(0,n)}_{\Cal T}=
\sigma_{\M}^{n}\circ  \Theta_{\Cal T}.$
\item"(ii)"
If ${\T}^*$ is 1-1 and surjective, 
then 
$
\Theta_{\Cal T}^* :
{\Cal U}_{\Cal T} \rightarrow \Lambda_{\N}
$
is a homeomorphism such that
$
\Theta_{\Cal T}^*\circ \sigma^{(k,0)}_{\Cal T}=
\sigma^{k}_{\N}\circ  \Theta_{\Cal T}^*.$
\endroster
\endproclaim
\demo{Proof}
(i)
If $\T$ is 1-1 and surjective, 
then 
$
\Theta_{\Cal T} :
{\Cal U}_{\Cal T} \rightarrow \Lambda_{\M}
$
is one-to-one and surjective so that it is a homeomorphism.
(ii)
The assertion is symmetric to (i).
\qed
\enddemo
\proclaim{Proposition 4.8}
If $\T$ is $1-1$ and surjective, 
then
$(\Lambda_{\M},\varphi_{\Cal T}^k \sigma_{\M}^n)$
is conjugate to
$({\Cal U}_{\Cal T},\sigma_{\Cal T}^{(k,n)})$
for all
$k,n \in \Bbb Z.$
\endproclaim
\demo{Proof}
We note that 
$
\varphi_{\Cal T}^k(({J_{lower}(\alpha_{0,j}))}_{j\in \Bbb Z})
=({J_{lower}(\alpha_{k,j}))}_{j\in \Bbb Z})
$ 
for ${(\alpha_{k,j})}_{(k,j) \in {\Bbb Z}^2}\in {\Cal U}_{\Cal T}$.
Since 
$\T$ is $1-1$ and surjective,
the map
$
\Theta_{\Cal T} :
{\Cal U}_{\Cal T} \rightarrow \Lambda_{\M}
$
is a homeomorphism that 
gives rise to a conjugacy between
$(\Lambda_{\M},\varphi_{\Cal T}^k \sigma_{\M}^n)$
 and 
$({\Cal U}_{\Cal T},\sigma_{\Cal T}^{(k,n)})$.
\qed
\enddemo

Now we reach the following theorem.
\proclaim{Theorem 4.9}
Suppose that $\T$ and $\T^*$ are both $1-1$ and surjective.
Then there exists
a homeomorphism
$$\chi_{\Cal T}:
\Lambda_{\M}\rightarrow \Lambda_{\N}
$$ such that 
the diagrams 
$$
\CD
\Lambda_{\M}  @>\varphi_{\Cal T}>> 
                                          \Lambda_{\M} \\
@VV{\chi_{\Cal T}}V   @VV{\chi_{\Cal T}}V \\   
\Lambda_{\N}  @>\sigma_{\N}>> 
\Lambda_{\N}, \\
\endCD
\qquad
\qquad
\CD
\Lambda_{\M}  @>\sigma_{\M}>> 
                                          \Lambda_{\M} \\
@VV{\chi_{\Cal T}}V   @VV{\chi_{\Cal T}}V \\   
\Lambda_{\N}  @>\varphi_{{\Cal T}^*}>> 
\Lambda_{\N} \\
\endCD
$$
are both commutative.
\endproclaim
\demo{Proof}
We set 
$$
\chi_{\Cal T} = \Theta_{\Cal T}^*\circ {\Theta_{\Cal T}}^{-1}:
\Lambda_{\M}\rightarrow \Lambda_{\N}.
$$
It satisfies
$$
\chi_{\Cal T}\circ \varphi_{\Cal T} 
= \sigma_{\N}\circ \chi_{\Cal T},\qquad
\chi_{\Cal T}\circ \sigma_{\M} 
= \varphi_{{\Cal T}^*}\circ \chi_{\Cal T}.
$$
\qed
\enddemo

There are various metrics on $\Lambda_\M$ by which the product topology on $\Lambda_\M$ 
is given. 
Any such metric makes the homeomorphism $\sigma_\M$ on $\Lambda_\M$   
expansive.
We may fix the previously defined metric on $\Lambda_\M$. 
By the metric, 
$\sigma_\M$ has $1$ as its expansive constant. 
Theorem 4.9 is generalized as folows:
\proclaim{Theorem 4.10(cf.[N2;Theorem 4.1])}
Assume that 
$\T$ is 1-1 and surjective. 
\roster
\item"(i)"
If  
$
\varphi_{\Cal T}
$
 is expansive and
its expansive constant $c$ satisfies
$c \ge \frac{1}{k}$ for some $k\in \Bbb N$,  
then
${{{\Cal T}_{{\K}^{\M}_{\N}}}^{[2k]}}^*$ is 1-1.
Hence if there is no $n\in \Bbb N$ such that 
${{{\Cal T}_{{\K}^{\M}_{\N}}}^{[n]}}^*$ is 1-1,
then 
$
\varphi_{\Cal T}: \Lambda_{\M} 
\rightarrow \Lambda_{\M}
$
is not expansive.
\item"(ii)"
If there is  $n\in \Bbb N$ such that 
${({{\Cal T}_{{\K}^{\M}_{\N}}}^{[n]})}^*$ is 1-1 and surjective,
then one has topological conjugacies:
$$
(\Lambda_{\M}, \varphi_{\Cal T})
\simeq 
(\Lambda_{{\N_{\Cal T}^{[n]}}}, \sigma_{{\N_{\Cal T}^{[n]}}}),
\qquad
(\Lambda_{{\N_{\Cal T}^{[n]}}}, \varphi_{{{\Cal T}^{[n]}}^*})
\simeq 
(\Lambda_{\M}, \sigma_{\M}),
$$
where ${\N_{\Cal T}^{[n]}}$ is the $n$-heigher block of ${\frak L}^\N$
relative to $\T$.
Hence 
the topological dynamical system
$
(\Lambda_{\M}, \varphi_{\Cal T})
$
is realized as the subshift
$ 
(\Lambda_{{\N_{\Cal T}^{[n]}}}, \sigma_{{{\N_{\Cal T}^{[n]}}}}).
$
 If in particular 
$
\varphi_{\Cal T}
$
 is expansive and
its expansive constant $c$ satisfies
$c \ge \frac{1}{k}$ for some $k\in \Bbb N$
and
${{{\Cal T}_{{\K}^{\M}_{\N}}}^{[2k]}}^*$ 
is surjective,
then
the topological dynamical system
$
(\Lambda_{\M}, \varphi_{\Cal T})
$
is topologically conjugate to the subshift
$ 
(\Lambda_{{\N_{\Cal T}^{[2k]}}}, \sigma_{{{\N_{\Cal T}^{[2k]}}}})
$
presented by the $\lambda$-graph system 
$
{\frak L}^{{\N_{\Cal T}^{[2k]}}}
$.
\endroster
\endproclaim
\demo{Proof}
The proofs bellow are essentially similar to the proofs of [N2;Theorem 4.1].
We will give the proofs for the sake of completeness.
(i)
Assume that
$
\varphi_{\Cal T}
$
 is expansive and
its expansive constant $c$ satisfies
$c \ge \frac{1}{k}$ for some $k\in \Bbb N$.
Suppose that 
${{{\Cal T}_{{\K}^{\M}_{\N}}}^{[2k]}}^*$ is not 1-1.
There are distinct textile labels
$ s = {(\beta_{i,j})}_{(i,j) \in {\Bbb Z}^2}$
and 
$ s' = {(\beta'_{i,j})}_{(i,j) \in {\Bbb Z}^2}$
in 
${\Cal U}_{\Cal T}$
such that 
$
\beta_{i,j} = \beta'_{i,j}
$
for $i \in \Bbb Z, -(k-1) \le j \le k-1$.
Now $\T$ is $1-1,$
by putting
$y ={(y_j)}_{j \in \Bbb Z}= \Theta_{\Cal T}(s),
y'={(y'_j)}_{j \in \Bbb Z} = \Theta_{\Cal T}(s'),
$
we have
$ y \ne y' \in \Lambda_{\M}.
$
Since
one has
$y_j = y'_j$ for 
$-(k-1)\le j \le k-1$,
one sees that
${\varphi_{\Cal T}^i(y)}_j = {\varphi_{\Cal T}^i({y'})}_j$ 
for $ i\in \Bbb Z$ and 
$-(k-1)\le j \le k-1$.
Hence we have
$$
d(\varphi_{\Cal T}^i(y),\varphi_{\Cal T}^i(y')) <\frac{1}{k}
$$ 
for all $ i\in \Bbb Z,$ a contradiction. 

(ii)
Since
$
(\Lambda_{\M^{[n]}},\varphi_{{\Cal T}^{[n]}})
$
is topologically conjugate to 
$ 
(\Lambda_{\M},\varphi_{\Cal T})
$ 
and
$
(\Lambda_{\M^{[n]}},\sigma_{\M^{[n]}})
$
is topologically conjugate to
$ 
(\Lambda_{\M},\sigma_{\M}),
$
The assetion holds from Theorem 4.9.
\qed
\enddemo

\medskip

Following Nasu's consideration as in [N2;Section 2], 
we will define bias shifts on textile $\lambda$-graph systems.
For a symbolic matrix system $(\M,I)$,
we set for $k\in \Bbb N$
$$
\M_{l,l+k}  = \M_{l,l+1}\cdots \M_{l,l+k},
\qquad
I_{l,l+k}  = I_{l,l+1}\cdots I_{l,l+k}, \qquad l \in \Zp.
$$
Let
$(\M,I^{\M})$ and 
$(\N,I^{\N})$ 
form squares.
Then for $k,n\in \Zp$,
we set
$$
\align
(\N^k\M^n)_{l,l+1} & = \N_{l(k+n),(l+1)k + ln}
\M_{(l+1)k+ ln,(l+1)(k + n)},\\
(I^{\N^k\M^n})_{l,l+1} & = I^{\N}_{l(k+n),(l+1)k + ln}
I^{\M}_{(l+1)k+ ln,(l+1)(k + n)}, \qquad l \in \Zp.\\
\endalign
$$
As
$I^\N_{l,l+1} = I^\M_{l,l+1}$,
one sees that 
$(\N^k\M^n, I^{\N^k\M^n})$
becomes a symbolic matrix system over ${(\Sigma^\N)}^k{(\Sigma^\M)}^n.$
Similarly we have a symbolic matrix system 
$(\M^n\N^k, I^{\M^n\N^k})$
over ${(\Sigma^\M)}^n{(\Sigma^\N)}^k.$
For 
$\alpha = {(\alpha_{i,j})}_{(i,j)\in {\Bbb Z}^2}\in {\Cal U}_{\Cal T},$
we set
$$
\align
{(\check{c}^{(k,n)}(\alpha)_{(ik,in)})}
&=b_{ik,in }b_{ik+1,in}\cdots b_{ik+k-1,in}
a_{ik+k,in}a_{ik+k,in+1}\cdots a_{ik+k,in+n-1}\\
{(\hat{c}^{(k,n)}(\alpha)_{(ik,in)})}
&= a_{ik,in}a_{ik,in+1}\cdots a_{ik,in+n-1}
b_{ik,ik+in }b_{ik+1,ik+in}\cdots b_{ik+k-1,ik+n}
\endalign
$$
where
$
a_{i,j} = J_{upper}(\alpha_{i,j}),
b_{i,j} = J_{left}(\alpha_{i,j}).
$
Define 
$$
\check{\Theta}_{\Cal T}^{(k,n)}:{\Cal U}_{\Cal T} \rightarrow 
\Lambda_{\N^k\M^n},\qquad
\hat{\Theta}_{\Cal T}^{(k,n)}:{\Cal U}_{\Cal T} \rightarrow 
\Lambda_{\M^n\N^k}
$$
by setting
$$
\align
\check{\Theta}_{\Cal T}^{(k,n)}(\alpha )
&= {(\check{c}^{(k,n)}(\alpha)_{(ik,in)})}_{i\in \Bbb Z} 
\in \Lambda_{\N^k\M^n}\\
\hat{\Theta}_{\Cal T}^{(k,n)}(\alpha )
&= {(\hat{c}^{(k,n)}(\alpha)_{(ik,in)})}_{i\in \Bbb Z} \in \Lambda_{\N^k\M^n}.
\endalign
$$
We set 
$$
\align
{\check{\Cal U}}_{\Cal T}^{(k,n)}
 = \check{\Theta}_{\Cal T}^{(k,n)}({\Cal U}_{\Cal T}),&
\qquad
{\hat{\Cal U}}_{\Cal T}^{(k,n)}
 = \hat{\Theta}_{\Cal T}^{(k,n)}({\Cal U}_{\Cal T}),\\
\check{\sigma}_{\Cal T}^{(k,n)}
({(\check{c}^{(k,n)}(\alpha)_{(ik,in)})}_{i\in \Bbb Z})
&=
{(\check{c}^{(k,n)}(\alpha)_{((i+1)k,(i+1)n)})}_{i\in \Bbb Z},\\
\hat{\sigma}_{\Cal T}^{(k,n)}
({(\hat{c}^{(k,n)}(\alpha)_{(ik,in)})}_{i\in \Bbb Z})
&=
{(\hat{c}^{(k,n)}(\alpha)_{((i+1)k,(i+1)n)})}_{i\in \Bbb Z}.
\endalign
$$
We have subshifts
$$
({\check{\Cal U}}_{\Cal T}^{(k,n)},\check{\sigma}_{\Cal T}^{(k,n)})
\quad
\text{ and }
\quad
({\hat{\Cal U}}_{\Cal T}^{(k,n)},\hat{\sigma}_{\Cal T}^{(k,n)})
$$
over 
${(\Sigma^\N)}^k{(\Sigma^\M)}^n$
and
over ${(\Sigma^\M)}^n{(\Sigma^\N)}^k$
respectively.
\proclaim{Lemma 4.11}
$
({\check{\Cal U}}_{\Cal T}^{(k,n)},\check{\sigma}_{\Cal T}^{(k,n)})
$
is topologically conjugate to
$({\hat{\Cal U}}_{\Cal T}^{(k,n)},\hat{\sigma}_{\Cal T}^{(k,n)}).
$
\endproclaim
\demo{Proof}
Define 
$
\psi:{\hat{\Cal U}}_{\Cal T}^{(k,n)}\rightarrow
{\check{\Cal U}}_{\Cal T}^{(k,n)}
$
by setting
$$
\psi(\hat{\Theta}_{\Cal T}^{(k,n)}(\alpha) ) 
= \check{\Theta}_{\Cal T}^{(k,n)}(\sigma^{(0,n)}_{\Cal T}(\alpha))
$$
for
$\alpha \in {\Cal U}_{\Cal T}.$
It is direct to see that 
 $\psi$ is a topological conjugacy between
$
({\check{\Cal U}}_{\Cal T}^{(k,n)},\check{\sigma}_{\Cal T}^{(k,n)})
$
and
$({\hat{\Cal U}}_{\Cal T}^{(k,n)},\hat{\sigma}_{\Cal T}^{(k,n)}).
$
\qed
\enddemo

We  call
the subshift
$
({\check{\Cal U}}_{\Cal T}^{(k,n)},\check{\sigma}_{\Cal T}^{(k,n)})
$
the $(k,n)$-{\it bias shift} \, defind by $\T.$

\heading 5. LR textile systems on $\lambda$-graph systems
\endheading
In this section, we formulate LR textile $\lambda$-graph systems, that are generalization of sofic LR textile systems defined in Nasu [N2].

\proclaim{Proposition 5.1}
Assume that  $\lambda$-graph systems 
${\frak L}^{\M}$ and 
${\frak L}^{\N}$
form squares.
If there exists a specification 
$\kappa$ between 
$\Sigma^{\M}\Sigma^{\N}$
and
$\Sigma^{\N}\Sigma^{\M}$
 that gives specified equivalences 
$$
\M_{l,l+1}\N_{l+1,l+2} \overset{\kappa}\to{\simeq} \N_{l,l+1}\M_{l+1,l+2},
\qquad l \in \Zp, \tag 5.1
$$
then there exists a 
$\lambda$-graph system
${\frak L}^{\K}$ 
and 
a textile $\lambda$-graph system
${\Cal T}_{{\K}^{\M}_{\N}} 
=(p,q:{\frak L}^{\K}\rightarrow{\frak L}^{\M}).
$
\endproclaim
\demo{Proof}
We identify the vertex sets
$
V_l^{\M}
$
and
$ 
V_l^{\N}
$ 
for 
$l \in \Zp$.
Put
$$
\align
V_l^\K
  = & E^\N_{l,l+1},\\
E^\K_{l,l+1} 
= & \{
(f',f, e, e')\in E_{l,l+1}^\N \times E_{l+1,l+2}^\N \times 
E_{l,l+1}^\M \times E_{l+1,l+2}^\M \mid \\
& s^\M(e) = s^\N(f'), t^\M(e) = s^\N(f),
 t^\N(f') = s^\M(e'), \\
& t^\M(e') = t^\N(f),
\kappa(\lambda^\M(e)\lambda^\N(f)) = 
\lambda^\N(f')\lambda^\M(e') \}, \qquad l \in \Zp
\endalign
$$
and 
$$
V^\K = \cup_{l\in \Zp}V^\K_l,\qquad
E^\K = \cup_{l\in \Zp}E^\K_{l,l+1}.
$$
Each element $(f',f, e, e')\in E^\K_{l,l+1}$ is visualized as a square:
$$
\CD
  \cdot  @>e>>  \cdot\\
@VV{f'}V       @VV{f\quad.}V \\   
\cdot @>e'>> \cdot\\
\endCD
$$
We define
$
s^\K:E_{l,l+1}^\K \rightarrow V_l^\K,
t^\K:E_{l,l+1}^\K \rightarrow V_{l+1}^\K
$
by setting
$$
s^\K(f',f, e, e') = f',\qquad 
t^\K(f',f, e, e') = f
\qquad
\text{ for } 
(f',f, e, e') \in 
E_{l,l+1}^\K,
$$
and
$
\iota_{l,l+1}^\K:V_{l+1}^\K \rightarrow V_{l}^\K
$
by setting
$
\iota_{l,l+1}^\K = \varphi_{l+1}^\N.
$
We put
$$
\align
\Sigma^\K = 
\{
(\lambda^\N(f'), \lambda^\N(f), 
\lambda^\M(e), \lambda^\M(e'))\in \Sigma^\N \times \Sigma^\N \times 
\Sigma^\M \times \Sigma^\M & \mid \\ 
(f', f, e, e' ) \in E_{l,l+1}^\K, l\in \Zp & \}
\endalign
$$
and
$$
\lambda^\K: E^\K 
\ni (f',f, e, e') \longrightarrow \Sigma^\K
\ni (\lambda^\N(f'), \lambda^\N(f), 
\lambda^\M(e), \lambda^\M(e')).
$$
Then we will show that 
${\frak L}^\K =
(V^\K, E^\K, \lambda^\K, \iota^\K)$
is a $\lambda$-graph system over 
$\Sigma^\K$.
For 
$u \in V_{l-1}^{\K}, v \in V_{l+1}^{\K},$
put
$
w=\iota_{l,l+1}^{\K}(v). 
$
 One sees
$$
E^{\K^{l-1,l}}_{\iota^{\K}}(u,v) 
 = \{ (u, w, e, e') \in E^{\K}_{l-1,l} \mid 
e \in E^\M_{l-1,l}, e'\in E^\M_{l,l+1} \}.
$$
As $w \in V_l^\K=E_{l,l+1}^\N$ is fixed, 
if we choose $e \in E^\M_{l-1,l}$
such 
that 
$(u, w, e, e') \in
E^{\K^{l-1,l}}_{\iota^{\K}}(u,v)
$
for some $e'$,
the label
$
\lambda^\M(e)\lambda^\N(w) \in \Sigma^\M\Sigma^\N
$
determines the labels
$\lambda^\N(u)$ and $\lambda^\M(e')$ 
of 
$u$ and $e'$ through the specification $\kappa.$
Since the label
$ \lambda^\M(e')$
and the terminal $t^\M(e') = t^\N(w)$ are determined,
the edge $e'$ is uniquely determined 
because ${\frak L}^\M$ is left-resolving.
Hence under fixing 
$u \in V_{l-1}^{\K}$ and $v \in V_{l+1}^{\K},$
the edge $e' \in E_{l,l+1}^\M$ is  
uniquely determined 
by 
$e \in E_{l-1,l}^\M$,
so that 
$
E^{\K^{l-1,l}}_{\iota^{\K}}(u,v) 
$
is identified with 
$
\{  e \in E^{\M}_{l-1,l} \mid 
s^\M(e) = s^\N(u), t^\M(e) = s^\N(w) \}.
$
Now $\iota^\N = \iota^\M$ 
so that 
$
s^\N(w) = \iota_{l,l+1}^\N(s^\N(v)) = 
\iota_{l,l+1}^\M(s^\N(v)).
$
Hence 
$
E^{\K^{l-1,l}}_{\iota^{\K}}(u,v)
$ 
is identified with
$
{E^\M}^{l-1,l}_{\iota^{\M}}(s^\N(u), s^\N(v)).
$
On the other hand, one sees 
$$
E^{\K^{\iota^{\K}}}_{l,l+1}(u,v) 
= \{ (w', v, g, g') \in E^{\K}_{l,l+1} \mid 
\iota^{\K}_{l-1,l}(w')= u \}.
$$
Similarly to the discussion of 
$
E^{\K^{l-1,l}}_{\iota^{\K}}(u,v)
$,
 if we choose $g \in E^\M_{l,l+1}$
such that 
$
(w', v, g, g') \in E^{\K^{\iota^{\K}}}_{l,l+l}(u,v)
$
for some $g'$,
the label
$\lambda^\M(g)\lambda^\N(v) \in\Sigma^\M\Sigma^\N
$
determines  the labels 
$\lambda^\N(w')$ and $\lambda^\M(g')$
of 
$w'$ and $g'$ through the specification $\kappa$.
Since the label
$ \lambda^\M(g')$
and the terminal $t^\M(g') = t^\N(v)$ are determined,
the edge $g'$ is uniquely determined 
because ${\frak L}^\M$ is left-resolving,
so that the source vertex $s^\M(g') \in V_{l+1}^\M$
of $g'$ is dertermined.
Since $t^\N(w') = s^\M(g')$, 
the edge $w' \in E^\N_{l,l+1}$ is uniquely determined.
Hence under fixing the vertices 
$u \in V_{l-1}^{\K}, v \in V_{l+1}^{\K},$
both the edges 
$g' \in E^\M_{l+1,l+2}$ and $w' \in E_{l,l+1}^\N$ 
are  
uniquely determined 
by 
$g \in E_{l,l+1}^\M$.  
Now one has
$s^\M(g) = s^\N(w'), \varphi_{l+1}^\N(w') =u$
so that
$$
\iota^\M_{l-1,l}(s^\N(w') = \iota^\N_{l-1,l}(s^\N(w')) = s^\N(u) \in V_{l-1}^\M = V_{l-1}^\N.
$$ 
It then follows that 
$
E^{\K^{\iota^{\K}}}_{l,l+1}(u,v) 
$
is identified with 
$ \{  g \in E^{\M}_{l,l+1} \mid 
\iota^\M_{l-1,l}(s^\M(g))  = s^\N(u), t^\M(g) = s^\N(v) \},
$ 
that is
$
E^{\M^{\iota^{\M}}}_{l,l+1}(s^\N(u), s^\N(v)).
$
By the local property of ${\frak L}^\M$,
one has a label preserving bijection between 
$$
E^{\M^{\iota^{\M}}}_{l,l+1}(s^\N(u), s^\N(v))
\quad
\text{ and }
\quad
E^{\M^{l-1,l}}_{\iota^{\M}}(s^\N(u), s^\N(v))
$$
that yields a
a label preserving bijection between 
$$
E^{\K^{\iota^{\K}}}_{l,l+1}(u,v)
\quad
\text{ and }
\quad
E^{\K^{l-1,l}}_{\iota^{\K}}(u,v).
$$
This means 
${\frak L}^\K =
(V^\K, E^\K, \lambda^\K, \iota^\K)$
is a $\lambda$-graph system over 
$\Sigma^\K$.

Define 
a $\lambda$-graph system homomorphism 
$p: {\frak L}^\K \longrightarrow {\frak L}^\M$
by setting
$$
\gather
p^E:(f',f, e, e' ) \in E^\K_{l,l+1}  \longrightarrow e \in E^\M_{l,l+1},\\
p^V:f \in V^\K_{l} = E^\N_{l,l+1}  \longrightarrow s^\N(f) 
\in V^\N_l, \\
p^\Sigma:(\lambda^\N(f'), \lambda^\N(f), \lambda^\M(e), \lambda^\M(e') ) 
\in \Sigma^\K  \longrightarrow \lambda^\M(e) \in \Sigma^\M.
\endgather
$$
Define a one-shift $\lambda$-graph system homomorphism 
$q: {\frak L}^\K \longrightarrow {\frak L}^\M$
by setting
$$
\gather
q^E:(f',f, e, e' ) \in E^\K_{l,l+1}  \longrightarrow e' \in E^\M_{l+1,l+2},\\
q^V:f' \in V^\K_{l} = E^\N_{l,l+1}  \longrightarrow t^\N(f') 
\in  V^\N_{l+1}, \\
q^\Sigma:(\lambda^\N(f'), \lambda^\N(f), \lambda^\M(e), \lambda^\M(e') ) 
\in \Sigma^\K  \longrightarrow \lambda^\M(e') \in \Sigma^\M.
\endgather
$$
Since for 
$\alpha = (f',f, e, e' ) \in E^\K_{l,l+1}$,
one has
$$
(s^\K(\alpha),t^\K(\alpha), p^E(\alpha),q^E(\alpha)) =(f',f, e, e' )
$$
the square  
$
(s^\K(\alpha),t^\K(\alpha), p^E(\alpha),q^E(\alpha)) 
$ determines
$\alpha$,
and
quadruple
$$
(\lambda^\N(s^\K(\alpha)), \lambda^\N(t^\K(\alpha)), 
\lambda^\M(p^E(\alpha)), \lambda^\M(q^E(\alpha))) 
$$ 
determines
$\lambda^\K(\alpha)$.
Hence one has 
a textile  $\lambda$-graph system
${\Cal T}_{{\K}^{\M}_{\N}} 
=(p,q:{\frak L}^{\K}\rightarrow{\frak L}^{\M})
$
through the specified equivalence (5.1).
\qed
\enddemo

We call this textile $\lambda$-graph system 
an LR textile  $\lambda$-graph system, following Nasu's terminology
for sofic textile systems ([N2]).
\proclaim{Lemma 5.2(cf.[N2:Fact.6.14])}
An LR textile $\lambda$-graph system $\T$
is nondegenerate.
\endproclaim
\demo{Proof}
Let $\T$ be an LR textile $\lambda$-graph system defined by (5.1).
Keep the notation as in the previous proposition.
We will prove that 
$p_X :X_{{\frak L}^\K} \longrightarrow
X_{{\frak L}^\M}
$
is surjective.
We set for $l\in \Zp, n \in \Bbb N$
$$
\align
E^\M_{l,l+n}=
& \{ (e_1,\dots,e_n) \in 
E^\M_{l,l+1}\times E^\M_{l+1,l+2}\times \cdots \times
E^\M_{l+n-1,l+n} \mid \\
& t^\M(e_i) = s^\M(e_{i+1}),i=1,2,\dots,n-1\}
\endalign
$$
and similarly $E^\K_{l,l+n}$.
It suffices to show that for 
$(e_1,\dots,e_n) \in E^\M_{l,l+n}$
there exists
$(g_1,\dots,g_n)\in  E^\K_{l,l+n}$
such that
$p_X(g_i) =e_i,i=1,2,\dots,n.$
Take $f_n \in E^\N_{l+n,l+n+1}$
such that 
$t^\M(e_n) = s^\N(f_n).$
Since $\T$ is LR,
there uniquely exists
$f_{n-1} \in E^\N_{l+n-1,l+n}$
and
${e'}_n \in E^\M_{l+n,l+n+1}$
such that the quadruple
$(f_{n-1}, f_n, e_n, {e'}_n)$ denoted by $g_n$
gives rise to an element of $E^\K_{l+n-1,l+n}.$
One may inductively find 
$f_{k} \in E^\N_{l+k,l+k+1}$,
${e'}_k \in E^\M_{l+k,l+k+1}$
for $k=1,2,\dots,n$
and
$f_{0} \in E^\N_{l,l+1}$
 such that the quadruple
$(f_{k-1}, f_k, e_k, {e'}_k)$ denoted by $g_k$
gives rise to an element of $E^\K_{l+k-1,l+k}$
for $k=1,2,\dots,n$.
They satisfy 
$(g_1,\dots,g_n) \in E^\K_{l+1,l+n+1}$
and
$p_X(g_i) =e_i,i=1,2,\dots,n.$
One also sees that
 $q_X: X_{{\frak L}^\K} \longrightarrow
X_{{\frak L}^\M_0}
$
is surjective in a smilar way.
\qed
\enddemo

\proclaim{Proposition 5.3}
Let $\T$ be an LR textile $\lambda$-graph system.
For $k,n \ge 1$, we have
$$
({\check{\Cal U}}_{\Cal T}^{(k,n)},\check{\sigma}_{\Cal T}^{(k,n)})
= (\Lambda_{\N^k \M^n},\sigma_{\N^k \M^n}),
\qquad
({\hat{\Cal U}}_{\Cal T}^{(k,n)},\hat{\sigma}_{\Cal T}^{(k,n)})
=
(\Lambda_{ \M^n\N^k},\sigma_{\M^n\N^k}).
$$
\endproclaim
\demo{Proof}
We will prove that
the  map
$
\check{\Theta}_{\Cal T}^{(k,n)}:{\Cal U}_{\Cal T} \rightarrow 
\Lambda_{\N^k\M^n}
$
is surjective so that the first equality holds.
Take an arbitrary sequence ${(a_i)}_{i\in \Bbb Z} \in \Lambda_{\N^k\M^n}$.
Since $\T$ is LR, there exists a two dimensional configuration
${(\alpha_{i,j})}_{(i,j)\in {\Bbb Z}^2}\in 
\prod \Sb {(i,j)\in {\Bbb Z}^2} \endSb \Sigma^\K$ 
such that
by putting
$$
\alpha_i^h= {(\alpha_{i,j})}_{j\in \Bbb Z},\qquad
\alpha_j^v= {(\alpha_{i,j})}_{i\in \Bbb Z}
$$
$\alpha_i^h$ belongs to $\Lambda_\K$ for all $i \in \Bbb Z$ and
$\alpha_j^v$ belongs to $\Lambda_{\K^*}$ for all $j \in \Bbb Z$, 
that satisfy
$$
\xi(\alpha_i^h )= \eta(\alpha_{i-1}^h)\quad\text{ for }i\in \Bbb Z,
\qquad 
\xi^*(\alpha_j^v ) = \eta^*(\alpha_{j-1}^v)\quad\text{ for }j\in \Bbb Z,
$$
and
$$
a_i = (p^*(\alpha_{ki,ni}),\dots,p^*(\alpha_{(k+1)i-1,ni}, 
q(\alpha_{(k+1)i-1,ni},\dots,q(\alpha_{(k+1)i-1,(n+1)i-1})
$$
for
$
 i \in \Bbb Z.
$
For $m \in \Bbb N$ and $((k+1)m-1,(n+1)m-1)\in \triangle$, 
we may take
$$
\align
& (e_{ki,ni},\dots, e_{(k+1)i-1,ni}, e_{(k+1)i-1,ni},
\dots,
 e_{(k+1)i-1,(n+1)i-1}) \tag 5.2\\
\in & E^{\K}_{ki + ni,ki + ni+1}\times\cdots\times
E^{\K}_{(k+1)i + ni-1,(k+1)i + ni}\times
E^{\K}_{(k+1)i + ni-1,(k+1)i + ni}\times \\
& \qquad\qquad \qquad \cdots\times
E^{\K}_{(k+1)i + (n+1)i-2, (k+1)i + (n+1)i-1},\quad
i=0,1,\dots,m
\endalign
$$
such that
$$
\align
&(\lambda^\K(e_{ki,ni}),\dots, \lambda^\K(e_{(k+1)i-1,ni}), 
\lambda^\K(e_{(k+1)i-1,ni}),
\dots,
\lambda^\K( e_{(k+1)i-1,(n+1)i-1}))\\
=
&(\alpha_{ki,ni},\dots, \alpha_{(k+1)i-1,ni}, \alpha_{(k+1)i-1,ni},
\dots,
 \alpha_{(k+1)i-1,(n+1)i-1}),\quad
i=0,1,\dots,m.
\endalign
$$
 Since ${\frak L}^\K$ and ${\frak L}^{\K^*}$
 are both left-resolving,
 edges of ${\frak L}^\K$ located in $\triangle_{((k+1)m-1,(n+1)m-1)}$
 are uniquely determined 
 by the edges (5.2) and the labels 
$
{(\alpha_{i,j})}_{(i,j)\in \triangle_{((k+1)m-1,(n+1)m-1)}}.
$
Hence we know that
$$
{(\alpha_{i,j})}_{(i,j)\in \triangle_{((k+1)m-1,(n+1)m-1)}}
\in {\Cal U}_{{\frak L}^\K}(\T,\triangle_{((k+1)m-1,(n+1)m-1)}).
$$
Since $m$ is arbitrary,
we have
${(\alpha_{i,j})}_{(i,j)\in \triangle} \in {\Cal U}_{{\frak L}^{\K}}.
$ 
By applying the above discussion to the sequences
${(a_{i-K})}_{i\in \Bbb Z} \in \Lambda_{\N^k \M^n}$ for
$K \in \Bbb Z$,
one has that
$
{(\alpha_{i,j})}_{(i,j)\in \triangle} \in {\Cal U}_{{\frak L}^{\K}}^\infty,
$ 
and
$
{(\alpha_{i-K,j})}_{(i,j)\in \triangle} \in {\Cal U}_{{\frak L}^{\K}}^\infty
$ 
and
$
{(\alpha_{i,j-K})}_{(i,j)\in \triangle} \in {\Cal U}_{{\frak L}^{\K}}^\infty
$ 
for
$K \in \Bbb Z$.
Thus we have
$
{(\alpha_{i,j})}_{(i,j)\in {\Bbb Z}^2} \in {\Cal U}_{\Cal T}
$ 
and
$
\check{\Theta}_{\Cal T}^{(k,n)}({(\alpha_{i,j})}_{(i,j)\in{\Bbb Z}^2}) = 
{(a_i)}_{i\in \Bbb Z}.
$
Therefore we conclude that 
$$
{\check{\Cal U}}_{\Cal T}^{(k,n)}
=\Lambda_{\N^k \M^n}\qquad \text{ and } \qquad 
{\check{\sigma}}_{\Cal T}^{(k,n)}
=\sigma_{\N^k \M^n}.
$$
The other equality is similarly  proved.
\qed
\enddemo

\proclaim{Proposition 5.4 (cf.[N2;Lemma 6.2])}
Let $\T$ be an LR textile $\lambda$-graph system.
Then for $k,n \ge 1$ the map
$\check{\Theta}^{(k,n)}_{\Cal T}:{\Cal U}_{\Cal T}\rightarrow
\check{\Cal U}^{(k,n)}_{\Cal T}
$
is injective. 
Hence it gives rise to a topological conjugacy 
between
$({\Cal U}_{\Cal T},\sigma^{(k,n)}_{\Cal T})$ and
$(\check{\Cal U}^{(k,n)}_{\Cal T},\check{\sigma}^{(k,n)}_{\Cal T}).$
Similarly we have a topological conjugacy between
$({\Cal U}_{\Cal T},\sigma^{(k,n)}_{\Cal T})$ and
$(\hat{\Cal U}^{(k,n)}_{\Cal T},\hat{\sigma}^{(k,n)}_{\Cal T}).$
\endproclaim
\demo{Proof}
By a similar way to the proof of 
[N2;Lemma 6.2], 
we can show that
for ${(a_i)}_{i\in \Bbb Z} \in \check{\Cal U}_{\Cal T}^{(k,n)}$
there uniquely exists
${(\alpha_{i,j})}_{(i,j)\in {\Bbb Z}^2} \in {\Cal U}_{\Cal T}$
such that
$\check{\Theta}_{\Cal T}^{(k,n)}({(\alpha_{i,j})}_{(i,j)\in {\Bbb Z}^2}) = 
{(a_i)}_{i\in \Bbb Z}.$
\qed
\enddemo

We note that if $\T$ is LR, then ${\T}^*$ is LR.
We provide the following lemma.
\proclaim{Lemma 5.5}
A $1-1$ LR textile $\lambda$-graph system is surjective.
\endproclaim
\demo{Proof} 
Let
$\T$
be a $1-1$ LR textile $\lambda$-graph system.
Since ${\T}^*$ is LR, the both $\T$ and ${\T}^*$ are nondegenerate.
We will prove that 
the map
$$
\Phi_{\Cal T}:{\Cal U}_{\Cal T}\rightarrow \Lambda_{\K}
$$
is surjective.
For ${(a_j)}_{j\in \Bbb Z} \in \Lambda_\K$,
take
${(e_j)}_{j\in \Zp}\in 
X_{{\frak L}^\K} = \{{( e_j)}_{j\in \Zp} \mid 
e_j\in  E^\K_{j,j+1}, t(e_j) = s(e_{j+1}), j\in \Zp \}$
such that
$a_j = \lambda^\K(e_j), j\in \Zp$. 
Recall
 ${\triangle} = \{ (i,j) \in {\Bbb Z}^2 \mid i+j \ge 0 \}.$
We set
$$
\align
\triangle_{r,u} & =
\{(i,j) \in \triangle \mid  i \le 0, 0 \le j \},\\
\triangle_{l,d} & =
\{(i,j) \in \triangle \mid  1 \le i, j\le -1 \},\\
\square_{r,d} & =
\{(i,j) \in \triangle \mid  1 \le i, 0 \le j \}.
\endalign
$$
Now $\T$ is $1-1$  so that
there uniquely exists
$\alpha_{i,j} \in \Sigma^\K$ 
for $ (i, j) \in {\Bbb Z}^2$
such that
by putting
$\alpha_i = {(\alpha_{i,j})}_{j\in \Bbb Z}$
one has $\alpha_i \in \Lambda_\K$,
$\alpha_0 = {(a_j)}_{j\in \Bbb Z}$
and
$\eta( \alpha_i) = \xi(\alpha_{i+1})$ for $i \in \Bbb Z$.
Take an arbitrary $(n,k) \in \triangle$.
We set
$$
\align
\triangle_{r,u}(k) & =
\{(i,j) \in \triangle_{r,u} \mid  j \le k \},\\
\triangle_{l,d}(n) & =
\{(i,j) \in \triangle_{l,d} \mid  i \le n \},\\
\square_{r,d}(n,k) & =
\{(i,j) \in \square_{r,d} \mid  i \le n,  j \le k \}.
\endalign
$$
Take 
$f_k(i) \in E_{i+k+1,i+k+2}^\N, i=1,2,\dots,n$
such that
there exist
$
e_{i,j}\in E^\K_{i+j,i+j+1}
$
for 
$
(i,j) \in \square_{r,d}(n,k) \cup \triangle_{l,d}(n)
$
satisfying
$$
\align
{(e_{i,j})}_{(i,j) \in \square_{r,d}(n,k) \cup \triangle_{l,d}(n)} 
&\in
X({\Cal T}_{{\Cal K}^\M_\N}; \square_{r,d}(n,k) \cup \triangle_{l,d}(n)),\\
p^E(e_{1,j}) = & q^E(e_j)\qquad \text{ for } j=1,2,\dots,k,\\
t(e_{i,k}) =  & f_k(i)\qquad \text{ for } i=1,2,\dots,n,\\
\alpha_{i,j} = & \lambda^\K(e_{i,j})\qquad \text{ for } (i,j) 
\in \square_{r,d}(n,k) \cup \triangle_{l,d}(n)).
\endalign
$$
As ${\frak L}^\K$ is left-resolving, such edges
$
e_{i,j}\in E^\K_{i+j,i+j+1}
$
for 
$
(i,j) \in 
\square_{r,d}(n,k) \cup \triangle_{l,d}(n)
$
are unique for
$f_k(i) \in E_{i+k+1,i+k+2}^\N, i=1,2,\dots,n.$
We set $e_{0,j} = e_j$ for $j=0,1,\dots,k$.
Since ${\frak L}^{\K^*}$ is left-resolving,
the vertices
${p^E(e_j)}, j= 0,1,\dots k$
of ${\frak L}^{\K^*}$ and labels
${(\alpha_{i,j})}_{(i,j)\in \triangle_{r,u}(k)}$
uniquely determine
edges
$e_{i,j}\in E^\K_{i+j,i+j+1}$
for
$(i,j) \in \triangle_{r,u}(k)$
such that
$t^{\K^*}(e_{1,j}) = s^{\K^*}(e_j), j=1,\dots,k$
and
$\lambda^\K(e_{i,j}) =\alpha_{i,j}$ 
for
$(i,j) \in \triangle_{r,u}(k)$.
Hence we have
$$
{(e_{i,j})}_{(i,j)\in \triangle_{r,u}(k)} \in X(\T;\triangle_{r,u}(k))
\quad
\text{ and hence }
\quad
{(e_{i,j})}_{(i,j)\in \triangle} \in X(\T;\triangle)
$$
so that
$$
{(\alpha_{i,j})}_{(i,j)\in \triangle(n,k)} 
\in {\Cal U}_{{\frak L}^\K}(\triangle(n,k))
$$
By Proposition 4.3,
the configuration
${(\alpha_{i,j})}_{(i,j) \in \triangle}$
belongs to
${\Cal U}_{{\frak L}^\K}$.
By applying this argument to the configurations
${(a_{i+k,j})}_{(i,j)\in {\Bbb Z}^2}$
and 
${(a_{i,j+k})}_{(i,j)\in {\Bbb Z}^2}$
for
$k \in \Bbb Z$,
we know that
${(\alpha_{i,j})}_{(i,j) \in \triangle}$
belongs to
${\Cal U}_{{\frak L}^\K}^\infty$ and to ${\Cal U}_{\Cal T}$.
Since
$\Phi_{\Cal T}({(\alpha_{i,j})}_{(i,j) \in {\Bbb Z}^2}) 
= {(a_j)}_{j \in \Bbb Z}$,
the map
$
\Phi_{\Cal T}:{\Cal U}_{\Cal T}\rightarrow \Lambda_{\K}
$
is surjective.
\qed
\enddemo
Therefore we obtain   
\proclaim{Theorem 5.6}
Let
$\T$
be a $1-1$ LR textile $\lambda$-graph system defined by a
specified equivalence:
$$
\M_{l,l+1}\N_{l+1,l+2} \overset{\kappa}\to{\simeq} \N_{l,l+1}\M_{l+1,l+2},
\qquad l \in \Zp.
$$
Then the dynamical system 
$(\Lambda_{\M}, \varphi_{\Cal T}^k\sigma_{\M}^n),$
$ k\ge 0,n\ge 1$
is topologically conjugate to the subshift
$(\Lambda_{\N^k\M^n},\sigma_{\N^k\M^n})$
 presented by the symbolic matrix system 
$(\N^k\M^n,I^{\N^k\M^n}),$
defined by
$$
\align
{(\N^k\M^n)}_{l,l+1} 
=& 
\N_{l(k+n),l(k+n)+1}\cdots \N_{l(k+n)+n-1,l(k+n)+n}\cdot\\
& \cdot \M_{l(k+n)+n,l(k+n)+n+1}\cdots\M_{(l+1)(k+n)-1,(l+1)(k+n)},\\
{I^{\N^k\M^n}}_{l,l+1} 
=& 
I^\N_{l(k+n),l(k+n)+1}\cdots I^\N_{l(k+n)+n-1,l(k+n)+n}\cdot\\
& \cdot I^\M_{l(k+n)+n,l(k+n)+n+1}\cdots I^\M_{(l+1)(k+n)-1,(l+1)(k+n)},\qquad l \in \Zp.
\endalign
$$
\endproclaim
\demo{Proof}
Since $\T$ is nondegenerate,
for the case when $k=0$
the assertion is clear.
We may assume that $k\ge 1$.
Since $\T$ is $1-1$ and LR, it is surjective by Lemma 5.4 
so that
$(\Lambda_{\M},\varphi_{\Cal T}^k\sigma_{\M}^n)$
is conjugate to 
$({\Cal U}_{\Cal T},\sigma_{\Cal T}^{(k,n)})$ 
by Proposition 4.8.
As $\T$ is LR and $k,n\ge 1$, one has that
  $({\Cal U}_{\Cal T},\sigma_{\Cal T}^{(k,n)})$
  is conjugate to
$(\check{\Cal U}_{\Cal T}^{(k,n)},\check{\sigma}_{\Cal T}^{(k,n)})$
by Proposition 5.4.
Hence by Proposition 5.3,
we obtain the assertion.
\qed
\enddemo

\heading 6. LR textile systems and properly strong shift equivalences
\endheading

Let $(\M,I)$ and $(\M^{\prime},I')$ be symbolic matrix systems over 
$\Sigma$ and 
$\Sigma^{\prime}$ respectively.

\noindent
{\bf Definition ([Ma2]).}
 $(\M,I)$ and $(\M^{\prime},I')$
 are said to be {\it properly strong shift equivalent in 1-step} 
 if
 there exist alphabets 
$C,D$ and
 specifications 
$
 \kappa: \Sigma \rightarrow C D,
$
$
 \kappa': \Sigma' \rightarrow D  C
$
 and increasing sequences $n(l),n'(l)$ on $l \in \Zp$
 such that
 for each $l\in \Zp$, there exist
an $n(l)\times n'(l+1)$ matrix ${\P}_l$ over $C,$ 
an $n'(l)\times n(l+1)$ matrix ${\Q}_l$ over $D,$
an $n(l)\times n(l+1)$ matrix $X_l $ over $ \{0,1\}$
and 
an $n'(l)\times n'(l+1)$ matrix $Y_l$ over $ \{0,1\}$
 satisfying the following equations: 
$$
\align
{\M}_{l,l+1} 
\overset{\kappa}\to {\simeq}&  {\P}_{2l}{\Q}_{2l+1},\qquad 
\M'_{l,l+1} 
\overset{\kappa'}\to {\simeq} {\Q}_{2l}{\P}_{2l+1},\\ 
I_{l,l+1} = & X_{2l}X_{2l+1}, \qquad 
I'_{l,l+1} = Y_{2l}Y_{2l+1} 
\endalign
$$
and
$$
X_l{\P}_{l+1}= {\P}_{l}{Y}_{l+1},\qquad
Y_l{\Q}_{l+1}= {\Q}_{l}{X}_{l+1}. 
$$
\medskip
We write this situation as
$
(\P,\Q,X,Y): (\M,I)\underset{1-pr}\to {\approx} (\M',I').
$

We in particular consider the case when 
$(\M',I') = (\M,I)$.
\proclaim{Lemma 6.1} 
Suppose that 
$
(\P,\Q,X,Y): (\M,I)\underset{1-pr}\to {\approx} (\M,I).
$
Put
$$
\align
\P_{l,l+1} & = \P_{2l} Y_{2l+1} (= X_{2l}\P_{2l+1}),
\qquad
I^\P_{l,l+1} = I_{l,l+1},\\
{\Q}_{l,l+1} & = \Q_{2l} X_{2l+1} (= Y_{2l}\Q_{2l+1}),
\qquad
I^{\Q}_{l,l+1}  = I_{l,l+1},\qquad l \in \Zp.
 \endalign
$$
 Then we have
\roster
\item"(i)"
$(\P,I^\P)= {(\P_{l,l+1},I^\P_{l,l+1})}_{l\in \Zp}$ and
$({\Q},I^{\Q})= {({\Q}_{l,l+1},I^{\Q}_{l,l+1})}_{l\in \Zp}$ 
are symbolic matrix systems over $C$ and $D$ respectively. 
\item"(ii)"
The pair
${\frak L}^\M$ and ${\frak L}^\P$,
and the pair  
${\frak L}^\M$ and ${\frak L}^{\Q}$ both form squares such 
that 
$$
\align
\M_{l,l+1}\P_{l+1,l+2} \overset{\kappa^\P}\to {\simeq} & 
\P_{l,l+1}\M_{l+1,l+2},\qquad l \in \Zp, \tag 6.1\\
\M_{l,l+1}{\Q}_{l+1,l+2} \overset{\kappa^\Q}\to {\simeq} & 
{\Q}_{l,l+1}\M_{l+1,l+2},\qquad l \in \Zp \tag 6.2
\endalign
$$
for some specifications 
$\kappa^\P : \Sigma C \longrightarrow C \Sigma$ and
$\kappa^\Q : \Sigma D \longrightarrow D \Sigma$.
\endroster
\endproclaim
\demo{Proof}We will prove the assertions for  $(\P,I^\P)$.
The assertions for the other one is symmetric. 

(i) The equality
$
\P_{l,l+1} I^\P_{l+1,l+2}
 =I^\P_{l,l+1}\P_{l+1,l+2}
$
is easily shown.
Hence 
$(\P,I^\P)= {(\P_{l,l+1},I^\P_{l,l+1})}_{l\in \Zp}$ 
is a symbolic matrix system over $C$. 

(ii) One has 
$$
\align
\M_{l,l+1} \P_{l+1,l+2}
& \overset{\kappa}\to {\simeq} \P_{2l} \Q_{2l+1} X_{2l+2} \P_{2l+3}\\
& = \P_{2l} Y_{2l+1} \Q_{2l+2} \P_{2l+3} \overset{{\kappa'}^{-1}}\to {\simeq} \P_{l,l+1} \M_{l+1,l+2}.
\endalign
$$
For $\alpha \in \Sigma, c \in C$, by putting 
$\kappa^\P(\alpha c) = c_\alpha {\kappa'}^{-1}(d_\alpha c )$
where $\kappa (\alpha ) = c_\alpha d_\alpha \in C D$,
the specification $\kappa^\P :\Sigma C \longrightarrow C \Sigma$ yields the desired specified equivalence (6.1).
\qed
\enddemo
By this lemma with Proposition 5.1, 
the relations (6.1) and (6.2) yield  LR textile $\lambda$-graph systems
$
 {\Cal T}_{{\K}^{\M}_{\P}}
$
 and
$
 {\Cal T}_{{\K}^{\M}_{\Q}}
$
respectively.
\proclaim{Lemma 6.2}
Suppose that 
$
(\P,\Q,X,Y): (\M,I)\underset{1-pr}\to {\approx} (\M,I).
$
Keep the notations as in the preceding lemma.
The LR textile $\lambda$-graph systems
$
 {\Cal T}_{{\K}^{\M}_{\P}}
$
 and
$
 {\Cal T}_{{\K}^{\M}_{\Q}}
$
are both $1-1$ and hence surjective.
\endproclaim
\demo{Proof}
We will prove that $
 {\Cal T}_{{\K}^{\M}_{\P}}
$
 is $1-1.$
The LR textile system defined by the specified equivalence (6.1)
comes from the specified equivalence 
$$
\P_{2l}\Q_{2l+1}\cdot \P_{2l+2}Y_{2l+3} \overset{\kappa^\P}\to{\simeq}
\P_{2l}Y_{2l+1}\cdot \Q_{2l+2}\P_{2l+3}.\tag 6.3
$$
Let ${(\beta_l)}_{l\in \Bbb Z}\in \Lambda_{{\frak L}^\K}$
 be such that
$
\xi ({(\beta_l)}_{l\in \Bbb Z}) = {(\alpha_l)}_{l\in \Bbb Z}.
$
We put
$\kappa(\alpha_l) = c_l d_l \in CD$ 
for
$
l \in \Bbb Z.
$ 
By (6.3) $\beta_l$ is uniquely determined by the square:
$$
\CD
 \cdot  @>c_l>>\cdot  @>d_l>> \cdot\\
@VV{c_l}V   @.    @VV{c_{l+1}}V \\   
\cdot @>d_l>> \cdot @>c_{l+1}>>\cdot\\
\endCD
$$
That is $\beta_l$ is uniquely determined by the quadruple
$(c_l, c_{l+1}, c_l d_l, d_lc_{l+1}) \in \Sigma^{\K},$
that are determined by the sequence
${(c_l d_l)}_{l\in \Bbb Z}$.
Hence the code
$\xi : \Lambda_{\K}\longrightarrow \Lambda_{\M}$
is one-to-one.
We similarly see that 
$\eta : \Lambda_{\K}\longrightarrow \Lambda_{\M}$
is one-to-one.
Hence by Lemma 5.5, 
$
 {\Cal T}_{{\K}^{\M}_{\P}}
$
is surjective.
We symmetrically see that 
$
 {\Cal T}_{{\K}^{\M}_{\Q}}
$
 is  $1-1$ and surjective.
\qed
\enddemo

Following Nasu's notation [N],[N2], an automorphism $\phi$ 
of a subshift $\Lambda$ 
over
$\Sigma$ is called a forward bipartite automorphism if there exist
alphabets $C, D$ and specifications
$\kappa:\Sigma \rightarrow CD,
\kappa':\Sigma \rightarrow DC$
such that 
$\phi$ is given by 
$$
\phi({(\alpha_l)}_{l \in \Bbb Z}) = 
{({\alpha'}_l)}_{l \in \Bbb Z},
\qquad {(\alpha_l)}_{l \in \Bbb Z} \in \Lambda
$$
where
${\kappa}(\alpha_l)= c_l  d_l $ for some $c_l \in C$ 
and 
$d_l \in  D, l\in \Bbb Z$ and
$
{\alpha'}_l = {\kappa'}^{-1}(d_l c_{l+1}) \in \Sigma.
$
Hence a a properly strong shift equivalence 
$
(\P,\Q,X,Y): (\M,I)\underset{1-pr}\to {\approx} (\M,I)
$
in $1$-step gives rise to a forward bipartite automorphism on the subshift 
presented by $(\M,I)$. 
\proclaim{Lemma 6.3}
Let $(\Lambda,\sigma)$ be a subshift presented by $(\M,I)$.
Let $\phi$ be a forward bipartite automorphism on $(\Lambda,\sigma)$
defined by a properly strong shift equivalence 
$
(\P,\Q,X,Y): (\M,I)\underset{1-pr}\to {\approx} (\M,I)
$
in $1$-step.
Let ${\Cal T}^{\P}$ and ${\Cal T}^{\Q}$ be the LR textile $\lambda$-graph systems 
$
 {\Cal T}_{{\K}^{\M}_{\P}}
$
and
$
 {\Cal T}_{{\K}^{\M}_{\Q}}
$
defined by the relations (6.1) and (6.2) respectively.
Then we have
$$
\varphi_{{\Cal T}^{\P}} =\phi,
\qquad 
\varphi_{{\Cal T}^{\Q}} =\phi^{-1}\circ \sigma
$$
as automorphisms on $\Lambda = \Lambda_{\M}.$
\endproclaim
\demo{Proof}
We will prove that 
$\phi = \varphi_{{\Cal T}^\P}.$
For ${(\alpha_l)}_{l \in \Bbb Z} \in \Lambda,$
put
$c_l  d_l = \kappa(\alpha_l)\in C D, l\in \Bbb Z.$
By setting
$
{\alpha'}_l = {\kappa'}^{-1}(d_l c_{l+1}) \in \Sigma,
$
one has 
$\phi({(\alpha_l)}_{l \in \Bbb Z}) = 
{({\alpha'}_l)}_{l \in \Bbb Z}.$
Put
$
\beta_l = (c_l,c_{l+1}, \kappa^{-1}(c_l  d_l), {\kappa'}^{-1}(d_l c_{l+1})) 
\in \Sigma^\K,
l\in \Bbb Z$
so that
one has
${(\beta_l)}_{l\in \Bbb Z} \in 
\Lambda_{\K}$
and
$$
\xi ( {(\beta_l)}_{l\in \Bbb Z} ) = {(\kappa^{-1} (c_l d_l ))}_{l\in \Bbb Z},
\qquad
\eta ( {(\beta_l)}_{l\in \Bbb Z} ) ={({\kappa'}^{-1} (d_l c_{l+1}))}_{l\in \Bbb Z}.
$$
It then follows that
$$
\phi({(\alpha_l)}_{l \in \Bbb Z}) = 
{({\alpha'}_l)}_{l \in \Bbb Z}
= {({\kappa'}^{-1} (d_l c_{l+1} ))}_{l\in \Bbb Z}
=\eta ( {(\beta_l)}_{l\in \Bbb Z} )
= \eta \circ \xi^{-1}( {(\alpha_l)}_{l\in \Bbb Z} ).
$$
The equality
$\varphi_{{\Cal T}^{\Q}} =\phi^{-1}\circ \sigma
$
is similarly shown.
\qed
\enddemo

We assume that the previously defined metric is equipped with $\Lambda$.
Then the homeomorphism $\sigma$ has $1$ as its expansive constant.
Therefore we have
\proclaim{Theorem 6.4}
Let $(\Lambda,\sigma)$ be a subshift presented by a symbolic matrix system 
$(\M,I)$.
Let $\phi$ be a forward bipartite automorphism on $(\Lambda,\sigma)$
defined by a properly strong shift equivalence 
$
(\P,\Q,X,Y): (\M,I)\underset{1-pr}\to {\approx} (\M,I)
$
in $1$-step.
If $\phi$ is expansive with $\frac{1}{k}$ as its expansive constant for some $k \in \Bbb N$,
 the dynamical system
$(\Lambda,\phi)$
is topologically conjugate to the subshift 
$(\Lambda_{\P_{\Cal T}^{[2k]}},\sigma_{\P_{\Cal T}^{[2k]}})$
presented by the symbolic matrix system
$(\P_{\Cal T}^{[2k]},I^{\P_{\Cal T}^{[2k]}})$
where
$(\P_{\Cal T}^{[2k]},I^{\P_{\Cal T}^{[2k]}})$
is the $2k$-heigher block of the symbolic matrix system
$(\P,I^\P)$ relative to the LR textile $\lambda$-graph system 
$ {\Cal T}_{{\K}^{\M}_{\P}}$
defined by the specified equivalence 
$$
\M_{l,l+1}\P_{l+1,l+2} \overset{\kappa^\P}\to {\simeq} 
\P_{l,l+1}\M_{l+1,l+2},\qquad l \in \Zp,
$$
where
$$
\P_{l,l+1}  = \P_{2l} Y_{2l+1} (= X_{2l}\P_{2l+1}),\qquad 
I^\P_{l,l+1} = I_{l,l+1}, \qquad l \in \Zp.
$$
\endproclaim
\demo{Proof}
Consider the LR textile $\lambda$-graph system 
${\Cal T}^\P = {\Cal T}_{{\K}^{\M}_{\P}}$ defined by (6.1),
that is $1-1$ and surjective by Lemma 6.2.
Lemma 6.3 says that  $\varphi_{{\Cal T}^\P} = \phi$ on $\Lambda.$
By the assumption on $\phi$ and Theorem 4.10 (i),
${{\Cal T}_{{\K}^{\M}_{\P}} }^{{[2k]}^*}$ is $1-1$.
Since ${\Cal T}_{{\K}^{\M}_{\P}} $ is LR, 
${{\Cal T}_{{\K}^{\M}_{\P}} }^{[2k]}$ 
and
hence 
${{\Cal T}_{{\K}^{\M}_{\P}} }^{{[2k]}^*}$ 
are both LR and nondegenerate. 
By Lemma 5.5,
${{\Cal T}_{{\K}^{\M}_{\P}} }^{{[2k]}^*}$ is surjective.
By Theorem 4.10 (ii),
the topological dynamical system
$
(\Lambda, \varphi_{{\Cal T}^\P})
$
is realized as the subshift
$ 
(\Lambda_{{\P_{\Cal T}^{[2k]}}}, \sigma_{{{\P_{\Cal T}^{[2k]}}}}).
$Hence,
one concludes that
 the dynamical system
$(\Lambda,\phi)$
is topologically conjugate to the subshift 
$(\Lambda_{\P_{\Cal T}^{[2k]}},\sigma_{\P_{\Cal T}^{[2k]}})$.
\qed
\enddemo

We also have 
\proclaim{Theorem 6.5}
Let $(\Lambda,\sigma)$ be a subshift presented by a symbolic matrix system 
$(\M,I)$.
Let $\phi$ be a forward bipartite automorphism on $(\Lambda,\sigma)$
defined by a properly strong shift equivalence 
$
(\P,\Q,X,Y): (\M,I)\underset{1-pr}\to {\approx} (\M,I)
$
in $1$-step.
Then the dynamical system
$(\Lambda,\phi^k\sigma^n)$
is topologically conjugate to the subshift 
$(\Lambda_{\P^k \M^n},\sigma_{\P^k \M^n})$
presented by the symbolic matrix system
$({\P^k \M^n},I^{\P^k \M^n})$
 for
 $k\ge 0, n\ge 1$
defined by
$$
\align
{(\P^k\M^n)}_{l,l+1} 
=& 
\P_{l(k+n),l(k+n)+1}\cdots \P_{l(k+n)+n-1,l(k+n)+n}\cdot\\
& \cdot \M_{l(k+n)+n,l(k+n)+n+1}\cdots\M_{(l+1)(k+n)-1,(l+1)(k+n)},\\
{I^{k+n}}_{l,l+1} 
=& 
I_{l(k+n),l(k+n)+1}\cdots I_{(l+1)(k+n)-1,(l+1)(k+n)},\qquad l \in \Zp
\endalign
$$
where
$
\P_{l,l+1}  = \P_{2l} Y_{2l+1} (= X_{2l}\P_{2l+1})$
for
$
l \in \Zp.
$
And also 
$(\Lambda,{(\sigma\phi^{-1})}^k\sigma^n)$
is topologically conjugate to the subshift 
$(\Lambda_{\Q^k \M^n},\sigma_{\Q^k \M^n})$
presented by the similarly defined symbolic matrix system
$({\Q^k \M^n},I^{\Q^k \M^n})$
 for
 $k\ge 0, n\ge 1,$
\endproclaim
\demo{Proof}
By Lemma 6.2, the LR textile $\lambda$-graph system
${\Cal T}_{\K^\M_\P}$ is $1-1$ and surjective.
Hence by Theorem 5.6,
the dynamical system 
$(\Lambda, \varphi_{\Cal T}^k\sigma^n),$
$ k\ge 0,n\ge 1$
is topologically conjugate to the subshift
$(\Lambda_{\N^k\M^n},\sigma_{\N^k\M^n})$
presented  by the symbolic matrix system 
$(\N^k\M^n,I^{\N^k\M^n}).$
Now 
${\Cal T}_{\K^\M_\P}$
is LR so that it is nondegenerate.
By Lemma 6.3 one sees  that 
$\varphi_{{\Cal T}^{\P}} =\phi$ 
so that 
the dynamical system 
$(\Lambda, \phi_{\Cal T}^k\sigma^n)$
is topologically conjugate to the subshift
$(\Lambda_{\N^k\M^n},\sigma_{\N^k\M^n})$.
It is similarly shown that 
$(\Lambda,{(\sigma\phi^{-1})}^k\sigma^n)$
is topologically conjugate to the subshift 
$(\Lambda_{\Q^k \M^n},\sigma_{\Q^k \M^n})$
 for
 $k\ge 0, n\ge 1.$
\qed
\enddemo

\heading 7.  Subshift-identifications of automorphisms of subshifts 
\endheading

Two symbolic matrix systems
 $(\M,I)$ and $(\M^{\prime},I')$ 
 are said to be {\it properly strong shift equivalent} 
 if
there exists a finite sequence 
$(\M^{(1)}, I^{(1)}),\dots,(\M^{(N-1)}, I^{(N-1)})$
of symbolic matrix systems such that
$$
(\M,I) \underset{1-pr}\to {\approx} (\M^{(1)},I^{(1)})
       \underset{1-pr}\to {\approx} 
       \cdots 
       \underset{1-pr}\to {\approx} 
(\M^{(N-1)},I^{(N-1)}) \underset{1-pr}\to {\approx} 
( \M', I').
$$
In [Ma2], the following theorem has been proved:
\proclaim{Theorem 7.1([Ma2])}
If two symbolic matrix systems $(\M,I)$ and $(\M',I')$ are properly strong shift equivalent, then their respect  presented subshifts
$\Lambda_\M$ and $\Lambda_{\M'}$ are topologically conjugate. 
Furthermore, 
two  subshifts $\Lambda$ and $\Lambda'$ are topologically conjugate
if and only if  
their canonical symbolic matrix systems
$(\M^\Lambda, I^\Lambda)$ and
$(\M^{\Lambda'}, I^{\Lambda'})$  are properly strong shift equivalent.
\endproclaim
In particular,
an automorphism of a subshift
$\Lambda$ is given by  a 
 properly strong shift equivalence from a symbolic matrix system that presents the subshift to itself.
Let $(\M,I)$ be a symbolic matrix system over $\Sigma$. 
 Let us consider 
 a  properly strong shift equivalence 
 from $(\M,I)$ to itself.
 Hence we consider
 symbolic matrix systems
 $(\M^{(k)},I^{(k)})$
over
$\Sigma^{(k)}, k=0,1,\dots,N,$
where
$(\M^{(0)},I^{(0)}) =
(\M^{(N)},I^{(N)})
=(\M,I)$
and
$\Sigma^{(0)} = \Sigma^{(N)} = \Sigma$
such that
there exist alphabets 
$C^{(k)},D^{(k)}$ and
 specifications 
$
 \kappa_0^{(k)}: \Sigma^{(k-1)} \rightarrow C^{(k)} D^{(k)},
$
$
 \kappa_1^{(k)}: \Sigma^{(k)} \rightarrow D^{(k)}  C^{(k)}
$
 and increasing sequences $n_0^{(k)}(l),n_1^{(k)}(l)$ on $l \in \Zp$
 such that
 for each $l\in \Zp$, there exist
an $n_0^{(k)}(l)\times n_1^{(k)}(l+1)$ matrix $\P^{(k)}_l$ over $C^{(k)},$ 
an $n_1^{(k)}(l)\times n_0^{(k)}(l+1)$ matrix $\Q^{(k)}_l$ over $D^{(k)},$
an $n_0^{(k)}(l)\times n_0^{(k)}(l+1)$ matrix $X^{(k)}_l $ over $ \{0,1\}$
and 
an $n_1^{(k)}(l)\times n_1^{(k)}(l+1)$ matrix $Y^{(k)}_l$ over $ \{0,1\}$
 satisfying the following equations: 
$$
\cases
{\M}^{(k-1)}_{l,l+1} 
\overset{\kappa_0^{(k)}}\to {\simeq} \P^{(k)}_{2l}\Q^{(k)}_{2l+1},
& \qquad
{\M}^{(k)}_{l,l+1} 
\overset{\kappa_1^{(k)}}\to {\simeq} \Q^{(k)}_{2l}\P^{(k)}_{2l+1},\\ 
I^{(k-1)}_{l,l+1} = X^{(k)}_{2l}X^{(k)}_{2l+1},
& \qquad 
I^{(k)}_{l,l+1} = Y^{(k)}_{2l}Y^{(k)}_{2l+1},\\ 
X^{(k)}_l \P^{(k)}_{l+1}= \P^{(k)}_{l}Y^{(k)}_{l+1},
& \qquad
Y^{(k)}_l \Q^{(k)}_{l+1}= \Q^{(k)}_{l}X^{(k)}_{l+1}. 
\endcases
\tag 7.1
$$
The equations (7.1) are simply written as
$$
(\P^{(k)}, \Q^{(k)}, X^{(k)}, Y^{(k)}): 
(\M^{(k-1)},I^{(k-1)}) \underset{1-pr}\to {\approx} (\M^{(k)},I^{(k)}),
\qquad
k=1,\dots,N.
$$
\proclaim{Lemma 7.2}
Keep the notations as above.
Put $m(l)\times m(l+N)$ matrices
$$
\align
\P_{l,l+N}
&=\P^{(1)}_{2l}Y^{(1)}_{2l+1}\P^{(2)}_{2l+2}Y^{(2)}_{2l+3}\cdots \P^{(N)}_{2l+2N-2}Y^{(N)}_{2l+2N-1},\\
\Q_{l,l+N}
&=\Q^{(1)}_{2l}X^{(1)}_{2l+1}\Q^{(2)}_{2l+2}X^{(2)}_{2l+3}\cdots \Q^{(N)}_{2l+2N-2}X^{(N)}_{2l+2N-1},\\
I_{l,l+N} 
&=  I_{l,l+1}I_{l+1,l+2}\cdots I_{l+N-1,l+N},
\qquad l\in \Zp.
\endalign
$$
\roster
\item"(i)" The equalities 
$$
\align
\P_{l,l+N}I_{l+N,l+N+1} 
& = 
I_{l,l+1}\P_{l+1,l+N+1}, \\
\Q_{l,l+N}I_{l+N,l+N+1} 
& = 
I_{l,l+1}\Q_{l+1,l+N+1},\qquad l\in \Bbb Z
\endalign
$$
and hence
$$\align
\P_{l,l+N}I_{l+N,l+2N} 
& = 
I_{l,l+N}\P_{l+N,l+2N},\\
\Q_{l,l+N}I_{l+N,l+2N} 
& = 
I_{l,l+N}\Q_{l+N,l+2N},\qquad
l\in \Bbb Z
\endalign
$$ hold.
\item"(ii)"
There exist specifications
$$
\align
\kappa_\P: & \Sigma \cdot C^{(1)}C^{(2)}\cdots C^{(N)}  \longrightarrow
 C^{(1)}C^{(2)}\cdots C^{(N)}\cdot  \Sigma,\\
\kappa_\Q: & \Sigma \cdot D^{(1)}D^{(2)}\cdots D^{(N)} \longrightarrow
 D^{(1)}D^{(2)}\cdots D^{(N)}\cdot \Sigma
\endalign
$$
such that
$$
\align
\M_{l,l+1}\P_{l+1,l+N+1} & \overset{\kappa_\P}\to{\simeq}
\P_{l,l+N}\M_{l+N,l+N+1},\tag 7.2\\
\M_{l,l+1}\Q_{l+1,l+N+1} & \overset{\kappa_\Q}\to{\simeq}
\Q_{l,l+N}\M_{l+N,l+N+1},\qquad l\in \Zp.\tag 7.3
\endalign
$$
\endroster
\endproclaim
\demo{Proof}
(i)
We note that 
the equalities 
$$
I_{l,l+1} = I^{(0)}_{l,l+1} = X^{(1)}_{2l}X^{(1)}_{2l+1}
 = I^{(N)}_{l,l+1} =  Y^{(N)}_{2l}Y^{(N)}_{2l+1},\qquad l\in \Zp
$$
hold.
It then follows that
$$
\align
& \P_{l,l+N}I_{l+N,l+N+1} \\
=&  \P^{(1)}_{2l}Y^{(1)}_{2l+1}\P^{(2)}_{2l+2}Y^{(2)}_{2l+3}
 \cdots 
Y^{(N-1)}_{2l+2N-3} X^{(N)}_{2l+2N-2}\P^{(N)}_{2l+2N-1}Y^{(N)}_{2l+2N}
Y^{(N)}_{2l+2N+1}\\
=&  \P^{(1)}_{2l}Y^{(1)}_{2l+1}\P^{(2)}_{2l+2}Y^{(2)}_{2l+3}
 \cdots 
Y^{(N-1)}_{2l+2N-3} X^{(N)}_{2l+2N-2}X^{(N)}_{2l+2N-1}\P^{(N)}_{2l+2N}
Y^{(N)}_{2l+2N+1}\\
=&  \P^{(1)}_{2l}Y^{(1)}_{2l+1}\P^{(2)}_{2l+2}Y^{(2)}_{2l+3}
 \cdots 
 Y^{(N-1)}_{2l+2N-3}\cdot I^{(N-1)}_{l+N-1,l+N}\P^{(N)}_{2l+2N}
 Y^{(N)}_{2l+2N+1}\\
 \endalign
$$
and hence inductively
$$
\align
& \P_{l,l+N}I_{l+N,l+N+1} \\
=& \P_{2l}^{(1)} Y^{(1)}_{2l+1}I^{(1)}_{l+1,l+2}\cdot
\P^{(2)}_{2l+4}Y^{(2)}_{2l+5}\cdots
\P^{(N-1)}_{2l+2N-2}Y^{(N-1)}_{2l+2N-1}\P^{(N)}_{2l+2N}Y^{(N)}_{2l+2N+1}.\\
 \endalign
$$
Since
$$
 \P_{2l}^{(1)} Y^{(1)}_{2l+1}I^{(1)}_{l+1,l+2}
=  X^{(1)}_{2l}\P^{(1)}_{2l+1}\cdot Y^{(1)}_{2l+2}Y^{(1)}_{2l+3} 
=  I_{l,l+1} \P^{(1)}_{2l+2}Y^{(1)}_{2l+3}, 
$$
one has
$$
 \P_{l,l+N}I_{l+N,l+N+1} =I_{l,l+1} \P_{l+1,l+N+1}.
$$
One then inductively gets the equalities
$$
\P_{l,l+N}I_{l+N,l+2N} = 
I_{l,l+N}\P_{l+N,l+2N},\qquad l\in \Zp.
$$
The other equalities
$$
 \Q_{l,l+N}I_{l+N,l+N+1} =I_{l,l+1} \Q_{l+1,l+N+1},\qquad l\in \Zp
$$
are similarly proved.

(ii)
It follows that
$$
\align
& \M_{l,l+1}{\P}_{l+1,l+N+1} \\
 =& {\M}^{(0)}_{l,l+1} \P^{(1)}_{2(l+1)}Y^{(1)}_{2(l+1)+1}
\P^{(2)}_{2(l+1)+2}Y^{(2)}_{2(l+1)+3}
\cdots 
\P^{(N)}_{2l+2N}Y^{(N)}_{2l+2N+1}\\
\overset{\kappa_0^{(1)}}\to {\simeq} 
 & \P^{(1)}_{2l}\Q^{(1)}_{2l+1}
X^{(1)}_{2(l+1)}\P^{(1)}_{2(l+1)+1}
X^{(2)}_{2(l+1)+2}\P^{(2)}_{2(l+1)+3}
\cdots 
X^{(N)}_{2l+2N}\P^{(N)}_{2l+2N+1}\\
= & \P^{(1)}_{2l}Y^{(1)}_{2l+1}
\Q^{(1)}_{2(l+1)}\P^{(1)}_{2(l+1)+1}
X^{(2)}_{2(l+1)+2}\P^{(2)}_{2(l+1)+3}
\cdots 
X^{(N)}_{2l+2N}\P^{(N)}_{2l+2N+1}\\
\overset{{\kappa_1^{(1)}}^{-1}}\to {\simeq} 
&  \P^{(1)}_{2l}Y^{(1)}_{2l+1}
\cdot
\M^{(1)}_{l+1,l+2}
\P^{(2)}_{2(l+1)+2}Y^{(2)}_{2(l+1)+3}
\cdots 
\P^{(N)}_{2l+2N}Y^{(N)}_{2l+2N+1}\\
 \endalign
$$
and similarly
$$
\align
& {\M}^{(1)}_{l+1,l+2} \P^{(2)}_{2(l+1)+2}Y^{(2)}_{2(l+1)+3}
\P^{(3)}_{2(l+1)+4}Y^{(3)}_{2(l+1)+5}
\cdots 
\P^{(N)}_{2l+2N}Y^{(N)}_{2l+2N+1}\\
\overset{{\kappa_1^{(2)}}^{-1}\kappa_0^{(2)}}\to {\simeq} 
&  \P^{(2)}_{2l+2}Y^{(2)}_{2l+3}
\cdot
\M^{(2)}_{l+2,l+3}
\P^{(3)}_{2(l+1)+4}Y^{(3)}_{2(l+1)+5}
\cdots 
\P^{(N)}_{2l+2N}Y^{(N)}_{2l+2N+1}.
 \endalign
$$
Hence we
inductively have
$$
\align
& \M_{l,l+1}{\P}_{l+1,l+N+1} \\
\overset{{(\kappa_1^{(N)})}^{-1} \kappa_0^{(N)}\cdots 
{(\kappa_1^{(1)})}^{-1} \kappa_0^{(1)}}\to{\simeq}
&
\P^{(1)}_{2l}Y^{(1)}_{2l+1}
\cdots
\P^{(N)}_{2l+2N-2}Y^{(N)}_{2l+2N-1}\M^{(N)}_{l+N,l+N+1}\\
=
& 
\P_{l,l+N}\M_{l+N,l+N+1}.
 \endalign
$$
By putting
$$
\kappa_\P = 
{(\kappa_1^{(N)})}^{-1} \kappa_0^{(N)}\cdots {(\kappa_1^{(1)})}^{-1} \kappa_0^{(1)}
:  \Sigma C^{(1)}C^{(2)}\cdots C^{(N)} \longrightarrow
 C^{(1)}C^{(2)}\cdots C^{(N)} \Sigma
$$
one has
$$
{\M}_{l,l+1} \P_{l+1,l+N+1}  \overset{\kappa_\P}\to{\simeq}
 \P_{l,l+N}{\M}_{l+N,l+N+1}.
$$
\qed
\enddemo
We set,
$$
\M^{[1]}_{l,l+1} =\M_{l,l+1},\qquad
{I}^{[1]}_{l,l+1} ={I}_{l,l+1}
$$
and
for $N\ge 2$
$$
\M^{[N]}_{l,l+1} 
 =\M_{Nl,Nl+1}{I}_{Nl+1,Nl+2}\cdots {I}_{N(l+1)-1,N(l+1)},\qquad
{I}^{[N]}_{l,l+1} 
= {I}_{Nl,Nl+N}.
$$
Then the pair
$({\M}^{[N]}, {I}^{[N]})$ is a symbolic matrix system over $\Sigma$.
Since one has 
$$
\align
\M^{[N]}_{l,l+m} 
& =\M^{[N]}_{l,l+1}{\M}^{[N]}_{l+1,l+2}\cdots {\M}^{[N]}_{(l+m)-1,(l+m)} \\
& ={\M}_{Nl,Nl+m}{I}_{Nl+m,N(l+m)},
\endalign
$$
a word $(a_1,\dots a_m) \in \Sigma^m$ is admissible for the subshift
$\Lambda_{\M^{[N]}}$ presented by $({\M}^{[N]}, {I}^{[N]})$
if and only if
it
is admissible for the subshift
$\Lambda_{\M}$ presented by $({\M}, {I}).$
Hence the subshifts
$\Lambda_{\M^{[N]}}$
and
$\Lambda_{\M}$
coincide.
\proclaim{Lemma 7.3}
Keep the above notations.
Put 
$$
\P^{[N]}_{l,l+1}  = \P_{Nl,Nl+N}, \qquad
\Q^{[N]}_{l,l+1}  = \Q_{Nl,Nl+N}, \qquad l\in \Zp.
$$
Then both $(\P^{[N]},{I}^{[N]}) $ and $(\Q^{[N]},{I}^{[N]}) $
are symbolic matrix systems such that 
\roster
\item"(i)"
the pair
$(\M^{[N]}, {I}^{[N]}) $ and
$(\P^{[N]},{I}^{[N]}) $,
and
the pair 
$(\M^{[N]}, {I}^{[N]}) $ and
$(\Q^{[N]},{I}^{[N]}) $
both form squares, and
\item"(ii)"
they satisfy the relations:
$$
\align
\M^{[N]}_{l,l+1}\P^{[N]}_{l+1,l+2}
& \overset{\kappa_\P}\to{\simeq}
 \P^{[N]}_{l,l+1}\M^{[N]}_{l+1,l+2}, \tag 7.4\\
\M^{[N]}_{l,l+1}\Q^{[N]}_{l+1,l+2}
& \overset{\kappa_\Q}\to{\simeq} \Q^{[N]}_{l,l+1}\M^{[N]}_{l+1,l+2},\qquad l \in \Zp. \tag 7.5 \\
\endalign
$$
\endroster
Hence the pair
$(\M^{[N]}, {I}^{[N]}) $ and
$(\P^{[N]},{I}^{[N]}) $,
and the pair
$(\M^{[N]}, {I}^{[N]}) $ and
$(\Q^{[N]},{I}^{[N]}) $
both give rise to LR textile $\lambda$-graph systems.
\endproclaim
\demo{Proof}
The assertion (i) is clear.
We will show the assertion (ii).
 By (7.2), 
one sees that
$$
\M_{Nl,Nl+1}\P_{Nl+1,N(l+1)+1}\overset{\kappa_\P}\to{\simeq}
{\P^{[N]}}_{l,l+1}\M_{N(l+1),N(l+1)+1}
$$
so that
$$
\M_{Nl,Nl+1}\P_{Nl+1,N(l+1)+1}{I}_{N(l+1)+1, N(l+2)}
\overset{\kappa_\P}\to{\simeq}
{\P^{[N]}}_{l,l+1}\M_{N(l+1),N(l+1)+1}{I}_{N(l+1)+1, N(l+2)}.
$$
Hence we get 
$$
\M^{[N]}_{l,l+1}\P^{[N]}_{l+1,l+2}
\overset{\kappa_\P}\to{\simeq}
 \P^{[N]}_{l,l+1}\M^{[N]}_{l+1,l+2}.
$$
We similarly have
$
\M^{[N]}_{l,l+1}\Q^{[N]}_{l+1,l+2}
\overset{\kappa_\Q}\to{\simeq} \Q^{[N]}_{l,l+1}\M^{[N]}_{l+1,l+2}.
$
\qed
\enddemo
\proclaim{Lemma 7.4}
The LR-textile $\lambda$-graph systems
 $
 {\Cal T}^{\P^{[N]}} = 
 {\Cal T}_{{\K}^{\M^{[N]}}_{\P^{[N]}}}
$ 
and 
${\Cal T}^{\Q^{[N]}}= 
 {\Cal T}_{{\K}^{\M^{[N]}}_{\Q^{[N]}}}
 $  
 defined by the relations (7.4) and (7.5) are both 1-1.
\endproclaim
 \demo{Proof}
 For ${(a_i)}_{i\in \Bbb Z} \in \Lambda_{\M^{[N]}},$
suppose that $a_l$ appears in a component of $\M^{[N]}_{l,l+1}.$
Since 
$\M^{[N]}_{l,l+1} = I_{Nl,N(l+1)-1}\M^{(0)}_{N(l+1)-1,N(l+1)}$,
the symbol $a_l$ 
 appears in a component of $\M^{(0)}_{N(l+1)-1,N(l+1)}.$
By the specified equivalence
$\M^{(0)}_{N(l+1)-1,N(l+1)}
\overset{\kappa_0^{(1)}}\to{\simeq}
\P_{2\{N(l+1)-1\}}^{(1)}\Q_{2\{N(l+1)-1\}+1}^{(1)}$,
a symbol
$c_{2\{N(l+1)-1\}}^{(1)}d_{2\{N(l+1)-1\}+1}^{(1)}
=\kappa_0^{(1)}(a_{l})\in C D$ 
appears in a component of
$\P_{2\{N(l+1)-1\}}^{(1)}\Q_{2\{N(l+1)-1\}+1}^{(1)}$.
 For $a_{l+1}$, 
 the corresponding symbol
$c_{2\{N(l+2)-1\}}^{(1)}$ 
 appears in a component of
$\P_{2\{N(l+2)-1\}}^{(1)}$.
Since one has
$$
I^{(0)}_{N(l+1),N(l+2)-1}\P_{2\{N(l+2)-1\}}^{(1)}
=
\P_{2N(l+1)}^{(1)}X_{2N(l+1)+1}^{(1)}\cdots 
X_{2\{N(l+2)-2\} +1}^{(1)}
X_{2\{N(l+2)-1\} }^{(1)} 
 $$
the corresponding symbol to
$c_{2\{N(l+2)-1\}}^{(1)}$ 
 appears in a component of
$\P_{2N(l+1)}^{(1)}$, that is denoted by
$c_{2N(l+1)}^{(1)}$. 
As one sees that
$$
\align
& \M^{[N]}_{l,l+1}{\P}^{[N]}_{l+1,l+2} \\
 =& I_{Nl,N(l+1)-1}{\M}^{(0)}_{N(l+1)-1,N(l+1)} 
 \P_{N(l+1),N(l+2)}\\
 \overset{\kappa_0^{(1)}}\to {\simeq}
& I_{Nl,N(l+1)-1}\P_{2\{N(l+1)-1\}}^{(1)}\Q_{2\{N(l+1)-1\}+1}^{(1)}
  \P^{(1)}_{2N(l+1)}Y^{(1)}_{2N(l+1)+1}\\
& \cdot
\P^{(2)}_{2N(l+1)+2}Y^{(2)}_{2N(l+1)+3}
\cdots 
\P^{(N)}_{2N(l+1)+2N-2}Y^{(N)}_{2N(l+1)+2N-1}\\
\endalign
$$
the symbol
$d_{2N(l+1)-1}^{(1)}c_{2N(l+1)}^{(1)}$
appears in a component of
$\Q_{2N(l+1)-1}^{(1)}\P_{2N(l+1)}^{(1)}.$
As
$$
\align
& \Q^{(1)}_{2N(l+1)-1}\P^{(1)}_{2N(l+1)} Y^{(1)}_{2N(l+1) +1}\\
=
& Y^{(1)}_{2N(l+1)-1}\Q^{(1)}_{2N(l+1)} \P^{(1)}_{2N(l+1) +1}
\overset{{\kappa_1^{(1)}}^{-1}}\to {\simeq} 
Y^{(1)}_{2N(l+1)-1}\M^{(1)}_{N(l+1),N(l+1)+1}
\endalign
$$
$d_{2N(l+1)-1}^{(1)}c_{2N(l+1)}^{(1)}$
appears in a component of
$\Q_{2N(l+1)}^{(1)}\P_{2N(l+1)+1}^{(1)}$,
that is written as
$d_{2N(l+1)}^{(1)}c_{2N(l+1)+1}^{(1)}.$
Hence
${\kappa_1^{(1)}}^{-1}(d_{2N(l+1)}^{(1)}c_{2N(l+1)+1}^{(1)})$
appears in a component of
$\M^{(1)}_{N(l+1),N(l+1)+1}$.
This procedure shows that 
for a given
 ${(a_i)}_{i\in \Bbb Z} \in \Lambda_{\M^{[N]}},$
by starting from $a_l $ in a component of $\M^{(0)}_{N(l+1)-1,N(l+1)}$
a symbol
$c^{(1)}_{2N(l+1)}$ in a component of $\P^{(1)}_{2N(l+1)}$ 
is determined and also 
${\kappa_1^{(1)}}^{-1}(d_{2N(l+1)}^{(1)}c_{2N(l+1)+1}^{(1)})$
in a component of
$\M^{(1)}_{N(l+1),N(l+1)+1}$
is determined.
One may  next find a symbol in a component of
$\P^{(2)}_{2N(l+1)+2}$ 
and a symbol in 
a component of
$\M^{(2)}_{N(l+1)+1,N(l+1)+2}.$
One inductively finds corresponding symbols in  
$\P^{(1)}_{2N(l+1)},\P^{(2)}_{2N(l+1)+2},\dots,\P^{(N)}_{2N(l+1)+2N-3}
$.
Hence ons finds a symbol in 
$\P^{(N)}_{l+1,l+2}$.
That is, 
a given sequence
 ${(a_i)}_{i\in \Bbb Z} \in \Lambda_{\M^{[N]}}$
determines a symbol
in $\P^{(N)}_{l+1,l+2}, l\in \Zp$
so that through the relation (7.4)
the labeled squares in the LR textile $\lambda$-graph system
are determined.
Hence we conclude that $\xi$ is injective, 
and similarly see that $\eta$ is injective.
 \qed
 \enddemo
 As stated in the begining of this section,
 an automorphism of a subshift $\Lambda$ 
presented by $(\M,I)$  
is given by a properly strong shift equivalence 
$$
(\M,I) \underset{1-pr}\to {\approx} (\M^{(1)},I^{(1)})
       \underset{1-pr}\to {\approx} 
       \cdots 
       \underset{1-pr}\to {\approx} 
(\M^{(N-1)},I^{(N-1)}) \underset{1-pr}\to {\approx} 
( \M, I)
$$
in $N$-step for some $N$, 
and conversely   
a properly strong shift equivalence 
from $(\M,I)$ to itself
gives rise to an automorphism of the subshift.
Put for $k=1,\dots,N$
$$
\align
\Lambda_{C^{(k)}D^{(k)}} & = \{ {(c_i d_i )}_{i \in \Bbb Z} \mid
c_i \in C^{(k)}, d_i \in D^{(k)}, i\in \Bbb Z \},\\
\Lambda_{D^{(k)}C^{(k)}} & = \{ {(d_i c_i )}_{i \in \Bbb Z} \mid
c_i \in C^{(k)}, d_i \in D^{(k)}, i\in \Bbb Z \}.
\endalign
$$
Define
$
\zeta_+^{(k)} : \Lambda_{C^{(k)}D^{(k)}} \longrightarrow \Lambda_{D^{(k)}C^{(k)}}
$
and
$
 \zeta_-^{(k)} : \Lambda_{C^{(k)}D^{(k)}} \longrightarrow \Lambda_{D^{(k)}C^{(k)}}$
 by setting
$
\zeta_+^{(k)}({(c_id_i)}_{i\in \Bbb Z}) = {(d_i c_{i+1})}_{i\in \Bbb Z}
$ 
and
$
\zeta_-^{(k)}({(c_id_i)}_{i\in \Bbb Z}) = {(d_{i-1} c_i)}_{i\in \Bbb Z}
$
respectively.
Then 
$\zeta_+^{(k)}$ is called a forward bipartite conjugacy and 
$\zeta_-^{(k)}$ is called a backward bipartite conjugacy ([N2]).
Nasu's result ([N], [N2])
says that any automorphism $\phi$ is factorized as follows:
$$
\phi = ({(\kappa_1^{(N)})}^{-1}\circ \zeta_{\pm}^{(N)}\circ \kappa_0^{(N)})
\circ \cdots \circ
({(\kappa_1^{(1)})}^{-1}\circ \zeta_{\pm}^{(1)}\circ \kappa_0^{(1)})
$$
where
$\zeta_{\pm}^{(N)}, \dots, \zeta_{\pm}^{(1)}$ are forward or backward bipartite conjugacies.
Since properly strong shift equivalence  corresponds exactly to bipartite codes of Nasu, the above  factorization of $\phi$ is so called Nasu's $\kappa-\zeta$ factorization ([N], [N2]). 
Following Nasu, an automorphism $\phi$ is 
said to be forward if 
$\zeta_{\pm}^{(N)}, \dots, \zeta_{\pm}^{(1)}$ 
are all forward bipartite conjugacies 
$\zeta_{+}^{(N)}, \dots, \zeta_{+}^{(1)}$.
 
 \proclaim{Lemma 7.5}
 Let $(\Lambda,\sigma)$ be a subshift presented by 
 $(\M,I)$.
Let $\phi$ be a forward automorphism on $(\Lambda,\sigma)$
defined by a properly strong shift equivalence 
$$
(\M,I) \underset{1-pr}\to {\approx} (\M^{(1)},I^{(1)})
       \underset{1-pr}\to {\approx} 
       \cdots 
       \underset{1-pr}\to {\approx} 
(\M^{(N-1)},I^{(N-1)}) \underset{1-pr}\to {\approx} 
( \M, I)
$$
in $N$-step.
Let ${\Cal T}^{\P^{[N]}}$ and ${\Cal T}^{\Q^{[N]}}$ 
be the LR textile $\lambda$-graph systems defined by the relations (7.4) and 
(7.5) respectively.
Then we have
$$
\varphi_{{\Cal T}^{\P^{[N]}}} =\phi,
\qquad 
\varphi_{{\Cal T}^{\Q^{[N]}}} =\phi^{-1}\circ \sigma
$$
as automorphisms on $\Lambda = \Lambda_{\M}$
under the identification $\Lambda_{\M^{[N]}} = \Lambda_\M$.
\endproclaim
\demo{Proof}
Keep the notations as in the proof of the previous lemma.
For 
${(a_i)}_{i\in \Bbb Z}\in \Lambda_{\M^{[N]}}$,
by putting
$$
c_{2\{N(l+1)-1\}}^{(1)}d_{2\{N(l+1)-1\}+1}^{(1)}
=\kappa_0^{(1)}(a_{l})\in C D,\qquad l\in \Zp
$$
the symbol
$c_{2N(l+1)-2}^{(1)}d_{2N(l+1)-1}^{(1)}$
 is written as
$d_{2N(l+1)}^{(1)}c_{2N(l+1)+1}^{(1)}$
and
${(\kappa_1^{(1)})}^{-1}(d_{2N(l+1)}^{(1)}c_{2N(l+1)+1}^{(1)})$
defines a symbol
of a component of $\M_{N(l+1),N(l+1)+1}^{(1)}.$
This procedure is nothing but to apply the map
${(\kappa_1^{(1)})}^{-1}\circ \zeta_+^{(1)} \circ \kappa_0^{(1)}$.
We next do this procedure to the symbol
${(\kappa_1^{(1)})}^{-1}(d_{2N(l+1)}^{(1)}c_{2N(l+1)+1}^{(1)})$
that corresponds to apply 
the map
$
{(\kappa_1^{(2)})}^{-1}\circ \zeta_+^{(2)} \circ \kappa_0^{(2)}
$
and get the symbols
$
({(\kappa_1^{(2)})}^{-1}\circ \zeta_+^{(2)} \circ \kappa_0^{(2)})
\circ ({(\kappa_1^{(1)})}^{-1}\circ \zeta_+^{(1)} \circ \kappa_0^{(1)})
((a_l)).
$
We continue this procedures and finally get the element 
$
({(\kappa_1^{(N)})}^{-1}\circ \zeta_+^{(N)} \circ \kappa_0^{(N)})\cdots
({(\kappa_1^{(2)})}^{-1}\circ \zeta_+^{(2)} \circ \kappa_0^{(2)})
\circ ({(\kappa_1^{(1)})}^{-1}\circ \zeta_+^{(1)} \circ \kappa_0^{(1)})
((a_l))
$
in $\M^{(N)}_{N(l+1)+N-1,N(l+1)+N},l\in \Zp.
$
The elements lie in the bottoms of the squares arising from the relation (7.4),
and hence that are the element 
$\eta\circ \xi^{-1}((a_l))$.
Hence we have  
$
\varphi_{{\Cal T}^{\P^{[N]}}}
=({(\kappa_1^{(N)})}^{-1}\circ \zeta_+^{(N)} \circ \kappa_0^{(N)})\cdots
({(\kappa_1^{(2)})}^{-1}\circ \zeta_+^{(2)} \circ \kappa_0^{(2)})
\circ ({(\kappa_1^{(1)})}^{-1}\circ \zeta_+^{(1)} \circ \kappa_0^{(1)}).
$
\qed
\enddemo

We assume that the previously defined metric is equipped with $\Lambda$.
Then the homeomorphism $\sigma$ has $1$ as its expansive constant.
Therefore we have
\proclaim{Theorem 7.6}
 Let $(\Lambda,\sigma)$ be a subshift presented by a symbolic matrix system
 $(\M,I)$.
Let $\phi$ be a forward automorphism on $(\Lambda,\sigma)$
defined by a properly strong shift equivalence 
$$
(\M,I) \underset{1-pr}\to {\approx} (\M^{(1)},I^{(1)})
       \underset{1-pr}\to {\approx} 
       \cdots 
       \underset{1-pr}\to {\approx} 
(\M^{(N-1)},I^{(N-1)}) \underset{1-pr}\to {\approx} 
( \M, I)
$$
in $N$-step.
If $\phi$ is expansive with $\frac{1}{k}$ as its expansive constant for some $k \in \Bbb N$,
 the dynamical system
$(\Lambda,\phi)$
is topologically conjugate to the subshift 
$(\Lambda_{  {\P^{[N]}}_{\Cal T}^{[2k]}  },
\sigma_{{\P^{[N]}}_{\Cal T}^{[2k]}})$
presented by the
$2k$-heigher block
$({\P^{[N]}}_{\Cal T}^{[2k]},I^{{\P^{[N]}}_{\Cal T}^{[2k]}})$
 of the symbolic matrix system
$(\P^{[N]},I^{\P^{[N]}})$ relative to the LR textile $\lambda$-graph system 
$ {\Cal T}_{{\K}^{\M^{[N]}}_{\P^{[N]}}}$
defined by the specification
$$
\M_{l,l+1}^{[N]}\P_{l+1,l+2}^{[N]} \overset{\kappa^{\P^{[N]}}}\to {\simeq} 
\P_{l,l+1}^{[N]}\M_{l+1,l+2}^{[N]},\qquad l \in \Zp,
$$
where
$$
\align
\P^{[N]}_{l,l+1} 
& =\P^{(1)}_{2Nl}Y^{(1)}_{2Nl+1}\P^{(2)}_{2Nl+2}Y^{(2)}_{2Nl+3}
\cdots \P^{(N)}_{2Nl+2N-2}Y^{(N)}_{2Nl+2N-1},\\
{I}^{[N]}_{l,l+1} 
& ={I}_{Nl,Nl+1}{I}_{Nl+1,Nl+2}\cdots {I}_{N(l+1)-1,N(l+1)},\qquad l \in \Zp
\endalign
$$
and $\P^{(i)}_{2Nl+2(i-1)},Y^{(i)}_{2Nl+2i-1}, i=1,\dots,N$ are matrices 
appearing in the properly strong shift equivalence in (7.1).
\endproclaim
\demo{Proof}
Since the LR textile $\lambda$-graph system ${\Cal T}^{\P^{[N]}}$
is nondegenerate, 1-1, and surjective, by Lemma 7.5 the assertion is proved in a similar way to the proof of Theorem 6.4.
\qed
\enddemo
For
 $k\ge 0, n\ge 1,$
let
$(\P^k \M^n,I^{kN+n})$
be the symbolic matrix system 
defined by setting
$$
\align
& {({\P}^k{\M}^n)}_{l,l+1}\\ 
=& 
{\P}_{l(kN+n),l(kN+n)+N}{\P}_{l(kN+n) +N,l(kN+n)+2N}\cdots 
{\P}_{l(kN+n)+(k-1)N,l(kN+n)+kN}\cdot\\
& \cdot {\M}_{l(kN+n)+kN,l(kN+n)+kN+1}{\M}_{l(kN+n)+kN+1,l(kN+n)+kN+2}\cdots
{\M}_{(l+1)(kN+n)-1,(l+1)(kN+n)},\\
& I^{kN+n}_{l,l+1} \\ 
=& 
I_{l(kN+n),l(kN+n)+1}I_{l(kN+n)+1,l(kN+n)+2}\cdots 
I_{(l+1)l(kN+n)-1,(l+1)(kN+n)},\qquad l \in \Zp.
\endalign
$$
\proclaim{Lemma 7.7}
The subshift 
$(\Lambda_{{\P^{[N]}}^k {\M^{[N]}}^n},\sigma_{{\P^{[N]}}^k {\M^{[N]}}^n})$
presented by the symbolic matrix system
$({{\P^{[N]}}^k {\M^{[N]}}^n},I^{ {[N]}^{kN+n}})$
coincides with the subshift
$(\Lambda_{\P^k \M^n},\sigma_{\P^k \M^n})$
presented by the symbolic matrix system
$(\P^k \M^n,I^{kN+n}).$
\endproclaim
\demo{Proof}
It is easy to see that 
the set of all admissible words of the subshift 
$(\Lambda_{{\P^{[N]}}^k {\M^{[N]}}^n},\sigma_{{\P^{[N]}}^k {\M^{[N]}}^n})$
coincides with the set of all admissible words of the subshift
$(\Lambda_{\P^k \M^n},\sigma_{\P^k \M^n})$.
\qed
\enddemo

We reach our main theorem. 
 \proclaim{Theorem 7.8}
 Let $(\Lambda,\sigma)$ be a subshift presented by a symbolic matrix system
 $(\M,I)$.
Let $\phi$ be a forward automorphism on $(\Lambda,\sigma)$
defined by a properly strong shift equivalence 
$$
(\M,I) \underset{1-pr}\to {\approx} (\M^{(1)},I^{(1)})
       \underset{1-pr}\to {\approx} 
       \cdots 
       \underset{1-pr}\to {\approx} 
(\M^{(N-1)},I^{(N-1)}) \underset{1-pr}\to {\approx} 
( \M, I)
$$
in $N$-step.
Then the dynamical system
$(\Lambda,\phi^k\sigma^n)$
is topologically conjugate to the subshift 
$(\Lambda_{\P^k \M^n},\sigma_{\P^k \M^n})$
presented by the symbolic matrix system
$(\P^k \M^n,I^{kN+n})$
 for
 $k\ge 0, n\ge 1,$
defined by
$$
\align
& {({\P}^k{\M}^n)}_{l,l+1}\\ 
=& 
{\P}_{l(kN+n),l(kN+n)+N}{\P}_{l(kN+n) +N,l(kN+n)+2N}\cdots 
{\P}_{l(kN+n)+(k-1)N,l(kN+n)+kN}\cdot\\
& \cdot {\M}_{l(kN+n)+kN,l(kN+n)+kN+1}{\M}_{l(kN+n)+kN+1,l(kN+n)+kN+2}\cdots
{\M}_{(l+1)(kN+n)-1,(l+1)(kN+n)},\\
& I^{kN+n}_{l,l+1} \\ 
=& 
I_{l(kN+n),l(kN+n)+1}I_{l(kN+n)+1,l(kN+n)+2}\cdots 
I_{(l+1)l(kN+n)-1,(l+1)(kN+n)},\qquad l \in \Zp
\endalign
$$
where
$
\P_{l,l+1} 
=\P^{(1)}_{2l}Y^{(1)}_{2l+1}\P^{(2)}_{2l+2}Y^{(2)}_{2l+3}
\cdots \P^{(N)}_{2l+2N-2}Y^{(N)}_{2l+2N-1}
$
and $\P^{(i)}_{2l+2(i-1)},Y^{(i)}_{2l+2i-1}, i=1,\dots,N$ are matrices 
appearing in the properly strong shift equivalence in (7.1).

And also 
$(\Lambda,{(\sigma\phi^{-1})}^k\sigma^n)$
is topologically conjugate to the subshift 
$(\Lambda_{\Q^k \M^n},\sigma_{\Q^k \M^n})$
presented by the similarly defined symbolic matrix system
$(\Q^k \M^n, I^{kN +n})$
 for
 $k\ge 0, n\ge 1.$
\endproclaim
\demo{Proof}
By a similar discussion to the proof of Theorem 6.5,
the dynamical system
$(\Lambda,\phi^k\sigma^n)$
is topologically conjugate to the subshift 
$(\Lambda_{{\P^{[N]}}^k {\M^{[N]}}^n},\sigma_{{\P^{[N]}}^k {\M^{[N]}}^n})$,
that is 
$(\Lambda_{\P^k \M^n},\sigma_{\P^k \M^n})$
by Lemma 7.7.
The assertion for 
the dynamical system
$(\Lambda,{(\sigma\phi^{-1})}^k\sigma^n)$
is similarly shown.
\qed
\enddemo

\heading 8. An application 
\endheading
Let
$\phi$ be an automorphism of a subshift $\Lambda$ over $\Sigma$.
We say that $\phi$
is given by a specification $\pi$ of a symbolic matrix system $(\M,I)$
if $(\M,I)$ presents the subshift $\Lambda$ and 
there exists a specification 
$\pi:\Sigma \rightarrow \Sigma$
such that $\pi$ gives rise to a specified equivalence 
$$
\M_{l,l+1}
\overset{\pi}\to {\simeq}\M_{l,l+1}
\qquad
\text{ for } l\in \Zp,
$$
and $\phi$ is given by the symbolic automorphism of $\Lambda$ induced by $\pi$.
The automorphism $\phi$ is written as $\phi_{\pi}$.
We note that the induced automorphism of the $\lambda$-graph system $\frak L$
for $(\M,I)$ by the specification $\pi$ fixes the vertices of $\frak L$.

\noindent
{\bf Definition.}
An automorphism $\phi$ of a subshift $\Lambda$
is called a {\it simple}\ automorphism 
if there exist an
automorphism
$\phi_{\pi}$ of a subshift $\Lambda_{\Cal M}$   
that
is given by a specification $\pi$ of 
a symbolic matrix system $(\M,I)$,
and a topological conjugacy 
$\psi: \Lambda \rightarrow \Lambda_{\Cal M}$
such that
$$
\phi = \psi^{-1}\circ \phi_{\pi} \circ \psi.
$$
The notion of a simple automorphism of a sofic shift has been introduced 
by M. Nasu in [N2].

As an application of our result, we see the following proposition.
\proclaim{Proposition 8.1}
Let $(\Lambda,\sigma)$ be a subshift over $\Sigma$.
If an automorphism $\phi_{\pi}$ of $(\Lambda,\sigma)$ 
is given by
a specification $\pi$ of a symbolic matrix system $(\M,I)$,
 the topological dynamical system 
$(\Lambda,\phi_\pi \circ \sigma^n )$ 
is topologically conjugate to the $n$-th power 
$(\Lambda, \sigma^n )$
of 
$(\Lambda,\sigma)$
for $n \in \Bbb Z, n\ne 0.$
\endproclaim
\demo{Proof}
By assumption,
the automorphism
$\phi_\pi$ is given by the one-block code
$
\phi_\pi( {(x_i)}_{i \in \Bbb Z}) = 
       { (\pi(x_i))}_{i \in \Bbb Z}
$ for
${(x_i)}_{i \in \Bbb Z} \in \Lambda$.
We will realize $\phi_\pi$ to be a forward automorphism defined by a properly strong shift equivalence from $(\M,I)$ to itself.
Let $\Bbb I$ be an arbitrary fixed symbol.
Put the alphabets 
$$
  C = \{{\Bbb I}\}, \qquad  D = \Sigma.
$$
Define the specifications $\kappa$ from $\Sigma$ to $C\cdot D$ 
and $\kappa'$ from $\Sigma$ to $D\cdot C$
by setting
$$
\kappa(\gamma) = {\Bbb I}\cdot \pi(\gamma),\qquad
\kappa'(\gamma) = \gamma \cdot {\Bbb I},\qquad \gamma \in \Sigma.
$$   
Suppose that the both matrices $\M_{l,l+1}, I_{l,l+1}$
are $m(l) \times m(l+1)$ matrices.
Let                                               
$I_l({\Bbb I})$ and $I_l(1)$ be  the $m(l) \times m(l)$
diagonal matrices with diagonal entries ${\Bbb I}$ and $1$ 
respectively.
Put  
$n(2l) = n(2l-1) = m(l)$ for $l\in \Bbb N$, and
$n'(2l) = n'(2l+1) = m(l)$ for $l\in \Zp$.
Define matrices 
$\P_l, \Q_l, X_l,Y_l$ for $l\in \Zp$ 
by setting
$$
\P_{2l} = I_l({\Bbb I}),\qquad
\P_{2l+1} = I_{l+1}({\Bbb I}),\qquad
\Q_{2l} = \Q_{2l+1} = \M_{l,l+1}
$$
and
$$
X_{2l} = Y_{2l+1} = I_{l,l+1},\qquad
X_{2l+1} = I_{l+1}(1),\qquad
Y_{2l} = I_l(1)
$$
By noticing that the matrices  
$X_{2l+1}, Y_{2l}$ are identity matrices, 
the above matrices give rise to 
a properly strong shift equivalence
$$
(\P,\Q,X,Y): (\M,I)\underset{1-pr}\to {\approx} (\M,I)
$$
in $1$-step
from
$(\M,I)$ to itself.
It is then direct to see that the automorphism 
$\phi_\pi$ is the forward automorphism of 
the above properly strong shift equivalence.
Put
$$
\P_{l,l+1} = \P_{2l}Y_{2l+1} = I_l({\Bbb I}) I_{l,l+1},\qquad l \in \Zp.
$$  
For $ n \in \Bbb Z$ with $n > 0$,
we set 
$$
\align
(\P \M^n )_{l,l+1} 
= & \P_{(n+1)l,(n+1)l+1}\cdots \P_{(n+1)l+n-1,(n+1)l+n} \\
  & \cdot \M_{(n+1)l+n,(n+1)l+n+1}\cdots \M_{(n+1)(l+1)-1,(n+1)(l+1)},\\
I^{n+1}_{l,l+1}
 = & I_{(n+1)l,(n+1)l+1}\cdots I_{(n+1)(l+1)-1,(n+1)(l+1)}.
\endalign
$$
 Then by Theorem 6.5,
 the topological dynamical system
 $(\Lambda, \phi_\pi \circ \sigma^n)$
 is realized as the subshift
$(\Lambda_{\P \M^n}, \sigma_{\P \M^n})$
presented by the symbolic matrix system
$(\P \M^n, I^{n+1})$.
Since the symbolic matrix system
$(\P \M^n, I^{n+1})$
does not depend on the choice of the specification $\pi$ on $\Sigma$,
we have
$(\Lambda, \phi_\pi \circ \sigma^n)$
is topologically conjugate to
$(\Lambda, \phi_\id \circ \sigma^n)$
 where $\phi_\id$ is the automorphism coming from the identity permutation.
 Hence 
$(\Lambda, \phi_\pi \circ \sigma^n)$
is topologically conjugate to the $n$-th power
$(\Lambda, \sigma^n)$ 
of 
$(\Lambda, \sigma)$.
For $n \in \Bbb Z$ with $n<0$,
the above argument says that
the dynamical system
$(\Lambda, \phi_{\pi^{-1}} \circ \sigma^{-n})$
is topologically conjugate to
$(\Lambda,  \sigma^{-n})$.
This implies that 
$(\Lambda, \phi_\pi \circ \sigma^n)$
is topologically conjugate to
$(\Lambda,  \sigma^n)$.
\qed
\enddemo
Thanks to this proposition, one has the following theorem.
\proclaim{Theorem 8.2}
If an automorphism $\phi$ of a subshift 
$(\Lambda, \sigma)$ is a simple automorphism, 
the dynamical system
$(\Lambda, \phi \circ \sigma^n)$ 
is topologically conjugate to the $n$-th power 
$
(\Lambda, \sigma^n)
$
of the subshift  
$
(\Lambda, \sigma)
$
for $n \in \Bbb Z, n\ne 0$.
\endproclaim
\demo{Proof}
As $\phi$ is a simple automorphism of $\Lambda$, 
there exist an
automorphism
$\phi_{\pi}$ of a subshift $\Lambda_{\Cal M}$   
that
is given by a specification $\pi$ of 
a symbolic matrix system
$(\M,I)$,
and a topological conjugacy 
$\psi: \Lambda \rightarrow \Lambda_{\Cal M}$
such that
$$
\phi = \psi^{-1}\circ \phi_{\pi} \circ \psi.
$$
Hence 
$$ 
\psi \circ ( \phi \circ \sigma^n )\circ \psi^{-1}
= \phi_{\pi} \circ {\sigma_{\M}}^{n}.
$$
By Proposition 8.1,
$(\Lambda, \phi_{\pi} \circ \sigma^n)$ is topologically conjugate to
$(\Lambda_{\M},{\sigma_{\M}}^n)$.
Hence
$(\Lambda,\phi \circ \sigma^n)$ is  
topologically conjugate to
$(\Lambda,\sigma^n)$.
\qed
\enddemo

\medskip

We finally give an example.
Let us  consider the full shift $\Sigma^{\Bbb Z}$ with alphabet
$\Sigma$.
It is easy to see that any permutation $\pi$ on $\Sigma$ yields 
a simple automorphism
$\phi_{\pi}$
of $\Sigma^{\Bbb Z}$.
Hence we have
\proclaim{Corollary 8.3}
Let $\phi_{\pi}$ 
be the automorphism of the full shift 
$\Sigma^{\Bbb Z}$ defined by a permutation 
$\pi $ on the symbols $\Sigma$.
Then the topological dynamical system
 $(\Sigma^{\Bbb Z}, \phi_{\pi} \circ \sigma)$
 is realized as the original full shift
 $(\Sigma^{\Bbb Z},\sigma)$.
 That is, 
 $\phi_{\pi} \circ \sigma$ 
 is topologically conjugate to the original full shift
 $\sigma$.
 \endproclaim
    \medskip

{\it Acknowledgments.}\/
The author would like to thank Masakazu Nasu for his useful comments on the first draft of this paper, suggestions in Section 8  and series of lectures on textile systems at Yokohama City University, December 2001. 
He in particular kindly pointed out some errors in the first draft.
The author also would like to thank to Wolfgang Krieger
for his discussions and constant encouragements.

\bigskip

{\it e-mail} :kengo{\@}yokohama-cu.ac.jp


\Refs

\refstyle{A}
\widestnumber\key{DGSW}





\ref\key BK
\by M. Boyle and W. Krieger
\paper Periodic points and automorphisms of the shift   
\jour Trans. Amer. Math. Soc.
\vol 302 
\yr 1988
\pages 125--148
\endref


\ref\key BLR
\by M. Boyle, D. Lind and D. Rudolph
\paper The automorphism group of a subshift of finite type   
\jour Trans. Amer. Math. Soc.
\vol 306 
\yr 1988
\pages 71--114
\endref


\ref
\key Bra
\by O. Bratteli 
\paper Inductive limits of finite-dimensional $C^*$-algebras
\jour Trans. Amer. Math. Soc.
\vol 171
\yr 1972
\pages 195--234
\endref






\ref\key CK
\by J. Cuntz and W. Krieger
\paper A class of $C^*$-algebras and topological Markov chains
\jour Invent. Math.
\vol 56
\yr 1980
\pages 251--268
\endref



\ref\key Fi
\by R. Fischer
\paper Sofic shifts and graphs                     
\jour Monats. Math.                 
\vol 80
\yr 1975
\pages 179-186
\endref


\ref\key H  
\by G. A. Hedlund   
\paper Endomorphisms and automorphisms of the shift dynamical system
\jour Math. Systems Theory          
\vol 3 
\yr 1969
\pages 320-375
\endref


\ref\key HN
\by T. Hamachi and M. Nasu
\paper Topological conjugacy for $1$-block factor maps of subshifts and sofic  covers
\jour Proceedings of Maryland Special Year in Dynamics 1986-87, 
Springer -Verlag Lecture Notes in Math
\vol 1342
\yr 1988 
\pages 251--260
\endref 

\ref\key KRW
\by K. H. Kim, F. W. Roush and J. W. Wagoner
\paper Automorphisms of the dimension group and gyration numbers
\jour  J. Amer. Math. Soc. 
\vol 5  
\yr 1992
\pages 191--212
\endref




\ref\key Ki
\by B. P. Kitchens
\book Symbolic dynamics
\publ Springer-Verlag
\publaddr Berlin, Heidelberg and New York
\yr 1998
\endref




\ref\key Kr
\by W. Krieger
\paper On sofic systems I
\jour Israel J. Math.
\vol 48
\yr 1984
\pages 305--330
\endref




\ref\key Kr2
\by W. Krieger
\paper On sofic systems II
\jour  Israel J. Math
\vol 60
\yr 1987
\pages 167--176
\endref






\ref\key KMT
\by W. Krieger, B. Marcus and S. Tuncel
\paper On automorphisms of Markov chains 
\jour  Trans. Amer. Math. Soc. 
\vol 333
\yr 1992
\pages 531--565
\endref


\ref\key KM 
\by W. Krieger and K. Matsumoto
\paper A lambda-graph system for the Dyck shift and its K-groups
\jour Doc. Math.
\vol 8
\yr 2003
\pages 79-96
\endref




\ref\key LM
\by D. Lind and B. Marcus
\book An introduction to symbolic dynamics and coding
\publ Cambridge University Press.
\yr 1995
\endref



\ref\key Ma
\by K. Matsumoto
\paper On $C^*$-algebras associated with subshifts
\jour Internat. J. Math.
\vol 8
\yr 1997
\pages 357--374
\endref


\ref\key Ma2
\by K. Matsumoto
\paper Presentations of subshifts and their topological conjugacy invariants
\jour Doc. Math.
\vol 4
\yr 1999
\pages 285-340
\endref


\ref\key Ma3
\by K. Matsumoto
\paper  $C^*$-algebras associated with presentations of subshifts 
\jour Doc. Math.
\vol 7
\yr 2002
\pages 1-30
\endref

\ref\key Ma4
\by K. Matsumoto
\paper  On strong shift equivalence of symbolic matrix systems  
\jour   Ergodic Theory Dynam. Systems
\vol 23
\yr 2003
\pages 1551-1774
\endref


\ref\key N
\by M. Nasu
\paper Topological conjugacy for sofic shifts 
\jour Ergodic Theory Dynam. Systems
\vol 6
\yr 1986
\pages 265--280
\endref

\ref\key N2
\by M. Nasu
\paper Textile systems for endomorphisms and automorphisms of the shift
\jour Mem. Amer. Math. Soc. 
\vol 546
\yr 1995
\endref

\ref\key N3
\by  M. Nasu
\paper Topological conjugacy for sofic systems and extensions of automorphisms of finite subsystemss of topological Markov shifts
\jour Proceedings of Maryland Special Year in Dynamics 1986-87, 
Springer -Verlag Lecture Notes in Math
\vol 1342
\yr 1988 
\pages 564--607
\endref 


\ref\key Wa
\by J. B. Wagoner 
\paper Triangle identities and symmetries of a shift of finite type
\jour Pacific J. Math. 
\vol 144
\yr 1990
\pages 181--205
\endref



\ref\key We
\by B. Weiss
\paper Subshifts of finite type and sofic systems 
\jour Monats. Math.
\vol 77
\yr 1973
\pages 462--474
\endref

\ref\key Wi
\by R. F. Williams
\paper Classification of subshifts of finite type 
\jour Ann. Math.
\vol 98
\yr 1973
\pages 120--153
\finalinfo erratum, Ann. Math.
$ 99(1974), 380-381$
\endref
\bye